\documentclass[12pt]{article}
\usepackage{amssymb}
\input epsf

\def\bfb{{\bf b}}
\def\D{{\cal D}}
\def\gplus{ G^A_\oplus}

\def\lapplus{\lap^A_\oplus}

\def\normaax{(|a_X|+|\a_X|)}
\def\normaay{(|a_Y|+|\a_Y|)}

\def\sqll{|\log\l|^{1/2}}
\def\tx{\tilde{X}}
\def\tX{\tilde{X}}
\def\ty{\tilde{Y}}

\def\tw{\tilde{W}}

\def\bdotb{ |\bfb_X||\bfb_Y| }

\def\bdota{ |\bfb_X|\normaay }

\def\adotb{ |\bfb_Y|\normaax }

\def\adota{ \normaax\normaay }

\def\baab{ \left( \bdota+\adotb \right) }

\def\bwedgeb{ |\bfb_X\wedge\bfb_Y| }

\def\awedgea{\left( |a_X||\a_Y| + |\a_X||a_Y| + |\a_X\wedge\a_Y| \right)}

\def\D{{\cal D}}
\def\muloc{\mu_{\rm loc}}

\def\F{{\cal F}}
\def\R{{\cal R}}

\def\bfr{{\bf R}}
\def\bfz{{\bf Z}}

\def\ddxi{\frac{\partial}{\partial x^i}}
\def\ddxj{\frac{\partial}{\partial x^j}}

\def\dist{{\rm dist}}

\def\embed{\hookrightarrow}

\def\ident{\equiv}

\def\intersect{\bigcap}
\def\iso{\cong}
\def\lap{\Delta}
\def\lapa{\Delta^A}
\def\lessim{\lesssim}
\def\lb{\langle}
\def\rb{\rangle}
\def\lnorm{\left\|}
\def\rnorm{\right\|}
\def\na{\nabla} 
\def\Na{\nabla^A}

\def\reals{ {\bf R} }

\def\plus{\oplus}
\def\supp{{\rm supp}}
\def\tensor{\otimes}
\def\tr{{\rm tr}}

\def\vol{{\rm Vol}}

\def\be{\begin{equation}}
\def\ee{\end{equation}}
\def\bearray{\begin{eqnarray}}
\def\eearray{\end{eqnarray}}
\def\bestar{\begin{eqnarray*}}
\def\eestar{\end{eqnarray*}}
\def\ben{\begin{displaymath}}
\def\een{\end{displaymath}}

\def\non{\nonumber}

\def\a{\alpha}
\def\b{\beta}
\def\d{\delta}
\def\e{\epsilon}

\def\g{\gamma}
\def\i{\iota}
\def\k{\kappa}
\def\l{\lambda}
\def\th{\theta}
\def\w{\omega}

\def\ni{\noindent}
\def\ss{\vspace{.1in}}
\def\bs{\vspace{.2in}}
\def\ms{\vspace{.15in}}

\def\ssn{\vspace{.1in}\noindent}

\newtheorem{theorem}{Theorem}[section]

\newtheorem{prop}[theorem]{Proposition}
\newtheorem{lemma}[theorem]{Lemma}

\newtheorem{cor}[theorem]{Corollary}
\newtheorem{defn}[theorem]{Definition}

{\addtocounter{theorem}{1}%
{\bf\\[2ex]Remark \arabic{section}.\arabic{theorem} }}{\\[2 ex]}

\def\qed{{\hfill\vrule height10pt width10pt}\\ \vspace{3 ex}}

\def\pf{{\bf Proof}: }

\def\adp{\mbox {$Ad\: P$}}

\def\A{{\cal A}}
\def\B{{\cal B}}

\def\M{{\cal M}}

\def\bfb{{\bf b}}
\def\D{{\cal D}}
\def\gplus{ G^A_\oplus}

\def\lapplus{\lap^A_\oplus}

\def\normaax{(|a_X|+|\a_X|)}
\def\normaay{(|a_Y|+|\a_Y|)}

\def\sqll{|\log\l|^{1/2}}
\def\tx{\tilde{X}}
\def\tX{\tilde{X}}
\def\ty{\tilde{Y}}

\def\tw{\tilde{W}}

\def\bdotb{ |\bfb_X||\bfb_Y| }

\def\bdota{ |\bfb_X|\normaay }

\def\adotb{ |\bfb_Y|\normaax }

\def\adota{ \normaax\normaay }

\def\baab{ \left( \bdota+\adotb \right) }

\def\bwedgeb{ |\bfb_X\wedge\bfb_Y| }

\def\awedgea{\left( |a_X||\a_Y| + |\a_X||a_Y| + |\a_X\wedge\a_Y| \right)}

\def\natural{{\rm I\kern-.18em N}}

\def\integer{{\rm Z\kern-.32em Z}}
\def\chix{{\raise.5ex\hbox{$\chi$}}}

\def\ep{\epsilon}
\def\real{{\bf R}}
\def\complex{\kern.1em{\raise.47ex\hbox{
	    $\scriptscriptstyle |$}}\kern-.40em{\rm C}}
\def\quat{{\bf H}}

\def\vs#1 {\vskip#1truein}
\def\hs#1 {\hskip#1truein}
\def\s{\sin(\theta)}
\def\c{\cos(\theta)}
\def\t{\theta}
\def\nd{\noindent}
\def\calb{{\cal B}}

\def\calm{{\cal M}}
\def\tcalm{{\tilde\calm_{k+1}}}
\def\l{{\lambda}}

\def\std{{\rm std}}

\date{\today}
\title{\bf Simple Type and the Boundary of Moduli Space}
\author{David Groisser\thanks
{Research supported in part by National
Science Foundation grant DMS-9307648}
\\ {\small Department of Mathematics}
\\ {\small University of Florida}
\\ {\small Gainesville FL 32611--8105}
\\ {\small {\it groisser@math.ufl.edu}}
\and Lorenzo Sadun\thanks{Research supported in part by National
Science Foundation grant DMS-9626698, an NSF
Mathematical Sciences Postdoctoral Fellowship and Texas ARP 
Grant 003658-037}
\\ {\small Department of Mathematics
}
\\ {\small University of Texas}
\\ {\small Austin, TX 78712}
\\ {\small {\it sadun@math.utexas.edu}}}

\begin{document}

\maketitle

\vskip -0.5in

\begin{abstract}

We measure, in two distinct ways, the extent to which the boundary
region of moduli space contributes to the ``simple type'' condition of
Donaldson theory. Using a geometric
representative of $\mu($pt$)$, the boundary region of moduli space
contributes $6/64$ of the homology required for simple type,
regardless of the topology or geometry of the underlying 4-manifold.
The simple type condition thus reduces to the interior of the $k+1$st
ASD moduli space, intersected with two representatives of (4 times)
the point class, being homologous to 58 copies of the $k$-th moduli
space.  This is peculiar, since the only known embeddings of the
$k$-th moduli space into the $k+1$st involve Taubes gluing, and the
images of such embeddings lie entirely in the boundary region.

When using de Rham representatives of $\mu($pt$)$, the boundary region
contributes 1/8 of what is needed for simple type, again regardless of
the topology or geometry of the underlying 4-manifold.  The difference
between this and the geometric representative answer is surprising but
not contradictory, as the contribution of a fixed region to the
Donaldson invariants is geometric, not topological.

\end{abstract}


\vspace{.2in}





\section{Introduction}\label{sect_intro}

This paper is a study in the geometry and topology of anti-self-dual
Yang-Mills moduli spaces.  Although moduli spaces were studied
extensively for their own sake in the 1970s and early 1980s, in the
late 1980s and early 1990s such studies were primarily a means to an
end.  Moduli spaces were studied to compute Donaldson invariants, and
Donaldson invariants were computed for their applications to
classifying smooth 4-manifolds.  Seiberg-Witten theory has, of course,
made that last road obsolete.  It is believed that the Seiberg-Witten
invariants determine the Donaldson invariants, and the former are far
easier to handle.

However, Seiberg-Witten theory has opened up new uses for Donaldson
theory.  From Seiberg-Witten theory, we now have a much better understanding
of Donaldson invariants.  Instead of using moduli spaces as a tool for
computing Donaldson invariants, we can now use Donaldson invariants
as a tool for understanding moduli spaces.  This paper is an exercise 
along those lines.  

A basic problem in four-dimensional gauge theory is to understand the
``simple type'' condition.  In Donaldson theory, a manifold is said to
have simple type if its Donaldson invariants satisfy a certain
recursion relation ([KM]; see (\ref{eqno2}) below).  In Seiberg-Witten
theory, a manifold has simple type if it has no Seiberg-Witten classes
of nonzero index.  The two notions of simple type are believed to be
equivalent, so that theorems proved about one form of simple type
should yield information about the other.

In this paper we work with the Donaldson theory sense of simple type,
examining what simple type implies about the geometry of
anti-self-dual moduli spaces.  In two ways---with intersection theory
and with de Rham theory---we measure the extent to which the boundary
region of moduli space contributes to the simple type recursion
relation.  Our results imply that the anti-self-dual moduli spaces
associated to any manifold of simple type have a very surprising
interior geometric structure.  Widely satisfied sufficient conditions
are known for a manifold to be of simple type [KM], and it is
conjectured that indeed {\it all} 4-manifolds with $b_+>1$ are of
simple type. Hence our results apply to a great many manifolds.

Simple type says that the $k+1$st moduli space $\calm_{k+1}$,
intersected with certain varieties, has the homology of a certain
multiple of the $k$-th moduli space $\calm_k$.  Our intersection
theory approach is based on the construction in
\cite{sadun1} of a geometric representative of $\mu$ of a point (see
below).  Using this representative we show that the portion of (a
small perturbation of) $\calm_{k+1}$ near the boundary 
contributes $6/64$ of the homology required for simple type,
regardless of the topology or geometry of the underlying 4-manifold.
(For a quick, heuristic derivation of this 6/64, see \cite{sadun2}.)
Simple type thus reduces to a statement relating $\calm_k$ to
nontrivial structure in the {\it interior} of $\calm_{k+1}$ (unless
our small perturbation of $\M_{k+1}$ is drastically unfaithful
topologically, which seems highly unlikely).  This is surprising,
since the only known relations between $\calm_k$ and $\calm_{k+1}$
involve Taubes patching, and relate $\calm_k$ to the boundary of
$\calm_{k+1}$.

Our second approach is to use differential form representatives of the
images of the $\mu$ map.  One then takes the wedge product of these
forms and integrates over $\calm_{k+1}$.  If we restrict the domain of
integration to a neighborhood of the boundary of $\M_{k+1}$,
we can reinterpret the simple type condition in terms of the integral
of a certain 8-form over a submanifold that represents the space of
``bubble parameters'' in the neighborhood of a background connection
in $\M_{k}$.  We show that, again independently of the topology and
geometry of the base manifold, this integral has precisely 1/8 the
value of what one would naively expect if simple type were captured
purely by a neighborhood of the boundary.  Thus again simple type
becomes a statement about the nontrivial structure of the interior of
moduli space.

It is curious but no contradiction that the two approaches yield the
different numerical answers 6/64 and 1/8. While the Donaldson
polynomial is topological, hence independent of the choice of
geometric or de Rham representatives, the contribution of each region
of moduli space is geometric, and need not be the same for two
different representatives.  Indeed, the de Rham and geometric
representative calculations not only disagree on the contribution of
the boundary region, but also disagree on how close to the boundary
the essential contributions are.  In terms of the small parameter $L$,
described below, the geometric representative picks up contributions
from bubbles of size $O(L^2)$, while the bulk of the support of the de
Rham representative is on bubbles of size $O(L)$.

%

To state our results more precisely, we must review the definition of
the Donaldson invariants, and of simple type.  Let $N$ be an oriented
4-manifold, let $G=SU(2)$ or $SO(3)$, and let $\calb^*_k$ be the space
of irreducible connections (up to gauge equivalence) on $P_k$, the
principal $G$ bundle of instanton number $k$ over $N$. Let $\calm_k
\subset \calb_k$ be the space of irreducible connections on $P_k$ with
anti-self-dual curvature, modulo gauge transformations.  We will
frequently omit the index $k$.

Donaldson [D1, D2] defined a map $\mu: H_i(N,{\bf Q}) \to
H^{4-i}(\calb_k^*,{\bf Q})$, $i=$0, 1, 2, 3, whose image freely generates
the rational cohomology of $\calb_k^*$.  Donaldson invariants are then
defined by pairing the fundamental class of $\calm_k$ with products of
$\mu$ of the homology classes of $N$, where $k$ is chosen so that the
dimensions match.  Formally, for elements
$[\Sigma_1],\ldots,[\Sigma_n] \in H_*(N)$, we write
\be\label{eqno1}
D([\Sigma_1] \ldots [\Sigma_n]) = \mu([\Sigma_1]) \smile \cdots \smile
\mu([\Sigma_n]) [\calm_k].
\ee
Now let $x$ be the point class in $H_0(N)$, and let $\omega$ be any
formal product of classes in $H_*(N)$.  The simple type condition is
that, for all $\omega$,
\be\label{eqno2}
 D(x^2\omega) = 4 D(\omega). 
\ee

Of course, the ``fundamental class of $\calm_k$'' is usually not well
defined, as $\calm_k$ is typically not compact.  The usual way to make
sense of (\ref{eqno1}) and (\ref{eqno2}) is with geometric
representatives.  One finds finite-codimension varieties $V_\Sigma$ in
$\calb_k^*$ that are Poincar\'{e} dual to $\mu([\Sigma])$.  One then
counts points, with sign, in $V_{\Sigma_1} \cap \cdots \cap
V_{\Sigma_n} \cap \calm_k$.  To make a topological invariant one must
show that the number of intersection points is independent of
auxiliary data, such as the metric and the choice of representatives.
This requires careful analysis of the bubbling-off phenomena that make
$\calm_k$ noncompact.

To compute the left-hand side of (\ref{eqno2}) we need a variety that
represents $\mu$ of the point class $x$.  In general $\mu(x)$
is not an integral class in $H^4(\calb^*)$ and
strictly speaking has no geometric representative. 
However, $-4\mu(x)$ {\it is} an integral class.  For $p\in N$, let
\be\label{eqno3}
 \nu_p = \{ [A] \in \calb^*_{k+1} | F_A^- \hbox{ is reducible at }p \}. 
\ee
Here $F_A^- = (F_A - *F_A)/2$ is the anti-self-dual part of the
curvature $F_A$, and by ``reducible at $p$'' we mean that the
components $F_{ij}^-(p)$ are all colinear as elements of the Lie
algebra of $G$.  In \cite{sadun1} it was shown that $\nu_p$ is a geometric
representative of $-4 \mu([p])$.  The simple type condition can
therefore be rewritten as
\be\label{eqno4} \#(\calm_{k+1} \cap \nu_p \cap \nu_q
\cap V_\omega)  = 64 \#(\calm_k \cap V_\omega), 
\ee
where $p$ and $q$ are any two points in $N$, $\omega$ is an arbitrary
formal product of homology cycles of $N$, and $V_\omega$ is a
geometric representative of $\mu(\omega)$.

More formally, one can write
\be\label{eqno5}
[\calm_{k+1} \cap \nu_p \cap \nu_q] = 64 [\calm_k].  
\ee
Strictly speaking, the left hand side is an element of
$H_*(\calb_{k+1})$, while the right hand side is in $H_*(\calb_k)$.
However, $\calb_k$ and $\calb_{k+1}$ are homotopy equivalent spaces,
and their homology classes may be identified.

The first question studied in this paper is this: {\it Suppose $p$
and $q$ are extremely close points in $N$, separated by a distance
$2L$.  How many of the points on the left hand side of (\ref{eqno4}) lie
near the boundary of $\calm_{k+1}$?} The answer is quite simple:

\begin{theorem}\label{thm0.1}
Let $N$ be a compact oriented Riemannian four-manifold of arbitrary
topology and geometry and let $4k \geq 3 b_+ + 5$.  Fix a coordinate
patch on $N$, and let $p$ and $q$ be the points $(\pm L, 0,0,0)$. 
Fix $\omega\in{\rm
Sym}^*(H_*N)$, $K>0$, and $\a \in (0,2)$.  Let $\tcalm^0$ be the
portion of the (perturbed) moduli space $\tcalm$ consisting of a
background of charge $k$ and a charge-one bubble of size $\l < KL^\a$.
For generic choices of geometric representatives $V_\omega$ of
$\mu(\omega)$, and for all sufficiently small $L$,
the intersection number of $\tcalm^0$ with $V_\omega \cap \nu_p \cap
\nu_q$ is $6 D(\omega)$.
\end{theorem}

The perturbed moduli space $\tcalm$ is constructed, and the genericity
conditions specified, in \S\S \ref{sect_pert}--\ref{sect_donald}.
This theorem is restated, more precisely, as Theorem \ref{thm5.1}.  In
this theorem, and throughout \S\S 2--4, we assume $k$ is in the
indicated ``stable range'' to avoid contributions from lower strata of
the compactified moduli space.

Simple type thus reduces to a statement that the {\it interior} of
$\calm_{k+1} \cap \nu_p \cap \nu_q$ be homologous to 58  copies
of $\calm_k$.  This is striking, since in general very little is known
about the interior of $\calm_{k+1}$.  As noted earlier, 
the only known embeddings of
$\calm_k$ into $\calm_{k+1}$ involve Taubes patching, and have an
image near the boundary of $\calm_{k+1}$.

On the level of differential forms, the de Rham-theoretic version of
the $\mu$-map is represented by a map
\be
\mu_d:\Omega^i(N)\to\Omega^i(\B_{k+1}^*), \ \ \ i=0,\dots,4;
\ee
the argument of $\mu_d$ is a form representing the Poincar\'{e}
dual of the argument of $\mu$.  In particular, $\mu(x)$
is represented by a 4-form $\mu_d(\w)\in H^4(\B_{k+1}^*)$ for any $\w\in
\Omega^4(N)$ with $\int_M\w=1$.  

One can write down an explicit formula for such a representative
$\mu_d(\w)$ by appealing to Chern-Weil theory on the canonical
$SO(3)$-bundle ${\cal P}\to\B_{k+1}^*\times M$ (see \S
\ref{sect_formsintro}).  Furthermore given $p\in M$, if we replace
$\w$ by $\d_p$, a delta-form supported at a point $p$, then the
resulting form on $\B_{k+1}^*$ is still de Rham cohomologous to a form
obtained using smooth $\w$ (although there is an important difference
that we will discuss later). Let us write
$\mu_d(p):=\mu_d(\d_p)$. Note that for smooth $\w\in\Omega^4(N)$ we
have
\be\label{forthmstc}
\left.\mu_d(\w)\right|_A = \int_N \left.\mu_d(p)\right|_A\w(p).
\ee

If $Z$ is an eight-dimensional cycle in $\B_{k+1}^*$, then since the
cohomology class of $\mu_d(p)$ is independent of $p$, for any points
$p,q\in M$ we have $\int_Z
\mu_d(p)\wedge \mu_d(p)=\int_Z \mu_d(p)\wedge\mu_d(q)$, and moreover
this integral depends only on the homology class of $Z$.  

Pretend, for a moment, that the moduli spaces $\M_{k+1},\M_k$ are the
total space and base space of a compact, connected, oriented fiber
bundle $\pi:\M_{k+1}\to\M_k$; the fibers would then be mutually
cohomologous compact submanifolds $Z\subset \M_{k+1}$.  For any form
$\phi\in\Omega^{\rm top}(\M_k)$, we would have a product formula
\be\label{orientations}
\int_{\M_{k+1}}\mu_d(p)\wedge\mu_d(q) \wedge \pi^*\phi=
\left(\int_Z \mu_d(p)\wedge\mu_d(q) \right) \left(\int_{\M_k} \phi\right)
\ee
(assuming compatible orientations), so the simple-type condition
(\ref{eqno2}) would
be equivalent to
\be\label{st2}
\int_Z \mu_d(p)\wedge\mu_d(q) =4.
\ee
In reality the moduli spaces are not compact and
there is no such global fibration.  However, from the current
understanding of Donaldson invariants one might speculate that the
relevant integrals are supported in a region near the ideal boundary
of $\M_{k+1}$.  Of course a random de Rham representative of a
cohomology class can be supported wherever it likes, but
$\mu_d(\cdot)$ is not random, and there is evidence that its
properties near the boundary of moduli space do indeed capture a lot
of cohomological information.  For example, consider the
five-dimensional moduli spaces of 1-instantons over simply connected
manifolds with $b_+=0$.  In such cases the inverse of a collar map
gives embeddings $\tau_\l:N\to
\M_1\subset\B^*$ for $\l$ sufficiently small (the image of $\tau_\l$
consisting of instantons of scale $\l$), and one has Donaldson's
theorem that the composition $\tau_\l^*\circ \mu: H_2(N,\bfz)\to
H^2(N,\bfz)$ is precisely Poincar\'{e} duality ([DK Corollary 5.3.3]).
The corresponding assertion in de Rham cohomology would be that
$\tau_\l^*\circ\mu_d: \Omega^2(N)\to\Omega^2(N)$ induces the
identity on cohomology.  But in fact in this context one can show that
$\lim_{\l\to 0}\tau_\l^*\circ\mu_d$ is already the identity map {\em
on the level of forms}, in all degrees (\cite{notes}). This fact
extends to more general moduli spaces $\M_k$ near the ``maximal
bubbling'' ends consisting of instantons with $k$ distinct charge-one
bubbles. 

If we consider that portion $\M'_{k+1,\l_0}$ of $\M_{k+1,\l_0}$ near
the highest-dimensional boundary stratum $\M_{k}\times N$, there is
indeed a fibration $\M'_{k+1,\l_0}\to\M'_{k}$ whose fibers can be
identified with subsets of an 8-dimensional space of framed ASD
connections on $\bfr^4$.  (Here $\M'_k$ denotes the space of
non-concentrated irreducible $k$-instantons, and $\M'_{k+1,\l_0}$ the
space of $(k+1)$-instantons with only a single ``bubble'', of charge
one, and scale less than some small number $\l_0$.)  The typical fiber
$Z=Z_{\l_0}$ is itself a bundle over $(0,\l_0)\times N$ for some small
$\l_0$, whose fiber over $(\l,p)\in (0,\l_0)\times N$ is the space of
``gluing parameters'' ${\rm Hom}_{SO(3)}(\Lambda^2_+T^*N,\adp_k) \iso
SO(3)$ (see [DK, p. 324]).  Since $\M'_{k+1,\l_0}$ is such a large
portion of the end of $\M_{k+1}$, one might then expect that an
approximate version of (\ref{st2}) holds under the assumption of
simple type.

What we show below is that (\ref{st2}) fails in a very precise way:
independent of the topology and geometry of $N$, if $\l_0$ is small
enough and ${\rm dist}(p,q)$ is small compared to $\l_0$ (but
nonzero), then
\be\label{stfails}
\int_{Z_{\l_0}} \mu_d(p)\wedge\mu_d(q)\approx \frac{1}{2}
\ee
under certain technical but intuitively reasonable assumptions about
the fiber $Z_{\l_0}$.  Taking a limit as $q\to p$ and then as $\l_0\to
0$, the integral above approaches an integral over the space of
framed instantons on $\bfr^4$, and this latter integral has the
precise value 1/2.

At this stage the reader may wonder why we do not simply take $p=q$ in
(\ref{stfails}).  The reason is that for purposes of integration, the
$\mu_d(p)$ turn out to be more singular than the representatives
$\mu_d(\w)$ for smooth $\w$. Were we to set $p=q$ in (\ref{stfails}),
1/2 would be replaced by 0.  This discontinuity can be modeled by the
following two-dimensional example.  Let $H$ be the upper half-plane
$\{(x,\l) \in \bfr^2 \mid \l>0\}$ and for each $L\in
\bfr$ let $\theta_L: H\to (0,\pi)$ be the usual polar-coordinate angle
as measured from $(L,0)$ (so $d\theta_L = ((x-L)d\l - \l
dx)/((x-L)^2+\l^2)$).  As forms on $H$, the $d\theta_L$ are all
cohomologous (in fact cohomologous to zero).  However, $\int_H
d\theta_0\wedge d\theta_0 =0$, while for $L> 0$ we have we have
$\int_H d\theta_0\wedge d\theta_L = \pi^2/2$.  Essentially,
$\mu_d(p)\wedge\mu_d(q)$ behaves like a quaternionic version of this example.


The technical assumptions on the fiber $Z$ are enumerated as ({\bf
Z1}--{\bf Z5}) in section \ref{sectfiber}.  The first three of these
assumptions are known to be satisfied by the fiber constructed in
[DK], but we have not determined whether the latter two are satisfied.
These two are assumptions on the tangent space to $T_{[A]}Z$, where
$[A]\in Z$, and we prove that they are satisfied by a subspace of
$T_{[A]}\M$ (the ``approximate tangent space'') that we argue is close
to $T_{[A]}Z$.  Because this step is only a plausibility argument,
(\ref{stfails}) implies one of two things: either $\mu_d(p)\wedge
\mu_d(q)$ has most of its support in the interior of moduli space (or
near higher codimension boundary strata), or the intuitive picture of
the fiber $Z$ is significantly wrong.  Either way, the conclusion is
surprising.

Our second main theorem is then

\begin{theorem}\label{thmst}\label{thm0.2}
 Let $N$ be a compact oriented Riemannian four-manifold of arbitrary
topology and geometry and let $k\geq 1$.  Assume that a typical fiber
$Z_{\l_0}$ of the fibration $\M'_{k+1,\l_0}\to \M'_k$ satisfies ({\bf
Z1}--{\bf Z5}) of \S
\ref{sectfiber}.  Then for any $p,q\in N$, the form $\mu_d(p)\wedge\mu_d(q)$ is
integrable over $Z_{\l_0}$ for $\l_0$ sufficiently small, and
\be\label{thmsta}
\lim_{\l_0\to 0}\left(\lim_{q\to p} \int_{Z_{\l_0}} \mu_d(p)\wedge\mu_d(q)
\right) = \frac{1}{2}, 
\ee
while
\be\label{thmstb}
\lim_{\l_0\to 0} \int_{Z_{\l_0}} \mu_d(p)\wedge\mu_d(p)
= 0.
\ee
The convergence in (\ref{thmsta})--(\ref{thmstb}) is uniform in $p,q$.
Hence if $\tilde{\d}_{p,L}$ denotes a smooth 4-form on $N$ of total
integral 1, supported in a ball of radius $L$ about $p$, then (using
(\ref{forthmstc})) 
\be\label{thmstc}
\lim_{\l_0\to 0}\left(\lim_{L\to 0} \int_{Z_{\l_0}}
\mu_d(\tilde{\d}_{p,L} )\wedge\mu_d(\tilde{\d}_{p,L})
\right) = \frac{1}{2}.
\ee
\end{theorem}

By uniform convergence in (\ref{thmsta}) we mean that 
for all $\e>0$ there exist $\l_1, \d(\cdot)>0$ such if $0<\l_0<\l_1$
and $0<{\rm dist}(p,q)<\d(\l_0)$ then the integral in (\ref{thmsta})
differs from $1/2$ by less than $\e$.

It is not necessary to take the limits in (\ref{thmsta})
completely independently, as long as $q\to p$ much faster than
$\l_0\to 0$.  If, for example, we require that ${\rm dist}(p,q)={\rm
const}\ \l_0^{1+\a}$ for some $\a>0$, and then take a limit as
$\l_0\to 0$, we again get 1/2.  

Note that if we held $p$ and $q$ fixed rather than taking $\lim_{q\to
p}$ in (\ref{thmsta}), the limit as $\l_0\to 0$ would necessarily be
zero (since $\mu_d(p)\wedge\mu_d(q)$ is integrable).  It turns out
that for $q\neq p$ the integrand in (\ref{thmsta}) is supported in a
region in which $\l$ is of order ${\rm dist}(p,q)$.  Thus if we wish
to extend $\mu_d(p)$ and $\mu_d(\tilde{\d}_{p,L})$ to forms on the
Uhlenbeck compactification of $\M$, with $\lim_{L\to 0} \mu_d(\tilde{\d}_{p,L})
=\mu_d(p)$ in a distributional sense, then $\mu_d(p)\wedge\mu_d(p)$
should be viewed as the sum of a delta-form supported on the boundary
of moduli space and a smooth form supported away from the boundary.

Note also that Theorem \ref{thmst} does not require $k$ to be in the
``stable range'', unlike Theorem \ref{thm0.1}.  However, Theorem
\ref{thmst} is most interesting for $k$ in the stable range, since 
only then can the Donaldson invariant $D([\Sigma_1] \cdots
[\Sigma_n])$ be expressed as a topologically invariant integral
$\int_\calm \mu_d(\omega_1)\wedge \cdots \wedge \mu_d(\omega_n)$. 

The rest of this paper is organized into two main parts, with \S\S 
\ref{sect_model}--\ref{sect_donald} devoted to proving Theorem
\ref{thm0.1} and \S\S  \ref{sect_formsintro}--\ref{sectpfthmrem2}
devoted to proving Theorem \ref{thm0.2}.  The strategy of proof, and
the division of the paper, is as follows:

Let $A$ be a connection obtained by gluing a small bubble onto a
background connection $A_0$.  
It turns out that the curvature of $A$ is well approximated by
the sum of the curvature $F_0$ of $A_0$ and the curvature $F_{\std}$
of a standard $k=1$ instanton, viewed in the correct gauge.  We are
thus led to the following model problem: {\it Given a connection
$[A_0] \in\calm_k$ and two closely spaced points $p$ and $q$, for how
many triples $(x,\l,g)$ is the sum of the curvature $F_0$ of $A_0$ and
the curvature $F_\std$ of a standard instanton, centered at $x$ with
size $\l$ and gluing angle $g$, reducible at both $p$ and $q$?}  
In \S \ref{sect_model} we solve this model problem and show that, for
generic $A_0$, the answer is 6.  

In \S \ref{sect_pert} we construct a family of approximately ASD
connections, based on an explicit gluing formula.  This is a
perturbation, which we denote $\tcalm$, of the boundary region of
$\calm_{k+1}$.  We check explicitly that in this family the curvature
is well approximated by $F_0+F_\std$.  By linearly interpolating
between $F_0+F_\std$ and the actual curvatures of connections in
$\tcalm$, we show that corresponding to each generic $A_0 \in \calm_k$
there are exactly 6 points in $\tcalm\cap \nu_p \cap
\nu_q$ with $\l$ sufficiently small.

In \S \ref{sect_donald} we apply these results to show that if we consider
only the boundary region of the (perturbed) moduli space, we obtain
(\ref{eqno5}) with 6 on the right-hand side rather than 64, thereby
completing the proof of Theorem \ref{thm0.1}.

Ideally, one would then like to interpolate from $\tcalm$ to
$\calm_{k+1}$.  This is quite difficult, as $\nu_p$ and $\nu_q$ are
defined by pointwise conditions on the curvature.  We know of no
pointwise estimates relating the curvature of an almost-ASD connection
to that of a nearby ASD connection.  In order to make use of the
integral estimates available in the literature one would have to
replace $\nu_p$ and $\nu_q$ by geometric representatives defined by
integral conditions.  While certainly possible, this is beyond the
scope of this paper.

We prove Theorem \ref{thmst} by exhibiting $\mu_d(p)$ as a purely
local piece $\muloc(p)$ plus a nonlocal remainder. The local piece 
dominates in (\ref{thmsta}): as $q\to p$ the integral of 
$\muloc(p)\wedge\muloc(q)$ approaches a calculable integral on
$\bfr^8$, with value 1/2, independent of $\l_0$.  (However,
$\muloc(p)\wedge\muloc(p)\ident 0$.) We establish
(\ref{thmsta})--(\ref{thmstb}) by showing that the
integral of the remainder terms in $\mu(p)\wedge\mu(q)$ approaches
zero as $\l_0\to 0$, independent of $p$ and $q$.  Thus taking a limit
as $q\to p$ is relevant only to the purely local part of
$\mu_d(p)\wedge\mu_d(q)$ (and taking a limit as $\l_0\to 0$ is
relevant only to the nonlocal part); the delta-form behavior of
$\mu_d(p)\wedge\mu_d(p)$ is due solely to
$\muloc(p)\wedge\muloc(p)$.   The uniformity assertion in Theorem
\ref{thmst} follows from the proofs of (\ref{thmsta})--(\ref{thmstb}),
and the final assertion (\ref{thmstc}) then follows from
(\ref{forthmstc}).

   In \S \ref{sect_formsintro} we begin our work on Theorem
\ref{thm0.2} by constructing the de Rham
representatives $\mu_d(p)$.  The splitup $\mu_d(p)=\muloc(p)+$
remainder is based on the ``approximate tangent space'' mentioned
above.  This approximation is built by lifting the action of certain
vector fields on $N$ to $\B^*$.  In \S \ref{sectaction} we discuss
this lifted action (the ``canonical flow''), use it to define the
approximate tangent spaces ${\cal H}_A$, and discuss how close the
${\cal H}_A$ are to being tangent to $\M$. We then exhibit the
relation between a fiber constructed from the canonical flow (whose
tangent space is essentially the projection to $T_{[A]}\M$ of
approximate tangent space above) and the fiber constructed in [DK].
This digression is needed to motivate the technical assumptions ({\bf
Z1}--{\bf Z5}) given and discussed in \S \ref{sectfiber}.  In
\S \ref{sectlocal} we return to the main track, defining
$\muloc(p)$ and computing the limiting integral of
$\muloc(p)\wedge\muloc(q)$.  \S\S
\ref{sectnonlocal}--\ref{sectpfthmrem2}  
are devoted to a study of the remainder terms
$\mu_d(p)\wedge\mu_d(q)-\muloc(p)\wedge \muloc(q)$.  In \S
\ref{sectnonlocal} we state the main technical theorem that yields the
pointwise norm of these terms (Proposition \ref{thmrem2}), and use
this theorem to establish that the integral of the remainder terms
tends to zero as $\l_0\to 0$.  Finally, in \S
\ref{sectpfthmrem2}, we prove Proposition \ref{thmrem2}.  It is this
section that contains the core of the analysis underpinning the
validity of all the earlier calculations.  The estimates in \S
\ref{sectpfthmrem2} require a weighted Sobolev inequality, proven in the
appendix, that the authors have not seen elsewhere.

The authors thank the 1994 Park City Mathematics Institute, where this
work was begun, the National Science Foundation, and the Texas
Advanced Research Program. We also thank Dan Freed, Tom Parker, Cliff
Taubes, and Karen Uhlenbeck for helpful insights and criticism, and
Margaret Combs for assistance with the figures.



\section{The Model Intersection Theory Calculation}\label{sect_model}
\setcounter{equation}{0}

In this section we begin to compare the boundary region of
$\calm_{k+1} \cap \nu_p \cap \nu_q$ with $\calm_k$, by looking at a
model problem.  Pick a small neighborhood $\tilde U$ of our manifold
$N$ and give it a flat metric with corresponding Euclidean
coordinates.  Let $U$ be the corresponding ball in $\real^4$.  We will
denote points either by four real coordinates $(x^0,\ldots,x^3)$ or by
a single quaternionic coordinate $x^0 + i x^1 + j x^2 + k x^3$.  Let
$p$ and $q$ be the points $(\pm L, 0,0,0)$.  Let $A_0$ be an ASD
connection on $N$, expressed in a smooth gauge on $\tilde U$.

An important notational tool is the identification of ASD curvatures
with $3 \times 3$ real matrices.  Let $F_0$ be the pullback, to $U$,
of the curvature $F_{A_0}$ of an ASD connection on $\tilde U$.
Relative to the standard oriented basis of $\Lambda^2_-T^*\real^4$
($\omega_1=dx^0 dx^1 - dx^2 dx^3$, $\omega_2= dx^0 dx^2 - dx^3 dx^1$,
$\omega_3=dx^0dx^3 - dx^1 dx^2$), $F_0$ has at each point 3
Lie-algebra-valued components, and so can be viewed as a triple of
3-vectors.  We package this triple of vectors into a $3 \times 3$ real
matrix, which we denote $Mat(F_0)$.  More precisely, the first, second
and third columns of $Mat(F_0)$ are half the $\omega_1$, $\omega_2$
and $\omega_3$ components of $F_0$, while the first, second and third
entries of each column refer to the three directions in
$\mathfrak{su}(2)$, the Lie Algebra of $SU(2)$.  $A_0$ is reducible at
a point if and only if $Mat(F_0)$ has rank 1 (or 0) there.

Of course, this
construction is dependent on gauge and a choice of basis for $TN$. A
gauge transformation is a change of basis in $\mathfrak{su}(2)$, and
thus changes $Mat(F_0)$ by left-multiplication by an element of $SO(3)$.
A change of basis in $TN$ changes $Mat(F_0)$ by right-multiplication by
an element of $SO(3)$.  Thus the singular values of $Mat(F_0)$, and in
particular the rank of $Mat(F_0)$, are gauge- and basis-independent.  We
shall frequently be thinking of curvatures as $3 \times 3$ matrices in
this way.  When the context is clear, we will omit the explicit
function ``$Mat$''.

Now consider a standard $k=1$ instanton on $\real^4$ of scale size
$\lambda$ and center $y$, viewed in a radial gauge that is singular at
$y$ and regular at $\infty$.  There are many such gauges, parametrized
by a gluing angle $m \in SO(3)$.  For fixed ($y,\l,m$), let $F_\std$
be the curvature of this connection restricted to $U$.

Let $A$ be an ASD connection obtained by gluing in a bubble with data
($y,\l,m$) to the background $A_0$.  In \S \ref{sect_pert} we shall
see that $F_A$, in an appropriate gauge, is approximately equal to
$F_0 + F_\std$.  This reduces our main question to the following model
problem:

{\it When $L$ is small, for what values of $(y,\l,m)$, with $\l$
small, is $F_0 + F_\std$ reducible at both $p$ and $q$?}

Of course, to obtain sensible answers, we must define what we mean by
$\l$ being ``small''.  Pick constants $K>0$ and $\a \in (0,2)$.  We
say $\l$ is small (or that the corresponding bubble is small) if $\l <
K L^\a$.  The set of gluing data for small bubbles near $p$ and $q$ is
$B = U \times (0, K L^\a) \times SO(3)$.  Let $\tilde
\nu_p$ (resp.  $\tilde \nu_q$) be the set of points $(\l , y,m) \in B$
such that $F_0(p) + F_\std(p)$ (resp. $F_0(q) + F_\std(q)$) is
reducible.  We must count the intersection points of $\tilde \nu_p$
and $\tilde \nu_q$.

Recall that the singular values $\sigma_1 \ge \sigma_2 \ge \sigma_3
\ge 0$ of a $3 \times 3$ real matrix $M$ are the square roots of the
eigenvalues of $M^T M$.  For $M$ generic, these are distinct and
positive.  The non-generic cases are as follows: Matrices in a
codimension-1 set have $\sigma_3=0$. Matrices in a codimension-2 set
either have $\sigma_1=\sigma_2$ or $\sigma_2=\sigma_3$.  Matrices in a
codimension-4 set have $\sigma_2=\sigma_3=0$; these
matrices have rank 1 or 0. Matrices in a codimension-5 set have
$\sigma_1=\sigma_2=\sigma_3$; these are all scalar multiples of
$SO(3)$ matrices.  Only the zero matrix (codimension-9) has
$\sigma_1=\sigma_2=\sigma_3=0$.

\begin{theorem}\label{thm3.1}
 Fix $K>0$, $\a \in (0,2)$, and a background
connection $A_0$.    
If the singular values of $Mat(F_0(0))$ are all distinct, then, for all
sufficiently small $L$, $\tilde \nu_p$ and $\tilde \nu_q$ intersect at
exactly six points.
These six intersections are all transverse, and the local
intersection number is $+1$ at each point.
\end{theorem}

\nd Remark: We shall see that, under the assumptions of the theorem, 
the intersection points all have $\l=O(L^2)$.  However, when
two of the singular values of $Mat(F_0(0))$ are the
same, then there are only four intersection points with $\l=O(L^2)$.
In that case there are generically four additional 
intersection points with
$\l=O(L)$. 
The intersection number of $\tilde \nu_p$ and $\tilde \nu_q$ 
is then 4 if $\a>1$ and 8 if $\a<1$.

Before beginning the proof of Theorem \ref{thm3.1} we need some basic
facts about $k=1$ instantons on $\real^4=\quat$, we need to fix some
conventions, and we need a linear algebra lemma.  Think of $SU(2)$ as
the unit quaternions, with $\mathfrak{su}(2)$ as the imaginary
quaternions.  The connection form of a standard instanton of scale
size 1, centered at the origin, is $A_{\std_0} = Im(\bar x dx /
(1+|x|^2))$.  The curvature of this connection is
\be\label{eqno3.1}
F_{\std_0} = {d \bar x dx \over (1 + |x|^2)^2} = 
{2i \omega_1  +  
2j \omega_2  +  
2k \omega_3
\over(1 + |x|^2)^2}. 
\ee
Note that the matrix
$Mat(F_{\std_0})$ is $1/(1 + |x|^2)^2$ times the identity matrix.

That is in the usual regular gauge, in which $A \sim \phi^{-1} d \phi$ as
$|x| \to \infty$, where $\phi(x)=x/|x|$.  We do a gauge transformation
by $\phi^{-1}$, to get a radial gauge in which $A= O(|x|^{-3})$ as
$|x| \to \infty$ (and in which $A$ is singular at the origin).  
We then do a further gauge transformation by a
constant $g_0$ to get the most general radial gauge with this
property.  Let $F_\std$ be the curvature form
in this gauge.  We have $F_\std = g_0^{-1} \phi F_{\std_0} \phi^{-1}
g_0$.  In terms of matrices, $Mat(F_\std) = \rho(g_0^{-1}) \rho(\phi)
Mat(F_{\std_0})$, where $\rho$ is the standard double covering map from
$SU(2)$ to $SO(3)$; the three columns of $\rho(\phi)$ are $\phi i
\phi^{-1}$,  $\phi j \phi^{-1}$, and $\phi k \phi^{-1}$. The matrix
$\rho(g_0)$ is our gluing angle $m$.

Now suppose that we have a $k=1$ instanton, centered at a point $y$,
with scale size $\lambda$.  The curvature matrix, expressed in the
exterior radial gauge of gluing angle $m$, is
\be\label{eqno3.2}
Mat(F_\std(x)) = {\lambda^2 \over (\lambda^2 + |x-y|^2)^2}\;  
m^{-1} \rho \left ( {x-y \over |x-y|} \right ).
\ee
Note that the matrix $Mat(F_\std(x))$ is a positive multiple of an $SO(3)$
matrix.  The multiple is determined by $\l$ and $|x-y|$, while the
$SO(3)$ matrix is determined by $m$ and $(x-y)/|x-y|$.  
(We hencefore will not explicitly distinguish between a curvature and
its matrix.)  

Our problem is thus one of adding positive multiples of $SO(3)$
matrices to $F_0(p)$ and $F_0(q)$ to make them reducible.  The
following lemma is essential.

\begin{lemma}\label{lemma3.2}
Let $P$ be a 3 by 3 real matrix with singular
values $\sigma_1 \ge \sigma_2 \ge \sigma_3 \ge 0$.  If these singular
values are all distinct, then there are exactly two pairs $(s,M)\in
(0,\infty) \times SO(3)$ for which
$P+sM$ has rank 1 (and no pairs $(s,M)$ for which $P+sM=0$).  
In both cases $s =\sigma_2(P)$.  
If exactly two of the singular values of $P$ are the same and
nonzero, then the two solutions $(s,M)$ coalesce to a double root.
\end{lemma}

\nd\pf Let $W=-(P+sM)$. Adding $sM$ to $P$ to make it
reducible is the same as decomposing $-P$ as $sM + W$, with $W$
reducible. We therefore count the ways to decompose a matrix $-P$ into
the sum of a positive multiple of an $SO(3)$ matrix and a rank 1
matrix. First we show that the desired decompositions can occur {\it only}
with $s = \sigma_2$, by assuming a decomposition $-P= s M + W$ and
computing $\sigma_2(P)$.  Multiplying $P$ on the left and right by
$SO(3)$ matrices does not change the singular values, but does allow
us to set $M=I$ and put $W$ into the form
\be\label{eqno3.3}
 W = \left ( \matrix{a & b & 0 \cr 0 & 0 & 0 \cr 0 & 0 & 0} \right ). 
\ee
Then
\be\label{eqno3.4}
P^TP = \left ( \matrix{(s + a)^2 & (s+a)b & 0 \cr (s+a)b & s^2 +
b^2 & 0 \cr 0 & 0 & s^2} \right ).  
\ee
One of the eigenvalues of $P^TP$ is obviously $s^2$, with eigenvector
$(0,0,1)$.  Restricting to the upper left 2 by 2 block, we subtract $s^2
I$ and get a matrix whose determinant, $-s^2b^2$, is nonpositive.
Thus at most one eigenvalue of $P^TP$ is greater than $s^2$ and at
most one eigenvalue is less than $s^2$.  Since $s^2$ is the middle
eigenvalue of $P^TP$, $\sigma_2(P)=s$.

Next we show that $P+sM$ {\it can} have rank 1, with
$s=\sigma_2(P)$, in two ways.  By multiplying on the left and right by
$SO(3)$ matrices, we can take $P$ diagonal, with entries $P_{11} \ge
P_{22} \ge |P_{33}|$.  Next we look for orthogonal matrices of the
form
\be\label{eqno3.5}
 M_\t = \left ( \matrix {-\c & 0 & \s \cr 0 & -1 & 0 \cr \s & 0 & \c} 
\right ). 
\ee
We then have 
\be\label{eqno3.6}
P+sM_\t = P + P_{22}M_\t = \left ( \matrix 
{P_{11}-P_{22} \c & 0 & P_{22} \s \cr 0 & 0 & 0 \cr P_{22} \s
 & 0 & P_{33}+ P_{22} \c} \right ).
\ee
This matrix has an obvious null vector $(0,1,0)$.  $P+sM_\t$ has rank
one (or zero) if, and only if, there is a second null vector.  To see
if there is a second null vector, we restrict $P+sM_\t$ to the 1-3
plane and take its determinant, which equals
$-P_{22}^2 + P_{11}P_{33} + (P_{11}-P_{33})P_{22} \c$.  This is a
periodic function of $\theta$ with a single maximum of
$(P_{11}-P_{22})(P_{22}+P_{33})$ at $\theta=0$ and a single minimum of
$-(P_{11}+P_{22})(P_{22}-P_{33})$ at $\theta=\pi$.  If $P_{11} > P_{22}
> |P_{33}|$, the maximum and minimum values have opposite signs, so the
function must cross zero exactly twice, at the points $\theta=\pm
\cos^{-1}([P_{22}^2 - P_{11}P_{33}]/(P_{11}-P_{33})P_{22})$.  If
$P_{11}=P_{22}$ or $P_{33}=-P_{22}$, then the maximum value becomes zero,
while if $P_{22}=P_{33}$ then the minimum becomes zero.  In these
cases we have a double root at $\theta=0$ or $\pi$.  Finally, if
$P_{11}=P_{22}=P_{33}$, then the function is identically zero and we
have an infinite number of roots.  This corresponds to the original 
$P$ being a positive multiple of an $SO(3)$ matrix.
 
Finally, we show that these are the only possible decompositions with
$s=P_{22}$.  Suppose that $M$ is an $SO(3)$ matrix with $P+sM$ having
rank 1.  Then every 2 by 2 block of $P+sM$ has determinant 0, and in
particular the upper left 2 by 2 block has a null vector $v$.
However, $P_{11}$ and $P_{22}$ are both at least $s$, so $|Pv| \ge s$.
The upper left corner of $sM$ has operator norm at most $s$, so $|sMv|
\le s$.  Thus we must have $|Pv|=s|Mv|=s=P_{22}$.  If $P_{11}>P_{22}$
this means $v=(0,1,0)$, so $Mv = (0,-1,0)$, so $M$ must take the form
(\ref{eqno3.5}).  The case $P_{11}=P_{22}$ must be checked separately,
but leads only to the solution $M=diag(-1,-1,1)$.  \qed

The form of the explicit solutions found above also demonstrates the
continuous dependence of $M$ on $P$.  Expressed invariantly, 
$M$ is a rotation by $\pi$ about an axis.  This axis is orthogonal to
the second principal axis of $P^TP$, and makes an angle $\theta/2= 
(\pm 1/2) \cos^{-1}([\sigma_2^2 \pm
\sigma_1\sigma_3]/[(\sigma_1 \pm
\sigma_3)\sigma_2])$ with the third principal axis of $P^TP$, where
the $\pm$ is determined by the sign of the determinant of $P$.  A
small change in $P$ can only change $\theta$ by an amount of order
$|\delta P|/ \min(\sigma_1-\sigma_2, \sigma_2-\sigma_3)$, and, by
first order perturbation theory (integrated to get rigorous bounds),
can only change the principal axes of $P^TP$ by a similar amount.  Thus if
$\delta P$ is a small perturbation of $P$, the norm of the
corresponding $\delta M$ is bounded by a constant times $|\delta P|/
\min(\sigma_1-\sigma_2, \sigma_2-\sigma_3)$.

Not surprisingly, this stability breaks down when we approach the
double root.  If two of the singular values are equal, then a small
perturbation may change $M$ by as much as $O(\sqrt{|\delta P|})$.

\nd {\bf Proof of Theorem \ref{thm3.1}}:  Let $s_p$ be the second
singular value of $F_0(p)$, and let $M_p\in SO(3)$ be a matrix such
that $F_0(p) + s_p M_p$ is reducible (with similar definitions for
$s_q$ and $M_q$).  Let $s_0$ be the second singular value of
$F_0(0)$.  Note that $s_0>0$, since the three singular values of
$F_0(0)$ were assumed distinct.  Since $s_p$ and $s_q$ are within
$O(L)$ of $s_0$, we can bound $s_p$ and $s_q$ away from
zero.

We shall count the ways to simultaneously make
$F_\std(p) = s_p M_p$ and $F_\std(q) = s_q M_q$. The condition
for the standard curvature $F_\std$ to have magnitude $s_p$ at $p$ is
\be\label{eqno3.7} {\l^2 \over (|y-p|^2 + \l^2)^2} = s_p, \ee
or equivalently 
\be\label{eqno3.8}  \l^2 + |y-p|^2 = \l/\sqrt{s_p}. \ee
%
%
%
As long as $|y-p| < 1/2\sqrt{s_p}$ there are two solutions to
(\ref{eqno3.8}), while for $|y-p| > 1/2\sqrt{s_p}$ there are none.
When $|y-p|< 1/2\sqrt{s_p}$, one solution has $\l > 1/2 \sqrt{s_p}$,
which is greater than $K L^\a$ for $L$ small.  The other solution 
qualifies as small if $|y-p|$ is small enough, and, for $|y-p|\ll
1/\sqrt{s_p}$, is approximately $\l = |y-p|^2 \sqrt{s_p}$.  As a set
in $\real^5 =(N,\l)$ space, the solutions to (\ref{eqno3.8}) are a 4-sphere.
Projected onto $N$, they form (2 copies of) a 4-disk.  In either case,
only a small subset of solutions qualifies as ``small''.

The interesting question, of course, is how many times we can solve
the equations for $p$ and $q$ simultaneously.  We begin with equation
(\ref{eqno3.8}) and the corresponding equation for $q$.  The
intersection of two 4-spheres in $\real^5$ is a 3-sphere.  Projected
onto $N$ we get a 3-dimensional ellipsoid, possible degenerating to
two disks.  As before, only a small patch of the ellipsoid (or
alternatively part of one of the two disks) gives a small enough value
of $\l$.  It is this region that we consider.

Recall that $p$ and $q$ are at $\pm L$, where we are using
quaternionic coordinates.  For $L$ small, $s_q = s_p + O(L)$.  Let
$s_m$ be such that $2/\sqrt{s_m} = 1/\sqrt{s_p} + 1/\sqrt{s_q}$.  Let
$\Delta = (1/\sqrt{s_p} \! - \! 1/\sqrt{s_q})/L$.  As $L \to 0$, $s_m= s_0
+ O(L^2)$, while $\Delta$ approaches $-(ds_p/dL)|_{L=0}/s_0^{3/2}$.  Let
$y_0$ and $y_I$ be the real and imaginary parts of $y$.  Adding and
subtracting (\ref{eqno3.8}) to the corresponding equation for $q$ we obtain
\be\label{eqno3.9} 
-4y_0 = \lambda \Delta ; \qquad \lambda^2 + L^2 + y_0^2 + |y_I|^2 = 
\lambda /\sqrt{s_m}.
\ee
Plugging the first equation into the second we get
\be\label{eqno3.10}
\lambda^2 \left ( 1 + {\Delta^2 \over 16} \right ) - {\lambda \over 
\sqrt{s_m}} + L^2 + |y_I|^2 = 0 
\ee
This equation shows that $\l$, and thus $y_0$, may be viewed as
functions of $y_I$.  As long as $L^2 + |y_I|^2 \ll 1/\sqrt{s_m}$ there
are two solutions to (\ref{eqno3.10}), one of which has $\lambda
\approx (L^2 + |y_I|^2)\sqrt{s_m}$, the other of which has $\lambda
\approx ((1+\Delta^2/16)\sqrt{s_m})^{-1}$.  The first solution has $\l
< KL^\a$ if and only if $|y_I|$ is small enough, while the second
solution always has $\l > KL^\a$.  Let $R_{K,\a}$ be the largest
number such that $|y_I|<R_{K,\a}$ implies $\l\leq K L^\a$.  Henceforth
we consider only ``admissible'' $y$, i.e. those with $|y_I|<R_{K,\a}$.
For $L$ chosen small enough, as we assume henceforth it is,
$R_{K,\a}^2 \sim KL^\a/\sqrt{s_m} - L^2 \sim KL^\a/\sqrt{s_m}$, since
$\a < 2$.
Note that 
\be\label{displayy0}
y_0 = -\lambda \Delta/4
\approx - (L^2 + |y_I|^2)\sqrt{s_m}
\Delta /4.
\ee
Hence for admissible $y$, we have $|y_I|<{\rm const}L^{\a/2}$ and
$|y_0|<{\rm const}\ L^{\a}$.  Let $r = (y_0(0),0)$ be the unique
admissible point where the ellipsoid of solutions $(y_0,y_I)$ to
(\ref{eqno3.9}) hits the real axis.  Since $|r|=O(L^2)$, $r$ lies on the
line segment $\overline{pq}$, and the ellipsoid has curvature $O(1)$
at $r$.  See Figure 1. 

\vbox{\medskip
\centerline{\epsfysize=3truein\epsfbox{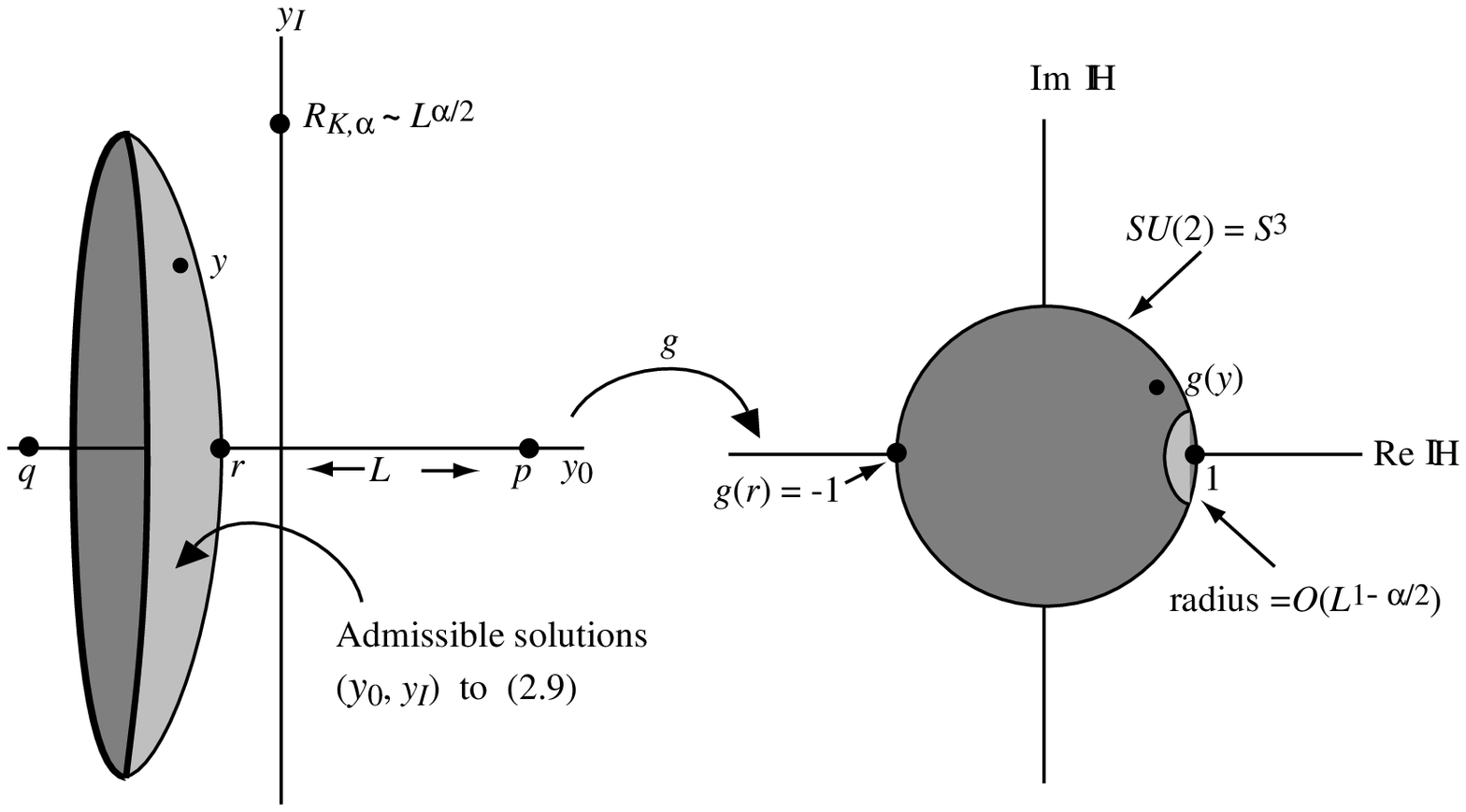}}
\medskip
\centerline{\bf Figure 1}
}
\bigskip

We still have to get the $SO(3)$ matrices right.  This means
simultaneously solving the equations $m^{-1} \rho((y\! - \!p)/|y\! -
\!p|)= M_p$ and $m^{-1} \rho((y\! - \!q)/|y\! - \!q|) = M_q$ for $m$.
If a solution exists it is obviously unique.  A solution exists if and
only if $\rho((y\! - \!p)/|y\! - \!p|)^{-1}\rho((y\! - \!q)/|y\! -
\!q|) = M_p^{-1} M_q$.  Let $g(y) = (\bar y\! - \! \bar p)(y \! - \!
q)/|(y\!  - \!p)(y\! - \!q)|$.  We must count the points on our 3-disk
(of small solutions to (\ref{eqno3.9}) and (\ref{eqno3.10})) for which
the $SO(3)$-valued function $\rho(g(y))$ equals $M_p^{-1} M_q$.  Note
that
\be\label{eqno3.11}
g(y) = -I+ 2\frac{y_I}{L}(1+O((|y_0|/L)^2)) + O((|y_I|/L)^2) 
\quad
\hbox{for }|y_I| \ll L,
\ee
while 
\be\label{eqno3.12}
g(y) = I + 2\frac{Ly_I}{|y_I|^2}(1+O((|y_0|/|y_I|)^2))
 + O((L/|y_I|)^2) \quad \hbox{for }|y_I| \gg L. 
\ee
In view of (\ref{displayy0}), we can replace $O((|y_0|/L)^2)$ in
(\ref{eqno3.11}) and $O((|y_0|/|y_I|)^2)$ in (\ref{eqno3.12}) 
by $O(L^2)$ and $O(L^{2\a})$, respectively.

Observe that $L/R_{K,\alpha}$ is $O(L^{1-\alpha/2})$ and hence goes to
zero as $L \to 0$.  On the disk of admissible $y_I$, the map $g$
covers all of $SU(2)$ except for a ball of radius $c L^{1-\alpha/2}$
around the identity, for some constant $c$.  Since $\rho$ is a 2-1 map,
$\rho(g(y))$ hits all of $SO(3)$ twice, except for a ball of radius $2
c L^{1-\alpha/2} + O(L^{2-\a})$ around the identity, which is only hit
once.  The number of solutions to our problem depends on whether, for
small $L$, $M_p^{-1}M_q$ is in this ball or not.

If the singular values of $F_0(0)$ are distinct, then, by Lemma
\ref{lemma3.2}, there are two distinct matrices $M_{1,2}(0)$ for which
$F_0(0) +\sigma_2(0) M$ has rank 1.  By the comment after the proof of
Lemma \ref{lemma3.2}, the two matrices for $p$ and $q$ satisfy
$M_{1,2}(p,q) = M_{1,2}(0) + O(L)$.  As $L \to 0$, $M_1(p)^{-1}
M_2(q)$ and $M_2(p)^{-1} M_1(q)$ are bounded away from the identity,
but $M_1(p)^{-1} M_1(q)$ and $M_2(p)^{-1} M_2(q)$ are within $O(L)$
(and hence within $2 c L^{1-\alpha/2} + O(L^{2-\a})$) of the identity.
Thus we have two configurations in $(y,\l,m)$ space that give $s_p
M_1(p)$ at $p$ and $s_q M_2(q)$ at $q$, two that give $s_p M_2(p)$ at
$p$ and $s_q M_1(q)$ at $q$, one that gives $s_p M_1(p)$ at $p$ and
$s_q M_1(q)$ at $q$ and one that gives $s_p M_2(p)$ at $p$ and $s_q
M_2(q)$ at $q$.  A total of six solutions in all.

On a codimension-2 set of background data, the background curvature at
the origin has two equal singular values, so $M_1(0)=M_2(0)$ and
$M_{1,2}(p,q) = M_1(0) + O(L^{1/2})$.  In that case all 4
possibilities have $M_p^{-1}M_q = 1 + O(L^{1/2})$.  If $\a>1$, this is
within $2 c L^{1-\alpha/2}$ of the identity for small enough $L$, and
so each possibility yields one solution.  If $\a<1$, and the
$O(L^{1/2})$ term in the expansion of $M_{1,2}(p,q)$ in powers of $L$
is nonzero, then each possibility yields two solutions. 
 
Finally we consider the orientations of our solutions.  It is not
immediately clear that all solutions have the same orientation, but in
fact they do.  The problem of matching amplitudes is the same in all
cases.  The problem of matching gluing angles reduces to the
intersection of two 3-cycles in a 3-disk $\times$ $SO(3)$ (i.e., all
possible pairs $(y_I,m)$), and is easily seen to be transverse.  The
intersection numbers are continuous functions of $M_p$ and $M_q$, as
long as a solution continues to exist.  Sending $M_p$ around a
noncontractible loop in $SO(3)$ interchanges the two solutions
associated to a given pair $(M_p,M_q)$, which shows that the two
solutions for any given $(M_p, M_q)$ have the same orientation.  Also
by continuity, this orientation is the same for all pairs $(M_p,M_q)$.

All that remains is to compute this orientation in one case.  Let
$s_p=s_q=1$, $M_p=M_q=I$, and look near the solution with $y=0$ and
$m=I$.  The varieties $\tilde \nu_p$ and $\tilde
\nu_q$ are just the zero sets of $F_\std(p)-I$ and $\tilde
F_\std(q)-I$, which we view as functions of $(y,\l,m)$. Taking
derivatives we find that changes in $(y,\l,m)$ give the following
first order changes in $F_\std(p)$ and $F_\std(q)$:
\begin{enumerate}
\item  Increasing $\l$ increases the magnitude of both $\tilde
F_\std(p)$ and $F_\std(q)$ without changing either direction.
\item  Increasing $y_0$ increases the magnitude of $F_\std(p)$
and decreases that of $F_\std(q)$, while keeping the directions
fixed.
\item  Increasing $y_1$ (resp. $y_2$, $y_3$) rotates $\tilde
F_\std(p)$ in the direction defined by the Lie algebra element $-i$
(resp. $-j$, $-k$), and rotates $F_\std(p)$ an equal amount in the
direction $+i$ (resp. $+j$, $+k$).
\item  Rotating $m$ in any direction rotates both $F_\std(p)$ 
and $F_\std(q)$ in the opposite direction.
\end{enumerate}

 From this we deduce that the Jacobian $|d(F_\std(p),
F_\std(q)) /d (y,\l,m)|$ is positive, and so that the local 
intersection number of $\tilde \nu_p$ and $\tilde \nu_q$ is $+1$ in
this case. Thus the local intersection number of $\tilde \nu_p$ and
$\tilde \nu_q$ is $+1$ in all cases.
\qed

Having proven Theorem \ref{thm3.1}, we consider the question of
stability.  How much do our intersection points move around if we
change $M_p$ or $M_q$ or $s_p$ or $s_q$ slightly?  Since $F_0 +
F_\std$ is only an approximation to the true curvature of a connection
in $\calm_{k+1}$, our results must be stable in order to be
meaningful.

Let $\chi$ be the map that takes $(y,\l,m)$ to $(F_\std(p),
F_\std(q))$.   Near our solutions, $d \chi$ is 
never close to singular.  By changing $\l$ and one component of $y$ we
can adjust $|F_\std(p)|$ and $|F_\std(q)|$ independently,
while by adjusting $m$ and the remaining three components of $y$ we
can adjust the directions of $F_\std(p)$ and $F_\std(q)$
independently.  It is not difficult to estimate the matrix elements of
$(d\chi)^{-1}$.  Some are $O(1)$, some are $O(L)$, and some are
$O(L^2)$.  If we know the required $F_\std(p$ or $q)$ to within
$\ep$, we know $m$ to within $O(\ep)$, $y$ to within $O(\ep L)$, and
$\l$ to within $O(\ep L^2)$.  In short, small errors in the input data
result in only small changes of the locations of our intersection
points in $(y,\l,m)$ space.

Finally, we consider a perturbation of our model problem that is more
directly applicable in the sequel.  Let $\tilde F_0(x)$ be the
curvature of the background connection in the standard radial gauge
about the gluing point $y$. (That is, use the original connection
$A_0$ to trivialize the fiber over $y$, and then use parallel
transport radially outwards from $y$ to trivialize the bundle over
$U$.) We wish to count the number of ways to make $\tilde F_0 +
F_\std$ reducible at both $p$ and $q$.  

\begin{theorem}\label{thm3.3}
 Under the assumptions of Theorem \ref{thm3.1},
the number of ways to make $\tilde F_0 + F_\std$ reducible at $p$
and $q$ (counted with sign) is the same as the number of ways to make
$F_0 + F_\std$ reducible at $p$ and $q$ (counted with sign),
namely +6.
\end{theorem}

\nd\pf We first put our background connection into a radial
gauge with respect to the origin.  This is a fixed gauge, and Theorem
\ref{thm3.1} applies.  Since $F_0$ and $\tilde F_0$ are related by a 
gauge transformation, the singular values of $F_0$ and $\tilde F_0$
are the same.  Thus we must solve (\ref{eqno3.8}), and the
corresponding equation for $q$, exactly as in Theorem \ref{thm3.1},
with the same values of $s_p$ and $s_q$.  We then solve $m^{-1}
\rho((y\! - \!p)/|y\! - \!p|)= M_p$ and $m^{-1} \rho((y\! - \!q)/|y\!
- \!q|) = M_q$ for $m$, as before.  The only difference in our
analysis is that $M_p$ and $M_q$ are now functions of $y$.  We compute
the extent to which they depend on $y$.

Let $A$ be a connection in radial gauge with respect to $y$, and let
$A'$ be the same connection in radial gauge with respect to $y'$.  The
gauge transformation that relates these, evaluated at the point $p$,
is the holonomy around a triangle from $p$ to $y$ to $y'$ to $p$, and
so its difference from the identity is bounded by the sup norm of
$|F_{A}|$ times the area of the triangle. See Figure 2.

\vbox{\medskip
\centerline{\epsfysize=2.5truein\epsfbox{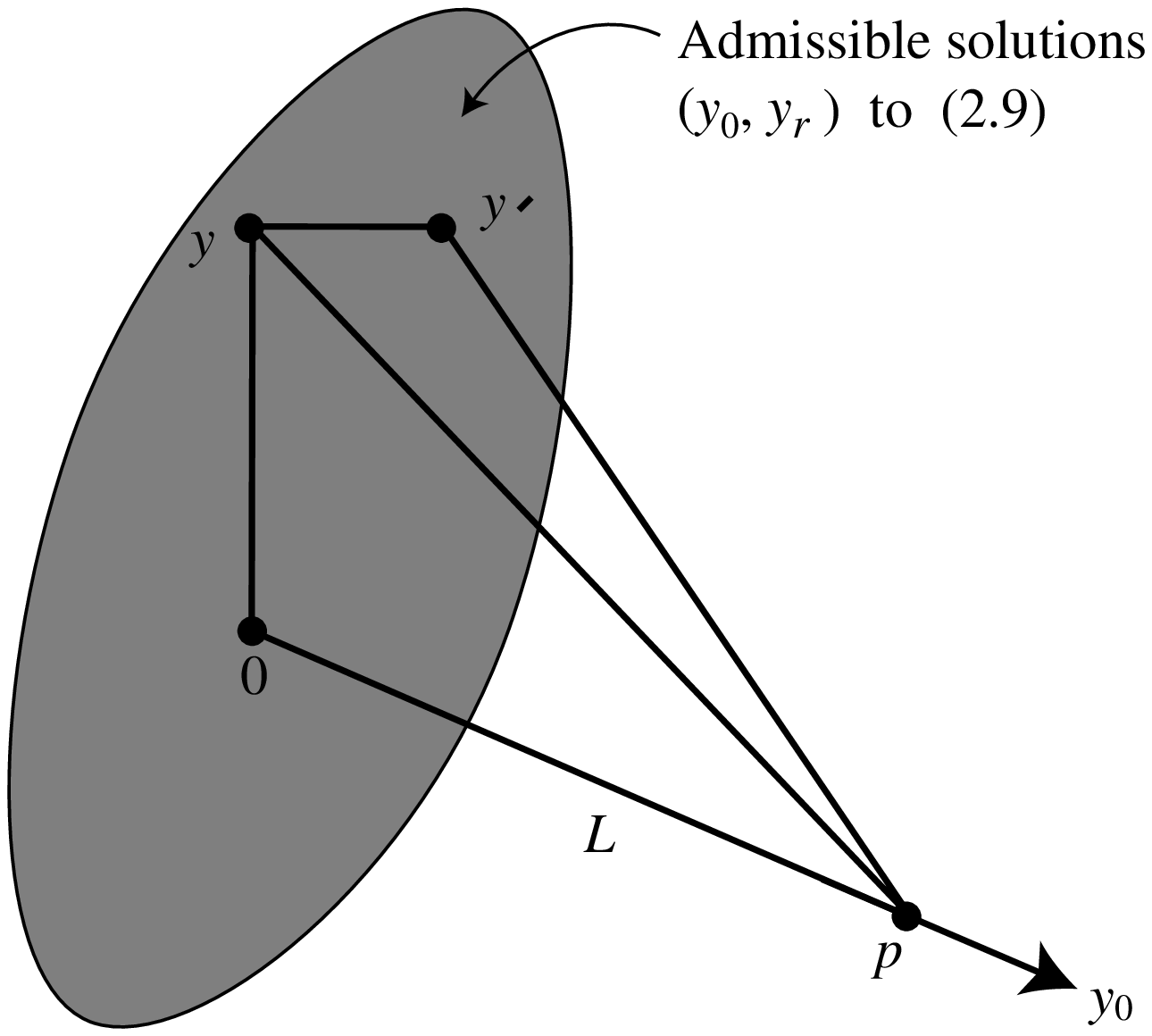}}
\medskip
\centerline{\bf Figure 2}
}
\bigskip

In our case, $A$ is the background connection, so $|F_A|$ is fixed and
bounded, and $y$ and $y'$ are restricted to lie on the ellipsoid of
solutions to (\ref{eqno3.9}) and (\ref{eqno3.10}), with $|y_I|$ and
$|y'_I|$ both less than $R_{K,\a}$.  Note that the area of a triangle
is bounded by half the product of the length of {\it any} two of its
legs. Because the curvature of the
ellipsoid of solutions $y=(y_0,y_I)$ is $O(1)$ at admissible points,
$|y_0-y'_0|$ is bounded by a constant times $|y_I-y'_I|$. As a result,
\be\label{eqno3.13}
|M_p(y)-M_p(y')| \le \hbox{const } \times \sqrt{L^2+|y|^2} |y_I -y'_I|,
\ee
while
\be\label{eqno3.14}
|M_p(y)-M_p(0)| \le \hbox{const } \times L|y_I|,
\ee
with similar estimates for $M_q$.
The second result is
an estimate on $M_p$ itself, while the first leads to a bound on the
derivative of $M$ with respect to $y_I$.  By (\ref{displayy0}),
$L^2+|y|^2 \leq {\rm const}\ (L^2+|y_I|^2)$, so we obtain 
\be\label{hideeho} 
|\partial M_p/\partial y_I|\leq {\rm const}\times \sqrt{L^2+|y_I|^2}.  
\ee

As before, we look for solutions to $\rho(g(y)) = M_p^{-1} M_q$, where
now the right hand side depends on $y$.  We break the disk of radius
$R_{K,\a} = O(L^{\a/2})$ into two pieces, a inner disk $D_1$ and an
annulus $D_2$.  The radii will be chosen such that on $D_1$ the
estimate (\ref{eqno3.13}) is strong enough to allow implicit function
theorem arguments to apply.  Here the solutions to $\rho(g(y)) =
M_p^{-1}(y) M_q(y)$ are but small perturbations of the solutions to
$\rho(g(y)) = M_p^{-1}(0) M_q(0)$.  On $D_2$ the estimate
(\ref{eqno3.14}) will be strong enough to show that there are no solutions
to $\rho(g(y)) = M_p^{-1}(y) M_q(y)$.  Taken together, this will prove
the theorem.

On the disk $D_1$, the implicit function theorem will apply as long as
the smallest singular value of $\partial(\rho\circ g)/\partial y_I$ is
at least twice the largest singular value of $\partial
(M_p^{-1}M_q))/\partial y_I$, which by (\ref{hideeho}) is bounded
above by a multiple of $(L^2 + |y_I|^2)^{1/2}$.  Computing the
derivative of $\rho \circ g$ is an easy geometrical calculation, and
one finds that all singular values are bounded below by a constant
times $L/(L^2+|y_I|^2)$. Comparing $L/(L^2 + |y_I|^2)$ to $(L^2 +
|y_I|^2)^{1/2}$ we see that the implicit function theorem applies
whenever $|y_I|$ is smaller than a constant times $L^{1/3}$, and in
particular whenever $|y_I| < L^{1/2}$ (and $L$ is sufficiently small).
We take the radius of $D_1$ (and the inner radius of $D_2$) to be
$L^{1/2}$.

Now consider $y_I \in D_2$.  If $\a>1$, then $D_2$ is empty, so we
assume $\a \le 1$.  By (\ref{eqno3.12}), $|I-\rho(g(y))| =
2L/|y_I|(1+O(L^\a)) + O(L^2/|y_I|^2)$.  Since $c_1 L^{1/2} < |y_I| <
c_2 L^{\a/2}$, $c_3 L^{1-\a/2} < |I-\rho(g(y))| < c_4 L^{1/2}$.  Now
recall that $M_p^{-1}(0)M_q(0)$ is either bounded away from the
identity or is within $O(L)$ of the identity (e.g. $M_1^{-1}(p)M_2(q)$
is bounded away from the identity while $M_1^{-1}(p)M_1(q)$ is within
$O(L)$ of the identity).  By (\ref{eqno3.14}), $M_p^{-1}(y)M_q(y)$ is
also either bounded away from the identity or within $O(L)$ of the
identity on $D_2$. Thus $|I-M_p^{-1}(y)M_q(y)|$ can never be between
$c_3 L^{1-\a/2}$ and $c_4 L^{1/2}$, so there are no solutions to
$\rho(g(y))=M_p^{-1}(y)M_q(y)$ on $D_2$. \qed



\section{The Perturbed Moduli Space}\label{sect_pert}
\setcounter{equation}{0}

In this section we show that the model problem of \S \ref{sect_model}
correctly describes the intersection of $\nu_p, \nu_q$, and a
perturbation (denoted $\tcalm$) of the boundary region of
$\calm_{k+1}$.  $\tcalm$ is parametrized by quadruples $(A_0,y,\l,m)$,
where $A_0 \in \calm_k$ is a background connection, and the glued-in
bubble has size $\l$, center $y$ and gluing angle $m$.  We construct
$\tcalm$ by an explicit gluing formula and show that, in the relevant
region, the curvature of a connection in $\tcalm$ is well approximated
by the sum of the background curvature $F_0$ and the curvature
$F_\std$ of a standard instanton of size $\l$, center $y$ and gluing
angle $m$.  Our model problem was essentially to make this sum
reducible at $p$ and $q$.  By interpolating between this sum and the
actual curvature of a connection in $\tcalm$, we show that the results
of \S \ref{sect_model} carry over almost word for word.

As before, we pick a background connection $A_0 \in \calm_k$ and
constants $K>0$ and $\a \in (0,2)$.  Let the neighborhood $\tilde U$
in $N$, and the corresponding neighborhood $U$ of the origin in
$\real^4$, be as in \S \ref{sect_model}. We now allow bubbles to be
glued in anywhere (not just in $\tilde U$), so the set $B$ of gluing
data is a $(0, K L^\a) \times SO(3)$ bundle over $N$, with local
coordinates $(y,\l,m)\in N \times (0,KL^\a) \times SO(3)$.  When the
center of the bubble is in $\tilde U$, we identify the center point in
$N$ with the corresponding coordinate in $U$, and call both points $y$.
For each $(y,\l,m)\in B$, let $F$ be the curvature of the connection
$(A_0,y,\l,m) \in \tcalm$.  The variety $\nu_p$ (resp.  $\nu_q$),
restricted to the fiber of $\tcalm$ over $A_0$, is the set of points
$(\l , y,m) \in B$ such that $F^-(p)$ (resp. $F^-(q)$) is reducible.
We must count the intersection points of $\nu_p$ and $\nu_q$.  In this
section we prove

\begin{theorem}\label{thm4.1}
Fix $K>0$, $\a \in (0,2)$, and $A_0 \in
\calm_k$.   If the singular values of $F_{A_0}(0)$ are all distinct, then, 
for all sufficiently small $L$, the intersection number of $\nu_p$,
$\nu_q$, and the fiber of $\tcalm$ over $A_0$ is +6.
\end{theorem}

We begin by constructing the space $\tcalm$.  For now, assume we are
gluing a bubble of size $\l$ in $\tilde U$, with the center point at
the origin.  There are three natural length scales determined by the
background connection $A_0$.  The first is the length scale
$|F_{A_0}(0)|^{-1/2}$ of the background curvature at the origin.  The
second is the length scale $|F_{A_0}(0)|/|\nabla^A F_{A_0}(0)|$ at
which this curvature varies.  Let $R_3$ be the smaller of these two
length scales.  Finally, let $s_0$ be the second singular value of
$F_{A_0}(0)$.  It is easy to see that $s_0 < 1/R_3^2$, but there is no
simple lower bound for $s_0$ (although, by assumption, $s_0$ is always
positive).  As we have seen, $s_0$ measures how far $F_{A_0}(0)$ is
from being reducible.

Now pick additional length scales $R_1$ and $R_2$, which can depend on
$\l$, $R_3$ and $s_0$, such that $R_1^2 < 10^{-6} \l R_3$ and $R_2^2 >
10^6 \l / \sqrt{s_0}$. When $\l \ll R_3$, which is the only case we
will consider, we want $\l \ll R_1 \ll R_2 \ll R_3$. The points of
interest $x$ will all have $R_1 < |x| < R_2$.  The number $10^6$ is of
course arbitrary.  It is just chosen large enough that we can safely
ignore small numerical factors.

Let $\beta(r)$ be a smooth monotonic function that equals zero for $r<
1/2$ and equals 1 for $r>2$, and such that $\beta' \le 1$.  We define
cutoff functions $\beta_1(x)=\beta(|x|/R_1)$ and
$\beta_2(x)=1-\beta(|x|/R_2)$.

Let $A_0$ be the background connection expressed in a smooth fixed
radial gauge with respect to the origin.  Let $A_\std$ be the
connection of a standard instanton of size $\l$ expressed in a radial
gauge that is {\it singular} at the origin and regular at $\infty$.  (This
gauge is not unique; it depends on a gluing angle $m$.  See the
discussion before expression (\ref{eqno3.2})).  Note that $|A_\std| \sim
\l^2/r^3$ for $r \gg \l$ while $|A_0| \sim r |F_{A_0}| = r/R_3^2$ for
$r \ll R_3$.

Our point $(A_0,0,\l,m) \in \tcalm$ is defined by the connection form 
\be\label{eqno4.1} A' = \beta_1 A_0 + \beta_2 A_\std. \ee
We compute
\bearray
F=F_{A'} = dA' + A' \wedge A' & = &\non  
\beta_1 F_{A_0} + \beta_2 F_{A_\std} \\
& + &\non 
(\beta_1^2 - \beta_1)A_0 \wedge A_0 + (\beta_2^2 - \beta_2)A_\std 
\wedge A_\std \\
& + & \non d\beta_1 \wedge A_0 +  d\beta_2 \wedge A_\std \\
& + & \label{eqno4.2}
\beta_1\beta_2(A_\std \wedge A_0 + A_0 \wedge A_\std),
\eearray
and the interpolating 2-form
\be\label{eqno4.3}
 F_t = t(F_{A_0}+ F_\std) + (1-t) F, 
\ee
where $0 \le t \le 1$.

If the bubble is to be glued in at a point $y$, rather than at the
origin, we must adjust the formulas as follows.  First suppose $y \in
\tilde U$.  Take $A_0$ as the connection of the background in radial gauge 
with respect to $y$ (not with respect to 0).  The quantities $s_0$,
$R_1$, $R_2$, and $R_3$ are computed from the curvature $F_{A_0}(y)$,
not $F_{A_0}(0)$.  The connection $A_\std$ is in a singular radial
gauge with respect to $y$ (not to 0).  The cutoff functions are
$\beta_1(x)=\b(|x-y|/R_1)$ and $\b_2(x)=1 - \b(|x-y|/R_2)$.  With
these modifications, we still have $A' = \b_1 A_0 + \b_2 A_\std$, and
formulas (\ref{eqno4.2}) and (\ref{eqno4.3}) still apply.  For $y \not
\in U$, just apply the same formulas, using geodesic normal
coordinates around $y$.  In this case the ``standard instanton''
$A_\std$ is no longer exactly anti-self-dual, but becomes
anti-self-dual in the $\l \to 0$ limit.  The gluing angle $m$ depends
on a local trivialization, but the set of gluing angles is
invariant. This defines the space $\tcalm$ for all $y$.

In \S \ref{sect_model} we distinguished notationally between radial
gauge with respect to 0 and radial gauge with respect to $y$, calling
the background curvature $F_0$ in the first case and $\tilde F_0$ in
the second case. Theorem \ref{thm3.1} discussed making $F_0 + F_\std$
reducible at $p$ and $q$, while Theorem \ref{thm3.3} discussed making
$\tilde F_0 + F_\std$ reducible at $p$ and $q$.  In this section the
background connection is {\it always} in radial gauge with respect to
the gluing point $y$. With only one case to consider, we always write
$F_0$, never $\tilde F_0$.
 
Note that we do {\it not} use the gluing formula found in standards
works such as [DK].  Traditionally, one takes
$A''=(1 - \b_2) A_0 + (1-\b_1)A_\std$, so that the resulting
connection is exactly flat in the annulus with radii $2R_1$ and
$R_2/2$ around $y$.  This makes identifying the bundles on which $A_0$
and $A_\std$ live conceptually easier.  However, such a procedure
makes for a perturbed moduli space on which $\nu_p$ and $\nu_q$
intersect non-transverely,  since $F^-_{A''}$ is reducible,
indeed zero, on the entire annulus $2R_1 < r < R_2/2$.  This makes the
intersection number effectively impossible to compute.  

Instead, we allow the supports of $\b_1 A_0$ and $\b_2 A_\std$ to
overlap, as in Taubes' work such as [T].  This allows us to
observe the interaction between the background connection and the
glued-in instanton.  In the [DK] method, the interaction only occurs
when we go from our explicit approximate ASD connection to the true
ASD connection (something we have relatively little analytic control
over).  In our method, the interaction is seen at the level of the
approximate connection $A'$, which we can calculate.  Moreover,
$F_{A'}^+$ is much smaller than $F_{A''}^+$ (in the $L^2$ norm), so
our method should give a closer approximation to the properties of the
true moduli space.

Let $\nu_{t,p}$ (resp. $\nu_{t,q}$) be the set of gluing data
$(y,\l,m)$ with $\l<KL^\a$ for which $F_t^-(p)$ (resp. $F_t^-(q)$) is
reducible.  If $y$ is not in $\tilde U$, then, for small enough $\l$,
the connection form near $p$ is exactly $A_0$.  By assumption, $F_0$ is
not reducible at the origin.  For small enough $L$, therefore, $F_0$
is not reducible at $p$, and $F_t(p)=F_0(p)$ is not reducible.  We
may therefore assume, without loss of generality, that our gluing
point $y$ is always in $\tilde U$.  Indeed, by picking $L$ small
enough we may assume that $y$ is in an arbitrarily small neighborhood
of the origin, and therefore that $F_0(y)$ is arbitrarily close to
$F_0(0)$.  Thus we may take the length scales $R_1$, $R_2$, and $R_3$
to be independent of $y$ (although $R_1$ and $R_2$ may depend on
$\l$).

We consider five possibilities:
\begin{enumerate}
\item {I.} $|p-y| \le  R_1/2$,  
($p$ is in the ``interior zone'', where $\b_1=0$ and $\b_2=1$),
\item {II.} $R_1/2 < |p-y| < 2R_1$, 
($p$ is in the interior ``shoulder''),
\item {III.} $2R_1 \le |p-y| \le R_2/2$, 
($p$ is in the ``plateau'', where $\b_1=\b_2=1$),
\item {IV.} $R_2/2 < |p-y| < 2R_2$, 
($p$ is in the exterior ``shoulder''), 
and
\item {V.} $|p-y| \ge 2 R_2$, ($p$ is in the ``exterior zone'', where
$\b_1=1$ and $\b_2=0$).
\end{enumerate}

As in \S \ref{sect_model}, we will be identifying curvatures with
$3\times 3$ real matrices.  The phrase ``the second singular value of
$F$'', for example, is shorthand for ``the second singular value of
$Mat(F^-)$''.

The problem of Theorem \ref{thm3.3} was to find $\nu_{1,p}$,
$\nu_{1,q}$ and count their intersection points.  In that problem
condition III always applied, with $|p-y|^2 \approx
\l/\sqrt{s_p}$.  We will show that $F_t^-(p)$ being reducible
with $\l < KL^\a$ also implies condition III, and that $\nu_{t,p}$ is 
a small perturbations of $\nu_{1,p}$. We establish condition III by
showing that the other conditions lead to contradictions.

We begin by considering condition I.  Where $|x-y| \le R_1/2$,
$A'=A_\std$ has ASD curvature, so $ F^-_t = F_t = tF_{A_0} + F_\std.$
For $F_t$ to be reducible at $p$ we need $|F_\std(p)|=ts_p$.  That is,
$\l^2 + |p-y|^2 = \l/\sqrt{ts_p}$.  This quadratic equation has two
solutions, one with $\l \approx |p-y|^2
\sqrt{ts_p}$, the other with $\l \approx 1/\sqrt{ts_p}$, but neither
is consistent with condition I.  Since $|p-y| < R_3$,
$s_p$ is close to $s_0$.  Since $|p-y|^2 < R_1^2 < 10^{-6} \l R_3$,
while $\sqrt{ts_0}
\le 1/R_3$, one cannot have $\l \approx |p-y|^2 \sqrt{ts_p}$.  The
second solution has $\l \approx 1/\sqrt{ts_p} > R_3$, which
contradicts $\l \ll R_3$.  Thus condition I is impossible.

If $p$ is in the interior shoulder we have additional terms to
consider:  
\be
F_t = F_{\std} + (t + (1-t)\beta_1)F_{A_0} +
(1-t)[(\beta_1^2-\beta_1) A_0 \wedge A_0 + d\beta_1 \wedge A_0 +
\beta_1 (A_\std\wedge A_0 + A_0 \wedge A_\std)].\ee
  The ASD part of the terms after $F_{\std}$ can be bounded in norm by
$1/R_3^2 + 4R_1/R_3^4 + 4/R_3^2 +
\l^2/R_1^2R_3^2 < 100/R_3^2$, and so the second singular value of
$F_t$ is within $100/R_3^2$ of the second singular value of $F_\std$.
For $F_t^-$ to
be reducible, $|F_{\std}|$ can be at most $100/R_3^2$.  Thus we need
$\l/(\l^2 + |p-y|^2) < 10/R_3$, which in turn means that either $\l >
R_3/100$ or $\l < 100 R_1^2/R_3$.  The first is not allowed, as $\l$
is assumed small.  The second contradicts the definition of $R_1$.  So
condition II is also impossible.

If $p$ is in the exterior zone, we have $F_t = F_{A_0} + (1-t)
F_{\std}$, so we need $\l/(\l^2 + |p-y|^2) = \sqrt{s_p/(1-t)}$, or
equivalently $\l = (\l^2 + |p-y|^2)
\sqrt{s_p/(1-t)}$.  But $|p-y|^2 > 2R_2^2 > 10^6 \l/\sqrt{s_0}$, 
so $\sqrt{s_p/(1-t)}$ always exceeds $\l/(\l^2 + |p-y|^2)$.  So again
we have a contradiction.

If $p$ is in the exterior shoulder, we have 
\be F_t = (t + (1-t)\beta_2)
F_\std + F_{A_0} + (1-t)[(\beta_2^2-\beta_2) A_\std
\wedge A_\std + d\beta_2 \wedge A_\std + \beta_2 (A_\std \wedge A_0 +
A_0 \wedge A_\std)].
\ee
The ASD parts of the terms other than $F_{A_0}$ have total magnitude
bounded by $\l^2/R_2^4 +
\l^4/R_2^6 + \l^2/R_2^4 + \l^2/R_2^2R_3^2 < 10 \l^2/R_2^4 <
 10^{-11} s_0< 10^{-10} s_p$.  But $F_{A_0}$ is a distance greater
than $s_p/2$ from the nearest reducible matrix, so $F_t^-(p)$ cannot
be reducible.

Thus for all points in $\nu_{t,p}$ condition III applies, and here the
analysis is relatively simple.  The cutoff functions are both 1, so
$F(p)=F_{A_0} + F_\std + (1-t)(A_\std \wedge A_0 + A_0 \wedge
A_\std)$.  This last term has magnitude bounded by $\l^2/R_1^2 R_3^2$,
and changes only slightly as $(y,\l,m)$ are varied.  It can thus be
treated as a perturbation of $F_{A_0}$.  We perturb $\nu_{1,p}$ to
$\nu_{t,p}$ iteratively (as in the standard proof of the inverse
function theorem): Given a point in $\nu_{1,p}$, compute $(1-t)(a
\wedge A_0 + A_0 \wedge a)$, use that to adjust $(y,\l,m)$, compute
the change in $(1-t)(A_\std \wedge A_0 + A_0
\wedge a)$, adjust $(y,\l,m)$, and so on.  The iteration converges
geometrically.  Similarly, a point in $\nu_{t,p}$ can be perturbed
to a point in $\nu_{1,p}$.  

Of course, the same analysis applies to $\nu_{t,q}$.

Now we consider the number of intersection points of $\nu_{t,p}$ and
$\nu_{t,q}$, as a function of $t$.  The only way the intersection
number can change is if intersection points appear or disappeared at
the ends of $\nu_{t,p}$ or $\nu_{t,q}$.  However, we have shown that
such intersection points can only occur when both $p$ and $q$ are in
the plateau.  In the proof of Theorem \ref{thm3.3} we saw that, for
$\l \gg L^2$ but $\l \ll 1$ (e.g., $\l \sim KL^\a$), the points of
$\nu_{1,p}$ are bounded away from $\nu_{1,q}$.  Since condition III
applies, for $\l \sim KL^\a$, $\nu_{t,p}$ and $\nu_{t,q}$ are close to
$\nu_{1,p}$ and $\nu_{1,q}$, respectively, and so are bounded away
from each other.  Thus intersection points between $\nu_{t,p}$ and
$\nu_{t,q}$ may not appear from or disappear to the boundary.  Thus
$\# (\nu_{0,p}\cap \nu_{0,q}) = \# (\nu_{1,p} \cap \nu_{1,q})$.  By
Theorem \ref{thm3.3}, the latter number is +6, regardless of $A_0$.
\qed




\section{Computing the Donaldson Invariants}
\label{sect_donald}
\setcounter{equation}{0}

In sections \ref{sect_model} and \ref{sect_pert} we saw that, for a
fixed generic background connection, there are six ways to glue in a
small bubble near $p$ and $q$ so as to make the curvature reducible at
$p$ and $q$.  In this section we demonstrate that this is sufficient
information to compute the contribution of the boundary region of
$\tcalm$ to the simple type condition.  For generic choices of
representatives (of the classes other than $\mu(p)$ and $\mu(q)$), and
for generic choice of the location of the origin of our coordinate
system, the boundary region contributes $6/64$ of what is needed for
simple type.

We continue the notation of \S\S \ref{sect_model}-\ref{sect_pert}.
$\tcalm$ is the perturbed moduli space and $\tilde U$ is a fixed ball
in $N$ with a Euclidean metric, which we identify with a neighborhood,
$U$, of the origin in $\real^4$. For fixed $K, \a, L$, let $\tcalm^0$
be the subset of $\tcalm$ with $\l < KL^\a$.  Let $\omega$ be a formal
product of cycles $[\Sigma_1], \ldots, [\Sigma_n] \in H_*(X)$ such
that $\hbox{deg}(\mu(\omega))= \hbox{dim}( \calm_k)$, so that the
Donaldson invariant $D(\omega)$ is computed on the $k$-th moduli space
$\calm_k$.

We assume that the classes $\{ [\Sigma_i] \}$ are represented by
smooth submanifolds $\{ \Sigma_i \}$ in general position. In
particular, a subset of the $\{ \Sigma_i \}$ can intersect only if
their codimensions add up to 4 or less.  Pick tubular neighborhoods
$\{ \tilde \Sigma_i \}$ of $\{ \Sigma_i \}$ small enough to have the
same property: a subset of the $\{ \tilde \Sigma_i \}$ can intersect
only if the codimensions of the corresponding $\Sigma_i$'s add up to 4
or less.  Similarly, we assume that the $\tilde
\Sigma_i$'s do not intersect our fixed ball $\tilde U$. Choose a 
geometric representative $V_i$ of each $\mu([\Sigma_i])$ that depends
only on the connection restricted to $\tilde \Sigma_i$.  This may be
done for the 1, 2, and 3-dimensional cycles as in [DK], and for the
0-dimensional $\Sigma$'s as in [DK] or \cite{sadun1}.  (This allows us to
identify the geometric representative of $\mu([\Sigma])$ on $\calb_k$
with the geometric representative of $\mu([\Sigma])$ on $\calb_{k+1}$.
In each case it is the set of connections whose restriction to $\tilde
\Sigma$ satisfies a certain condition).  Note that the codimension of
$V_i$ in $\calb$ is the codimension of $\Sigma_i$ in $N$.  Let
$V_\omega = \cap_i V_i$.  $V_\omega$ is a geometric representative of
$\mu(\omega)$.  Generically, $V_\omega$ will intersect $\calm_k$ at a
finite number of points (this number, counted with sign, {\it is} the
Donaldson invariant $D(\omega)$), and each of these points will
exhibit generic behavior.  In particular, for each such point $A_0$,
we can assume that $Mat(F_{A_0}(0))$ has three distinct singular
values.

\begin{theorem}\label{thm5.1}
Fix $\tilde U, \omega, K>0$, and 
$\a \in (0,2)$. For generic choices of $V_\omega$, as described above,  
and for all sufficiently small $L$, the intersection number of
$\tcalm^0$ with $V_\omega \cap \nu_p \cap \nu_q$ is $6
D(\omega)$.
\end{theorem}

\nd \pf
We need to show that the only way for the boundary region of $\tcalm$
to intersect $V_\omega \cap \nu_p \cap \nu_q$ is if a bubble is
pinching off near $p$ and $q$, while the background connection in
$\calm_k$ is contributing to $D(\omega)$.  We then must demonstrate
that, under these circumstances, the problem reduces to the counting
problems studied in \S\S \ref{sect_model}--\ref{sect_pert}.

Suppose we have a small bubble centered at a point $y$ that is not
in $\tilde U$.  The point $y$ can lie in at most 4 of the $\bar
\Sigma_i$'s, with the corresponding $\Sigma_i$'s having total
codimension 4 or less.  Recall that we are using the explicit formula
(\ref{eqno4.1}), and that outside a neighborhood of $y$ the new
connection is {\it identical} to the background connection.  For small
$\l$, therefore, the bubble inserted at $y$ has no effect on the
connection restricted to the remaining $\tilde \Sigma$'s (which we
index by $j$).  Therefore, for a connection $(A_0,\l,y,m)\in \tcalm$
to lie in $\cap_i V_i$, the background connection $A_0 \in \calm_k$
must lie in $\cap_j V_j$.  However, $\calm_k$ has dimension 8 less
than $\tcalm$, while $\cap_j V_j$ has dimension at most 4 more than
$\cap_i V_i$. Since the dimension of $\calm_k$ is less than the
codimension of $\cap_j V_j$, $\cap_j V_j \cap \calm_k$ is generically
empty.
  
Next we consider the case where a small bubble is centered in $\tilde
U$.  Then $\{\Sigma_j \}$ is all the cycles $\Sigma$ except the
two points $p$ and $q$.  For small $\l$, on each of the $\tilde
\Sigma_j$'s the connection form is equal to the background connection 
$A_0$, which must therefore be in $\cap_j V_j \cap \calm_k$.  However,
now the dimension of $\calm_k$ and the codimension of $\cap_j V_j$
match. $\cap_j V_j \cap \calm_k$ is, by our genericity assumption, a
discrete set of points, whose number (counted with sign) is the
Donaldson invariant $D(\omega)$.  For each of these points, the
singular values of $Mat(F_{A_0}(0))$ are distinct.  

By Theorem \ref{thm4.1}, for each such background $A_0$, and for $L$
small enough, there are exactly 6 values of $(\l,y,m)$ such that
$(A_0,\l,y,m) \in
\tcalm^0$ has reducible curvature at $p$ and $q$.
Furthermore, the intersection numbers for the local problem are all
$+1$.  Now the orientation of $\tcalm^0$ is the same as that of
$\calm_k \times U \times (0,KL^\a) \times SO(3)$ [D3, \S 3].

Thus the contribution of points $(A_0,\l,y,m)$ to $D([p]\cdot [q]\cdot
\omega)$, for fixed $A_0$, is exactly 6 times the contribution of
$A_0$ to $D(\omega)$.  Summing over the finite set $\{ A_0 \}$, we get
that the contribution of $\tcalm$ to $D([p]\cdot [q]\cdot \omega)$ is
$6 D(\omega)$.
\qed


\newpage


\section{Differential forms and the $\mu$-map:
Introduction}\label{sectintro}\label{sect_formsintro} 
\setcounter{equation}{0}

Theorem \ref{thm0.1}, restated precisely as Theorem \ref{thm5.1}, is
one of the two major results of this paper.  It quantifies the
contribution of the boundary region of moduli space to the geometric
representative computation of the Donaldson invariants that appear in
the simple type recursion relation.  The remainder of the paper is a
proof of Theorem \ref{thm0.2}, which quantifies the contribution of
the boundary region to a differential forms calculation of the same
Donaldson invariants.

In this section we construct a de Rham-theoretic version of
Donaldson's $\mu$-map using Chern-Weil theory.  Recall that there is a
canonical $SO(3)$-bundle ${\cal P}\to \B^*\times N$, and that the
$\mu$-map is defined by slant product with $-\frac{1}{4}p_1({\cal
P})$.  Using the $L^2$ metric one can produce a natural connection on
${\cal P}$, with curvature $\F$; see [DK \S\S 5.1-2].  By the
Chern-Weil formula one has
\be
-\frac{1}{4}p_1({\cal P}) =\frac{1}{8\pi^2}\tr(\F\wedge\F) \in
\Omega^4(\B^*\times N),
\ee
where the trace comes from the two-dimensional representation of
$\mathfrak{so}(3)\iso \mathfrak{su}(2)$. Let us write tangent vectors
to $\B^*\times N$ as pairs $(\a,X)$ with $\a\in
T\B^*$ and $X\in TN$, and
identify $T_A\B^*$ with $\ker((d^A)^*)\subset \Omega^1(\adp)$.  Further,
for $\a,\b\in\Omega^1(\adp)$, define $\{\a,\b\}=-\sum_{i=0}^4
[\a_i,\b_i] \in \Omega^0(\adp)$, where the local $\adp$-valued
functions $\a_i,\b_i $ are the components of $\a,\b$ relative to a
local orthonormal basis of $T^*N$.  If $A$ is irreducible, then
$\F((\a,0),(\b,0)) =-2G^A_0\{\a,\b\}$, where $G^A_0$ is the inverse of
the covariant Laplacian on $\Omega^0(\adp)$, and hence
\bearray
\left.\mu_d(\w)(\a,\b,\g,\rho)\right|_A
&=& \non \int_N \left(\i_{(\rho,0)} \i_{(\g,0)} \i_{(\b,0)}
\i_{(\a,0)} \frac{1}{8\pi^2}\tr(\F\wedge\F)\right)\w\\
&=& \non \frac{1}{\pi^2}
\int_N \tr (G^A_0\{\a,\b\}G^A_0\{\g,\rho\}
+G^A_0\{\a,\g\}G^A_0\{\rho,\b\}\\ &&
\mbox{\hspace{.5in}}+G^A_0\{\a,\rho\}G^A_0\{\b,\g\})\w.
\label{combinatorics}
\eearray
For our application it is crucial to get the combinatoric factors in
(\ref{combinatorics}) correct.

If we replace $\w$ by $\d_p$, a delta-form
supported at a point $p$, the resulting form on $\B^*$ is still de
Rham cohomologous to a form obtained using smooth $\w$.
Henceforth we write $\mu_d(p):=\mu_d(\d_p)$.
For any $p\in
N$, a 4-form representing $\mu_d(p)$ is given by
\bearray\non 
\lefteqn{\mu_d(\d_p)(\a,\b,\g,\rho)= } \\
&& \label{mud} 
\left. \frac{1}{\pi^2}
\tr \left(G^A_0\{\a,\b\}G^A_0\{\g,\rho\}
+G^A_0\{\a,\g\}G^A_0\{\rho,\b\}
+G^A_0\{\a,\rho\}G^A_0\{\b,\g\}\right)\right|_p.
\eearray

To make use of (\ref{mud}) we need some concrete formulas---with
calculable leading terms and small remainders---for $G^A_0\{\a,\b\}$.
We can obtain such formulas when $A$ is a concentrated instanton with
a ``charge-one bubble'' and $\a,\b$ come from infinitesimal changes in
the bubble parameters (center, scale, and gluing angle).  Tangent
vectors of this type span an ``approximate tangent space'' on which
very strong estimates are possible.  This space, its relation to the
action of the quaternionic affine group on $\bfr^4$, and its relation
to the gluing construction in [DK], are central to the proof of
Theorem \ref{thm0.2}.  In the next section, we define the approximate
tangent space precisely and study these relations in detail.



\section{Group actions and the approximate tangent space.}
\label{sectaction}\label{sectats}
\setcounter{equation}{0}

Let {\bf H} denote the quaternions and ${\bf H}^*$ the nonzero
quaternions.  The eight-dimensional approximate tangent space we
define later is obtained by an ``almost-action'' of ${\bf
H}^*\times{\bf H}\iso
\bfr_+\times SU(2)\times \bfr^4$ on $\B$, induced by an almost-action
on $P$ (what ``almost-action'' means is explained below).
Essentially, we lift from $N$ to $P$ cut-off versions of translations,
dilations, and ``self-dual rotations'' in a gauge-invariant way.

   To make this more precise, let $X$ be a vector field on $N$, and
$A$ a connection on $P$. The pair $(X,A)$ defines a
flow on $P$ obtained by lifting $X$ $A$-horizontally to $P$.  We
thereby obtain from $X$ the ``canonical flow of $X$
on $\A$'', with associated vector field $A\mapsto\tx:=\i_X F_A\in
\Omega^1(\adp)\iso T_A\A$ (see [GP1] Proposition 4.3). The canonical
flow is invariant under the gauge group, hence descends to $\B$.
Moreover any two lifts to $P$ of a diffeomorphism of $N$ differ
by a gauge transformation, and hence given an action on $N$ by any
connected Lie group $G$ on $N$, the canonical flow integrates to a
well-defined action of $G$ on $\B$, though in general not on $\A$.
Of interest to us later will be
the comparison of the canonical lift to that obtained by lifting $X$
horizontally with respect to a reference connection $A_0$.  In that
case the difference between the two tangent vectors in $T_A\A$ induced
by the two flows is $d^Au$, where $u=\i_X(A-A_0)$.

   Now let $G$ be a Lie group acting from the left on $N$.  Suppose
that for each $g\in G$, the action $\Phi_g$ of $g$ on $N$ lifts to a
bundle map $\tilde{\Phi}_g:P\to P$; if $G$ is connected we can obtain
such lifts by using the canonical flow.  (We do not require the lifts
to piece together to a $G$-action.) For later purposes we will need to
calculate the differential of the induced action of $G$ on $\B$, at
any $g\in G$.  This is not difficult, but it is easy to confuse the
roles of $g$ and $g^{-1}$ in this calculation, and this mistake would
be fatal for our application.

   For each connection $A\in \A$, let $\Theta_A\in \Omega^1(P,
\mathfrak{su}(2))$ denote the
connection form of $A$.  Given a lift $\tilde{\phi}_g$ as above,
define $g\cdot A$ to be the connection with connection form
$(\tilde{\Phi}_g^{-1})^*\Theta_A$.  If the lifts piece together into an
action of $G$ on $P$ (necessarily a left action), then $(g,A)\to
g\cdot A$ defines a left action of $G$ on $\A$.  
$\tilde{\Phi}_{g_1g_2}$ and $\tilde{\Phi}_{g_1}\circ
\tilde{\Phi}_{g_2}$ are gauge-equivalent, 
since both are lifts of $\Phi_{g_1g_2}$, so an element-wise
liftable $G$-action on $N$ {\em always} induces a $G$-action on $\B$,
whether or not it induces one on $P$.

Now fix $[A_0]\in \B$ and define $\overline{\rho}:G\to \B$ by
$\overline{\rho(g)}=[g\cdot A_0]$.  This is well-defined and is
independent of the choice of lifts.  On a small enough neighborhood
$U$ of any $g\in G$ we can always
choose the $\tilde{\Phi}_{h}$ to vary smoothly with $h$, so that on $U$
the map $\overline{\rho}$ factors
through a smooth map $\rho:G\to \A$ defined by $\rho(g)=g\cdot A_0$.
Let $v=\frac{d}{dt}g_t|_{t=0}\in T_gG$ and write $v={R_g}_*w$, where
$w\in T_eG=\mathfrak{g}$.  Then
\be
{\rho_*}_g v=\left.\frac{d}{dt}( (\exp(tw)g)\cdot A_0)\right|_{t=0}.
\ee
But $\tilde{\Phi}_{\exp(tw)g}= \g(t)\circ  \tilde{\Phi}_{\exp(tw)}
\circ \tilde{\Phi}_{g}$ for some gauge  transformation $\g(t)$
varying smoothly in $t$, and hence 
$(\exp(tw)g)\cdot A_0= 
((\exp(tw))\cdot g\cdot A_0)\cdot \g(t))$.
Thus 
\bearray\non
{\rho_*}_gv &=&\left.\frac{d}{dt}( \exp(tw)\cdot g \cdot
A_0)\right|_{t=0} \ \mbox{\rm mod Im} (d_{g\cdot A_0})\\
&=&\left.\frac{d}{dt}( (\tilde{\Phi}_{\exp(tw)}^{-1})^* \w_{g \cdot
A_0}\right|_{t=0} \ \mbox{\rm mod Im} (d_{g\cdot A_0}).
\eearray
Let $\hat{w}\in \Gamma(TN)$ and $\hat{w}_P\in \Gamma(TP)$
be the vector fields on $N$ and $P$ induced by the infinitesimal
action of $w$.  Then 
\be
\left.\frac{d}{dt}( (\tilde{\Phi}_{\exp(tw)}^{-1})^* \w_{g \cdot
A_0}\right|_{t=0} 
= -{\cal L}_{\hat{w}_P}\w_{g\cdot A_0}
= -\i_{\hat w}F_{g\cdot A_0} \ \mbox{\rm mod Im} (d_{g\cdot A_0})
\ee
(if $\tilde{\Phi}$ is defined by the canonical flow then the ``$
\mbox{\rm mod Im} (d_{g\cdot A_0})$'' can be erased in this line).
Note that $v$ directly defines a vector field on $N$ by
$\hat{v}|_{\Phi_g(x)}=\frac{d}{dt}(\Phi_{g_t}(x))|_{t=0}$.  Since we
can take $g_t=\exp(tw)g$, it immediately follows that
$\hat{v}|_{\Phi_g(x)}= \hat{w}|_{\Phi_g(x)}$ for all $x\in N$, so the
vector fields $\hat{v}$ and $\hat{w}$ are the same.  Hence
\be
{\rho_*}_g v = -\i_{\hat v}F_{g\cdot A_0} 
\ \ \ \mbox{\rm mod Im} (d_{g\cdot A_0}).
\ee
Thus if we identify $T_{[A]}\B$ with $\ker((d^A)^*)\subset
\Omega_1(\adp)$, then
\be\label{drho}
{\overline{\rho}_*}_g X = -\pi'_A\i_{\hat X}F_{g\cdot A_0}.
\ee
where $\pi'_A: \Omega^1(\adp)\to \ker((d^A)^*)$ is the $L^2$-orthogonal
projection.  (Here and below, for notational convenience we do not
distinguish between a tangent vector to $\B$ of $\M$ at $[A]$,
literally a gauge-invariant section of a certain vector bundle over
the gauge-orbit through $A$, and the representative of this section at
$A$.)

   We would like to apply these ideas to the situation of a local
action of ${\bf H}^*\times {\bf H}$ on a neighborhood of a point in
$N$.  Given a concentrated connection $A$, with scale $\l=\l(A)$ and
center point $p_A$, fix a positively oriented normal coordinate system
centered at $p_A$.  \footnote{The precise definitions of $\l$ and
$p_A$ are not important here.  There are several definitions in
the literature, leading to some arbitrariness in the definition of
``near'', ``bubble'', etc.  In all instances in which the differences
among these definitions have been carefully analyzed, it has been
found that these differences do not affect the estimates we need in
any material way (cf. \cite{decayest}, section 5).  We will simply
assume in this paper that the same is true here, and will freely quote
results proved using different definitions as if they had been proved
using consistent definitions of scale and center.}  Near $p_A$ it
makes sense to speak of the translation, dilation, and rotation vector
fields.  These are determined invariantly by data $(\bfb,a,\a)\in
T_{p_A}M\plus
\bfr
\plus\Lambda^2(TN)$ by setting
\be\label{defxhat}
\hat{X}_{(\bfb,a,\a)}=b^j \ddxj+(\sqrt{2}\l^{-1})(a x^i \ddxi
+ \a_{ij}x^i\ddxj)
\ee
where $\{ x^i \}$ are normal coordinates centered at $p_A$ and
$b^j,\a_{ij}$ are the associated components of $\bfb,\a$.  We include
the normalization factor $(\sqrt{2}\l)^{-1}$ to arrange
$\|\i_{\hat{X}} F_A\|_2\approx {\rm const}$ (independent of $\l$); see
Proposition \ref{lemmad1} below.  We call $\a_{ij}x^i \ddxj$ a
self-dual/anti-self-dual rotation vector field if $\a_{ij}dx^i\wedge
dx^j$ is a SD/ASD two-form at $p_A$.

   Since the $\hat{X}$ are only defined locally, we extend them to $N$
by cutting them off outside a small ball.  For this purpose $\hat{X}$
as above we define $X=\b \hat{X}$, where $\b$ is a cutoff of scale
\be\label{defeps}
\e=4n_0\l^{1/2}.
\ee
Here $n_0$ is a constant taken large enough to ensure that $\b$ can be
used in the gluing constructions of [DK], but for all of our other
applications we can ignore $n_0$.  For convenience we take $\b=\b_{\rm
std}(r_A/\e)$, where $r_A$ is distance to $p_A$ and $\b_{\rm std}$ is
a cutoff function with support in $[0,2]$, identically 1 on
$[0,1]$. (These cutoffs, which will be used for the rest of this
paper, are different from the ones in \S 3.)

We define
\be\label{frakha}
\mathfrak{h}_A = \{ X_{(\bfb,a,\a)}=\b \hat{X}_{(\bfb,a,\a)} \in 
\Gamma(TN) \mid (\bfb,a,\a)\in T_{p_A}M\plus \bfr\plus\Lambda^2_+(TN)
\} 
\ee

   It is worthwhile to observe that the definition of (A)SD rotation
vector fields is necessarily a local definition, since globally a
nontrivial exact 2-form cannot be SD or ASD on an orientable
compact manifold. In fact on $S^4$, rotations that are SD at one pole
are ASD at the other.  This is most easily seen by using stereographic
projection to identify $S^4-\{\infty\}$ with $\bfr^4$, then with ${\bf
H}$.  Left-multiplication by unit quaternions induces SD rotations
near 0, while right-multiplication induces ASD rotations near 0.  But
coordinates near ${\infty}$ on $S^4$ are related to those near 0 by
quaternionic inversion (the orientation-preserving map $x\mapsto
x^{-1}$), which interchanges the roles of left- and
right-multiplication.

When $A$ is ASD, we make the following definition.  

\begin{defn}\label{defn_of_ha}
{\rm
The {\em approximate tangent space} ${\cal H}_A$ at $A$ to
the moduli space is the space
\be \{ \tx_A:=\i_X F_A \mid X\in \mathfrak{h}_A \} . \ee
We usually write simply $\tx$ and leave the $A$-dependence implicit.
}
\end{defn}

\ss   
To justify this terminology, we consider the action induced by such
$X$ on an ASD connection.  Since the $X$'s are nearly conformal vector
fields, one expects the induced flow to map an ASD connection to a
nearly ASD connection.  Proposition \ref{lemmad1} below shows that
this is the case, and more--but first we need a definition and a
lemma.

\begin{defn}\label{restricted_M}
{\rm
Given $\k,\nu, \l_0 >0$, let 
$\M_{k+1,\l_0}^{\k,\nu}\subset \M'_{k+1,\l_0}$ denote the subset of
instantons $[A]$ obeying the conditions 

\ss
(i) the first eigenvalues of the Laplacians $(d^A)^*d^A$ on zero-forms,
$d^A_+ (d^A_+)^{\ *}$ on SD two-forms, are both greater than $\nu$,

\ss
\ni and
\ss

(ii) for all $p\in N$,
\be\label{ptwsfest}
|F_A(p)|\leq C\l^2/(\l^2+r_A(p)^2)^2 +\k, 
\ee
where $F_A$ is the curvature of $A$, $r_A(p)={\rm dist}(p_A,p)$, and $\l,p_A$
are the scale and center point of $A$, respectively. }
\end{defn}

\ms
The pointwise bound (\ref{ptwsfest}) essentially says that $|F_A|$ is
bounded by the curvature of a standard instanton plus a contribution
$\k$ from a background connection.  At small distances from $p_A$, the
latter term is negligible, but far from $p_A$ the background term
dominates.

When dealing with estimates for the approximate tangent space, one
must decide at what scale $\e$ to cut off the vector fields in
$\mathfrak{h}_A$. If one takes $\e$ too small, the derivatives of the
cutoff function become inconveniently large, while if one takes $\e$
too large, the contribution from the background connection swamps the
contribution from the concentrated curvature.  If we require that $\e$
scale as a power of $\l$, we get the optimal balance between these
undesirable features only if $\e\sim\l^{1/2}$.  Earlier we chose
$\e=4n_0\l^{1/2}$, and we now take $n_0$ 
large enough so that in the gluing construction of [DK] one is assured
of landing in the domain of Taubes' contracting mapping argument.
Once $n_0$ has been so chosen, it (like $\k$) is simply another
ignorable constant for the computations we need in this section.  In
particular, note that on the ball $B(p_A,\e)$ or radius $\e$ centered
at $p_A$, we have $\l^2/(\l^2+r_A^2)^2\geq
(\l_0+64n_0^2)^{-2}\geq {\rm const}\cdot\k$.  Hence, with a new constant
$c=c(\k)$,
\be\label{ptwsfest2}
|F_A|\leq c\l^2/(\l^2+r_A^2)^2 \mbox{\hspace{.3 in} \rm{on
supp}}(\b).
\ee
This enormously simplifies our computations.

The next lemma shows that we can always arrange the fiber $Z$ to lie
in some $\M_{k+1,\l_0}^{\k,\nu}$.   

\begin{lemma} Given small $\l_0$, let $T: \M'_{k+1,\l_0}\to \M_k$ be 
the projection sending a concentrated connection to a ``background''
connection.  Let $[A_0]\in
\M'_k$ and assume that the first eigenvalues of the Laplacians 
$d_{A_0}^*d_{A_0},d_{A_0}^+(d_{A_0}^+)^*$ on $\adp$-valued 0-forms and
SD 2-forms, respectively, are positive.
Then there exist $\l_0>0,\nu>0,\k>0,C>0$, and a neighborhood 
$\cal U$ of $[A_0]$ in $\M_{k+1}$ such that $T^{-1}({\cal U})\subset
\M_{k+1,\l_0}^{\k,\nu}$.
\end{lemma}

\pf That condition (i) in the definition of $\M_{k+1,\l_0}^{\k,\nu}$
can be satisfied follows from the proof of Lemma 7.1.24 in [DK]; that
condition (ii) can be satisfied follows from modifying several ideas
in [GP2] (Definition 4.1, Lemma 4.3d, and Propositions 4.4).  \qed

Henceforth we will always assume that instantons $[A]$ lie in a fixed
$\M^{\k,\nu}$.  For such connections we have the following.

\begin{prop}\label{lemmad1} Fix $\k,\nu$.  
Let $\pi_A:\Omega^1(\adp)\to H^1_A:= \ker((d^A)^*)\intersect
\ker(d^A_+)$ (naturally identified with the tangent space 
$T_{[A]}\M$)
be the $L^2$-orthogonal projection.
For all sufficiently small, positive $\d$
there exist $c,\e_1(\l)$, such that
if $A\in
\M_{k+1,\l_0}^{\k,\nu}$ then
\be\label{lemmad1a}
\left|\ \lnorm \tX_{(\bfb,a,\a)} \rnorm_2^2 -
8\pi^2(|\bfb|^2+|a|^2+|\a|^2) \right| \leq \e_1(\l)(|\bfb|^2+|a|^2+|\a|^2) 
\ee
where $\e_1(\l)\to 0$ as $\l\to 0$, and
\be
\lnorm \tx-\pi_A\tX \rnorm_2
\leq
c\ \l^{\d}\left(|\bfb| \l + (|a|+|\a|)\l^{1/2}\right).
\label{lemmad1b}
\ee
\end{prop}

\pf The proof of the first statement is similar to that of [GP1,
Proposition 3.6]; we omit the details.  We prove the second statement later
as Proposition \ref{dogprop2}(b). \qed

Thus by taking $\l$ small enough we can ensure that
$\pi_A:{\cal H}_A\to H^1_A$ is
injective. Let $O_0\subset\mathfrak{h}_A$ be an open neighborhood of
zero. For $X\in\mathfrak{h}_A$, let $A^X$ denote the connection that
results from acting on $A$ by the canonical flow of $X$ for unit time,
and let $O_A=\{ A^X \mid X\in O_0 \}$.  Proposition \ref{lemmad1} has two
implications, once we take $O_0$ small enough. First, $O_A$ lies in a
neighborhood of the ASD connections on which Taubes'
contracting-mapping argument lets us ``project'' the image of $A$ to
an ASD connection.  Second, by the implicit function theorem, the
image of $O_A$ in $\M_{k+1}$ is an 8-dimensional submanifold of
$\M_{k+1}$.

The quantity $\tx-\pi_A\tx$ will be central to the definition of the
remainder terms in $\muloc(p)$ and to the analysis in \S
\ref{sectpfthmrem2}.
We define 
\bearray\label{defxix}\non
\xi_X &=& \tx-\pi_A\tx \\
&=& d^A G^A_0 (d^A)^*\tX \ +\  (d^A_+)^*G^A_+ d^A_+\tX
\eearray
Here $G^A_0$ and $G^A_+$ are the inverses of the Laplacians $(d^A)^*d^A$
and $d^A(d^A_+)^*$ on $\Omega^0(\adp)$ and $\Omega^2_+(\adp)$, respectively. 

    We make three observations here.  First, if $A$ is ASD,
$\|(d^A_+)^*G^A_+d^A_+\tx\|_2/\|\tx\|_2$ is small for {\em any}
rotation vector field, not just SD ones.  This is to be expected since
any rotation vector field is an approximate isometry and hence should
approximately preserve anti-self-duality. However
$\|d^AG^A_0(d^A)^*\tx\|_2/\|\tx\|_2$ is small only for the rotation
vector fields of duality opposite to that of the connection.  Second,
to deduce from this smallness that the parameter space injects
(locally) into $\B$, one must know that the first eigenvalue of the
Laplacians on zero-forms does not tend to 0 as $\l$ tends to zero, as
it will if the ``background'' connection is flat (or merely
reducible).  Indeed on $\M_1(S^4)$, all rotation vector fields (not
cut off), lifted as above, preserve the standard instanton. (On
$\M_1(\bfr^4)$, if one writes instantons in the usual gauge and
instead lifts rotations using the flat connection, then ASD rotations
preserve the standard ASD instantons centered at the origin, while SD
vector fields induce the effect of a gauge transformation.)  Third,
because of the cutoff $\b$, $\mathfrak{h}_A$ is not a Lie subalgebra
of $\Gamma(TN)$, although in some sense it is close to being
one. Thus, while intuitively $\mathfrak{h}_A$ is associated with the
Lie algebra of an 8-dimensional group of translations, dilations, and
rotations, $O_A$ is not quite the orbit of an 8-dimensional local Lie
group, hence the term ``almost-action''.

We will return to this point at the end of this section, but first we
wish to relate $O_A$ to the gluing construction in [DK].  The
fibration of a region in $\M'_{k+1,\l_0}$ over $\M'_{k}$ is usually
viewed in terms of center point, scale, and gluing parameter.  We
claim that on an infinitesimal level these are essentially the eight
parameters used to define the approximate tangent space.  Indeed [GP2,
\S 5] it was shown that lifts using the translation and dilation
vector fields do correspond to infinitesimal changes in center point
and scale, up to an error that is essentially $O(\l)$.  ([GP2] dealt
only with $\M_1$, but under a suitable definition of ``concentrated''
the same argument works more generally). It remains to identify our
action of $SO(3)$ (the SD rotations) with the ``gluing parameters'' of
the construction in [DK,
\S 7.2].  As both constructions are non-canonical we content
ourselves with a somewhat heuristic correspondence.

Instantons in the subspace
$\M^{\kappa,\nu}_{k+1,\l_0}$ have a single ``charge-one bubble'' and
are otherwise not concentrated.  For any such ASD reference connection
$A_0=A$ there exists a gauge over ball $B(p_A,K\l)$ such that after
pulling back to $B(0,K\l)\subset\bfr^4$ by a positively-oriented
normal coordinate system $\{ x^i\}$, the connection form is close to
$A^{\rm std}_\l$, the standard instanton on $\bfr^4$ of scale $\l$ and
center the origin.  (See [DK], \S 8.2.1.) Here $K>1$ is any
fixed number and ``close'' means that after dilating by $\l$, the two
connections are $C^2$-close on $B(0,K)\subset\bfr^4$; the undilated
connections satisfy $|(\na^j(A-A^{\rm std}_\l)| \leq \e_1\l^{-1-j}, \
0\leq j\leq 2$, where by taking $\l_0$ small enough we can take $\e_1$
as small as we please.  After a choice of normal coordinate system,
the identification $\bfr^4\iso{\bf H}$, and an identification of
$SU(2)$ with the unit quaternions, the connection form for $A^{\rm
std}_\l$ on our ball is
\be\label{astd}
\w_0=\frac{{\rm Im}(\overline{x} dx)}{\l^2+|x|^2}.
\ee

For integers $j=1,\dots,10$ define the sets $U_j=B(p_A,jn_0\l^{1/2})$
and $V_j=N-U_j$.  Also let let $U_\infty$ denote the annulus
$U_{10}\intersect V_1$ and let $\Omega$ denote the smaller annulus
$U_9\intersect V_2$.  We choose gauges $s_0,s_\infty$ (local sections
of $P$) over $U_{10},U_\infty$ respectively such that the transition
function between $U_{10}$ and $U_\infty$ is
$g_{0\infty}(x):=\overline{x}/|x|$ (i.e.
$s_\infty=s_{0}g_{0\infty}$); furthermore we take $s_{0}$ to be the
radial gauge for $A$ with respect to $p_A$ in which (\ref{astd})
is written.  The connection form for $A^{\rm std}_\l$ with respect to
$s_\infty$ on $U_\infty$ is then
\be
\w_\infty=\frac{{\l^2\rm Im}(x\overline{dx})}{|x|^2(\l^2+|x|^2)}.
\ee

To make contact with the construction in [DK], we will pretend that on
the ball $U_8$, our connection $A$ is exactly standard (so that the
connection form relative to $s_0$ on $U_8$ is (\ref{astd})). Let
$\tilde{\b}$ be a function that is identically 1 on $N-U_\infty$ and
identically 0 on $\Omega$, with $|\na\tilde{\b}|\leq
cn_0^{-1}\l^{-1/2}$.  (Note that the ``interior'' part of the support of
$\tilde{\b}$ occurs where $\b\ident 1$.)  On $U_\infty$ let
$\w'_\infty$ denote the connection form $\tilde{\b}\w'$.  We then
define a new connection $A'$ on $P$ by declaring the connection form
for $A'$ in the gauge $s_\infty$ over $U_\infty$ to be $\w'_\infty$,
and declaring $A'= A \ {\rm on}\ N-U_\infty$.  We think of $A'|_{V_8}$
as a cut-off ``background connection''.  In fact, we can define a
bundle $P_k$ of Pontryagin index $k$ by replacing the transition
function $g_{0\infty}$ by the identity; $A'|_{V_8}$ extends to a
connection on $P_k$ that is flat near $p_A$.  Note that over $\Omega$
the connection $A'$ is flat; its connection form there, relative to
$s_0$, is $g_{0\infty}dg_{0\infty}^{-1} = {\rm
Im}(\overline{x}dx)/|x|^2.  $

Our choice of normal coordinates and identification $\bfr^4\iso{\bf
H}$ induces a Lie algebra isomorphism
$\theta:\Lambda^2_+T_{p_A}^*M\to\mathfrak{su}(2)={\rm Im}({\bf H})$,
mapping the standard basis of $\Lambda^2_+T^*_0\bfr^4$ to ${\bf \{
i,j,k\}}$. (Alternatively, $\theta^{-1}$ is given by mapping 
$v\in\mathfrak{su}(2)$ first to the vector field induced by
quaternionic left-multiplication on {\bf H}, then to the two-form
obtained by lowering an index using the metric.)  Let $v=\theta(\a)$
and assume $|v|$ is not too large.  We consider the canonical flow of
$X=\b\a_{ij}x^i\ddxj$ acting on the cut-off connection $A'$. After
integrating the flow for time 1, the action on the base is $x\mapsto
\phi(x)=h_1(r)x$, where $r=|x|$ and where $h_1(r)=\exp(\b v)$.  Let
$A'_v$ be the connection determined by this integrated canonical flow.

An alternative flow, the ``$s_0$-flat'' flow, is obtained by lifting
$X$ to $P|_{U_8}$ using the flat connection determined by $s_0$, and
extending this flow to the complement of $U_8$ by the identity (since
$X$ is supported in $U_8$).  If we integrate the $s_0$-flat flow for
time 1, the resulting connection form  on $\Omega$
with respect to $s_0$ is
\be\label{good}
{\w'}_v=
\frac{{\rm
Im}(\overline{h_1(r)x}d(h_1(r)x))}{|x|^2}
=g_{0\infty}h_1^{-1}dh_1g_{0\infty}^{-1}+g_{0\infty}dg_{0\infty}^{-1} .
\ee
By our earlier comments, the connections resulting from the canonical
flow and the $s_0$-flat flow are gauge equivalent (and in fact are
equal outside $U_8$).  Thus $A'_v$ is gauge-equivalent to a connection
$A''_v$ equaling $A'$ on $N-U_\infty$, and whose connection form in
$\Omega$ (in the gauge $s_0$) is (\ref{good}).

We claim that the connection $A''_v$ is the one constructed in [DK,
p. 296].  The latter essentially begins with the connection $A'$
(thought of as a cut-off connection on $P_k|_{V_1}$ glued to a
connection on a $k=1$-bundle $P_1|_{U_{10}}$) and modifies it on
$U_\infty$ as follows. Let $h_1(r)$ be as above, let
$h_2(r)=\exp(-(1-\b) v)$ and consider the two gauge transformations
$\tilde{h}_1, \tilde{h}_2$ over $U_\infty$ given by
$\tilde{h}_i(s_\infty(x))=s_\infty(x) h_i(r)$.  The gauge
transformation $\tilde{h}_1$ does not extend to all of $P$ (unless
$\exp(v)=\pm1$), but it does extend to the bundle $P_k$ defined
earlier, and for $r\leq 4n_0\l^{1/2}$ changes the trivialization
$s_\infty$ (extended to $P_k$) by the constant $\exp(v)$.  Similarly
the gauge transformation $\tilde{h_2}$ extends to $P|_{U_{10}}$,
changing the trivialization $s_\infty$ for $r\geq 8n_0\l^{1/2}$ by
$\exp(-v)$.  Because $h_1^{-1}dh_1=h_2^{-1}dh_2$, the two gauge
transformations have the same effect on the flat connection
$A'|_{\Omega}$.  Therefore we can define a new connection $A^{\rm
DK}_v$ by
\be
A^{\rm DK}_v=\left\{ \begin{array}{c} \tilde{h}_2(A') \ {\rm on}\ U_9 \\
\tilde{h}_1(A') \ {\rm on}\ V_2
\end{array}\right.
\ee

The connection form for $A^{\rm DK}_v$ with respect to $s_\infty$ on
$\Omega$ is $h_1^{-1}dh_1=h_2^{-1}dh_2$, so with respect to $s_0$ the
connection form is precisely (\ref{good}).  Thus $A''_v$ and $A^{\rm
DK}_v$ coincide on $\Omega$.  Since $\tilde{h}_1\ident 1$ on $V_9$ we
have $A^{\rm DK}_v=A'$ on this region, and since $X\ident 0$ on $V_9$
we have $A''_v=A'$ there as well.  Thus $A''_v=A^{\rm DK}$ on $V_2$.
It remains to consider only $U_2$. On this ball, a computation shows
that the connection form for $A'$ relative to $s_0$ is
\be\label{cutw}
\w'_0= \frac{{\rm Im}(\overline{x}dx)}{|x|^2}\left(
1-\tilde{\b}(r)\frac{\l^2}{\l^2+|x|^2}\right).
\ee
On $U_2$, we have $h_2\ident 1$, so the connection form for $A^{\rm
DK}_v$ remains $\w'_0$.  But the connection $A'$ is also preserved by
the ``$s_0$-flat'' flow of $X$; replacing $x$ by $h_1(r)x$ in
(\ref{cutw}) does not change $\w_0$, since $h_1$ is constant on
$U_2$.
Therefore $A''_v=A^{\rm DK}_v$ over all of $N$. 

   Now let $A_v$ be the connection obtained from applying the
canonical flow of $X$ for time 1 to $A$ (rather than to $A'$).  The
preceding shows that up to gauge equivalence, when $|v|$ is not too
large the only differences between $A_v$ and $A^{\rm DK}_v$ arise from
the facts that (i) $A$ is only approximately standard on a small ball
$B(p_A,K\l)$ rather than exactly standard on the larger ball
$B(p_A,10n_0\l^{1/2})$, and (ii) we do not cut off $A_v$ before applying
the flow.  

   It should also be noted that since the subspace $\mathfrak{h}^{\rm
rot}_A$ corresponding to the SD rotation vector fields is not closed
under Lie bracket, if we let $\mathfrak{h}^{\rm rot}_A$ act on $A$ by
the canonical flow for time 1, we should not expect to get a closed
``orbit''.  But the construction in [DK] shows that the space of
gluing parameters is a copy of $SO(3)$.

To address this discrepancy, first note that if $N$ were $\bfr^4$ we could
dispense with the cutoffs in the definition of $\mathfrak{h}_A$.  The vector
fields would be globally defined, and would generate a Lie algebra
exponentiating to the group of motions of $\bfr^4$, 
\be \{ x\mapsto ax+b \mid
(a,b)\in {\bf H}^*\times {\bf H} \}.  
\ee 
The stabilizer of the origin would be ${\bf H}^*\times\{0\}$, and if
the initial connection $A$ were standard, the set of connections
generated by letting $SU(2) \subset {\bf H}^*$ act via the canonical
flow would give two copies of the space obtained by the construction
in [DK], as $(a,0)$ and $(-a,0)$ yield the same connection.
(Alternatively, if $v\neq 0$ is small enough, the connections
$A'=A^{\rm DK}$ and $A''_v=A^{\rm DK}_v$ are gauge-inequivalent,
because the gauge transformation $\tilde{h}_1$ defined earlier---which
always extends to $P|_{N-\{p_A\}}$---extends to $P$ if and only if
$\exp(v)=\pm 1$.)  From our earlier discussion the action of $(a,b)$
on the standard instanton is given by pulling back the connection form
by the {\em inverse} of $(a,b)$, which results in a connection of
scale $|a|$ and center $b$ (cf. (\ref{whyb})). The unit quaternion
$a/|a|$ corresponds to gluing angle, doubly parametrized.

Intuitively, then, we have the following picture.  Fix a reference
connection $A=A_0\in\M^{\kappa,\nu}_{k+1,\l_0}$. Let
$B\subset\mathfrak{su}(2)$ be the ball centered at the origin that is
carried diffeomorphically to $SU(2)-\{ -1\}$ by the exponential map,
and let $B'\subset\mathfrak{h}^{\rm rot}_A$ be the corresponding set
of SD rotation vector fields.  If we let the canonical flow of
elements in $B'$ act for time 1 on $[A_0]$ we obtain a space that (for
purposes of integrating reasonably behaved differential forms)
approximates two copies of the fiber $Z_{\rm DK}$.  This
correspondence becomes sharper as $\l_0\to 0$: as we take the limit
and rescale the (local) metric and normal coordinates correspondingly,
the failure of $\mathfrak{h}_A$ to close under Lie bracket disappears
on any ball of fixed rescaled size.  Furthermore, because the rescaled
metric becomes flat, the limiting space of vector fields
$\mathfrak{h}_A$ is the same whether we center the rotations and
dilations at $p_A$ or at $p$.  Thus the limiting action of ${\bf
H}^*\times {\bf H}$ above appears to generate an immersed manifold
that we can treat ``homologically'' as two copies of $Z_{\rm DK}$.
This discussion motivates the assumptions we make on $Z$ in the next
section.

\section{The fiber Z.} \label{sectfiber}
\setcounter{equation}{0}

For our purposes we need only consider one fiber $Z=Z_{\l_0}$ of the
projection $\M'_{k+1,\l_0}\to \M'_k$; we do not need to construct the
whole fibration.  We will assume $Z$ has the following five
properties.  The first three are known to be satisfied by $Z_{\rm
DK}$, so the key assumptions are really the last two, which
require the tangent spaces of $Z$ and $Z_{\rm DK}$ to be
close in various norms.  The assumptions can almost certainly be
weakened from those below, at the cost of considerably more technical
work.

({\bf Z1}) $Z$ fibers over $N$ via the projection $Z\to N$ sending a
concentrated connection to its center.  Given $U\subset N$ we let
$Z|_{U}$ denote the inverse image of $U$ under the projection.  We
assume that $N$ can be covered by a finite number of normal coordinate
charts $U_i$ (which we may take to be geodesic balls) such that for
each $i$ there is a two-to-one fiber-preserving covering map
$\overline{\rho}_i:{\bf H}^*_{\l_0}\times U_i\to Z|_{U_i}$ having
additional properties listed below.  Here ${\bf H}^*_{\l_0}=\{ a\in {\bf
H}^* \mid |a|<\l_0 \} \iso (0,\l_0)\times SU(2)$, where the
isomorphism is $a\mapsto (|a|,a/|a|)$.  (Note that the center-point
and scale maps are defined globally on $Z$; it is only for the purpose
of handling gluing parameters that we need to chop up $N$.)

In general a normal coordinate system $\{ x^j \}$ on $U_i$ determines
an identification between $U_i$ and a ball in ${\bf H}$ centered at
the origin, and hence a local action of ${\bf H}^*_{\l_0}\times {\bf H}$
on $U_i$ given by $((a,b),x)\mapsto ax+b$. We assume that on each
$U_i$ there is a positively oriented normal coordinate system $\{ x^j
\}$ on $U_i$ such that $\overline{\rho}_i$ is approximately given by
the induced canonical flow of this ${\bf H}^*_{\l_0}\times {\bf
H}$-action, based at the standard instanton on $\bfr^4$, in the sense
that ({\bf Z2})-({\bf Z5}) below are true.  From [D, \S 3], the
orientation induced on the fiber $Z$ by the standard orientation of
${\bf H}\times {\bf H}$ as a complex vector space is then compatible
with the standard orientations of $\M_{k+1}$ and $\M_{k}$ (i.e. the
orientation of $\M'_{k+1,\l_0}$ is the product of the orientation of
$Z$ and the pullback of the orientation of $\M_k$).  These are the
orientations used in (\ref{orientations}).

({\bf Z2}) We assume that for each $i$, the scale and center point of
$A=\overline{\rho}_i(a,b)$ are $\l(A)=|a|$ and $p_A=b$ (in quaternionic normal
coordinates) respectively. 

({\bf Z3}) Given $i$, let $[A_{a,b}]=\overline{\rho}_i(a,b)$ and let
$F_{a,b}=F_{A_{a,b}}$ . We assume that
for any $K>0$, on the ball $B(p_A,K\l(A))$ we
have
\be\label{bfz3}
\left| |F_{a,b}|-
\frac{\sqrt{48}|a|^2}{(|a|^2+|x-b|^2)^2} \right| \leq \e_1(\l)\l^{-2}
\ee
where $\e_1(\l)\to 0$ as $\l\to 0$.  

({\bf Z4}) Let $B$ be the component of $\exp^{-1}(SU(2)-\{-1\})$
containing $0$.  A tangent vector $v\in T_{(a,b)} ({\bf
H}^*_{\l_0}\times B)$ gives rise to a vector field on a neighborhood of
$b\in B$ that determines an element $X_v\in \mathfrak{h}_A$.  Writing
$\tx_v=\i_{X_v}F_A$, we require that
\be\label{z4a}
\| \overline{\rho}_*v+\tx_v \|_{L^2} \leq \e_2(\l)|v|
\ee
where $\e_2(\l)\to 0$ as $\l\to 0$.  Observe that because of
(\ref{lemmad1a}), we can alternatively write (\ref{z4a})
as $\| \overline{\rho}_*v-(-\pi_A\tx_v) \|_{L^2} \leq \e_2(\l)|v|$;
cf. (\ref{drho}). 

({\bf Z5}) Letting $\xi'_v =\overline{\rho}_*v+\tx_v$ and
$\xi_v=\tx_v-\pi^A\tx_v$, we further require that $\xi'_v$ satisfy the
same weighted $L^4$ bounds as $\xi_v$ given in Proposition
\ref{pigcor11} ((\ref{pigcor11a}) and (\ref{pigcor11c})), and the
pointwise bound (\ref{dogprop2c}).

\bs
If not for ({\bf Z4}) and ({\bf Z5}), we would not need to assume
({\bf Z1})-({\bf Z3}). By itself, ({\bf Z1}) follows from the
description of the ends of moduli space in [DK, \S\S 7.2 and 8.2]; we
simply take the local diffeomorphism $(0,\l_0)\times SO(3)\times
U_i\iso Z|_{U_i}$, and pre-compose with the covering map $SU(2)\to
SO(3)$.  Similarly, ({\bf Z2}) and ({\bf Z3}) follow from [DK, \S 8.2.1].

What is not clear is whether the construction in [DK] yields a fiber
whose tangent space at $[A]$ is sufficiently close to $\pi_A({\cal
H}_A)$ in the norms required for our analysis.  If the subspace
$\mathfrak{h}_A\subset \Gamma(TN)$ were a Lie subalgebra, then by
(\ref{drho}) the canonical flow would generate a fiber whose
tangent space at $[A]$ would be precisely $\pi_A({\cal H}_A)$.
However, $\mathfrak{h}_A$ is not closed under Lie bracket, and
the canonical flow of vector fields in $\mathfrak{h}_{A_0}$, acting on
a single reference connection $[A_0]$, has no chance of generating an
orbit that reasonably approximates {\it all} of $Z_{[DK]}$; the cutoff in
the translation vector fields prevents the canonical flow from moving
the center point very far from $p_{A_0}$, whereas all points in $N$
can occur as center in $Z_{[DK]}$.  But the estimates relevant to
proving Theorem \ref{thmst} are much less sensitive to changing the
definition of of translations than to changing the definition of
rotations and dilations, so it seems plausible that, by a patching
argument, altering the definitions of only the translation vector
fields in any significant way, we can splice together canonical flows
based at connections with nearby center points.  Presumably by
splicing enough flows together we can obtain a fiber that is
$C^1$-close globally to $Z_{[DK]}$ and $C^1$-close locally to the
orbit of some canonical flow.  Even if the splicing construction
fails, there are two reasons why, for
purposes of integration, we may not need to define a true fiber (such as
$Z_{[DK]}$) in a topological sense.  First, when we integrate
$\mu(p)\wedge\mu(q)$ over an orbit of the canonical flow, only
connections with center point near $p$ and $q$ contribute
significantly to the integral.  It is likely that the same holds
for an integral over $Z_{[DK]}$, so that it suffices to approximate
only a region of $Z_{[DK]}$ consisting of connections with center
point in a fixed small ball.  Second, although the canonical flow of
the subspace $\mathfrak{h}_{A_0}^{\rm rot}$ acting on $[A_0]$ does not
generate a {\em closed} manifold, it does generate an immersed copy of
$SU(2)-\{-1\}$ lying in a small neighborhood of an $SO(3)$-orbit in
$Z_{[DK]}$, and which geometrically wraps twice around this orbit.  A
careful analysis may show that there is a homotopy from the immersed
punctured $SU(2)$ to a punctured double cover of the $SO(3)$ in
$Z_{[DK]}$, small enough in all relevant norms that
there is only a negligible difference between integrating over
$Z_{[DK]}$ and over the orbit of the canonical flow.

   Thus the idea behind ({\bf Z1}--{\bf Z5}) is basically that there
is fiber that interpolates between $Z_{[DK]}$ and the not-quite-fiber
generated by splicing together canonical flows.  The hypotheses ({\bf
Z4}--{\bf Z5}) amount to assuming that, in this interpolated fiber, the
bounds on $\xi'$ are as good as they would be if the tangent space to
the fiber were the one determined by the canonical flow.  We need such
an assumption because, when we pull $\mu_d(p)\wedge\mu_d(q)$ back to $Z$,
we need to insert true tangent vectors to $Z$ into (\ref{mud}); the
$\pi_A\tx$'s in the expansion (\ref{expansion0}) below should be
replaced by $\overline{\rho}_*v$'s---which has the effect of replacing
each $\xi$ in (\ref{defrem2}) with $\xi'$.

There is other evidence making the simultaneous satisfaction of at
least ({\bf Z1}--{\bf Z4}) very plausible.  On ${\bfr^4}$, if we
remove the cutoffs in the definition of $\frak{h}_A$ and define $Z$
from the canonical flow acting on the standard instanton, then the
spaces $T_{[A]}Z$ and ${\cal H}_A$ coincide.  In the case of
1-instantons over simply connected definite manifolds (where the
background connection is flat and there are no gluing parameters, so
$Z$ is five-dimensional), ({\bf Z4}) was shown in [GP2] to be true
with $\e_2(\l)\leq c\l^{1-\d}$ for small $\d>0$; in
\cite{Gppt} this was strengthened to $\l^{1+\d}$.

The technical hypothesis ({\bf Z5}) is more ad hoc, and stronger than
necessary, but is not without basis. In the setting of the
five-dimensional moduli spaces mentioned above, certain estimates of
this paper and \cite{Gppt} can be combined to show that
$\|r_p^{-\d}\xi'\|_4\lessim
\l^{-1+\d}(b\cdot \l+a\cdot \l^{1/2})$ and $\| r_A\xi'\|_4
\leq b\cdot \l^\d(b\cdot \l+a\cdot \l^{1/2})$, much stronger 
than the $L^p$ bounds assumed in ({\bf Z5}). (Here $r_p$ denotes
distance to an arbitrary point $p\in N$.)

    An important implication of ({\bf Z4}) and (\ref{lemmad1a}) is the
following.  Let $dvol_{Z}^{L^2}$ be the volume form on $Z$ induced by
the $L^2$ metric on $\B$.  (Hypothesis ({\bf Z1}) determines an
orientation on $Z$, so there is no sign ambiguity here.)  Let $a\in
\bfr^4$ denote the quaternionic variable in ${\bf H}^*_{\l_0}$.  Then
\be\label{comparevols}
\overline{\rho}_i^*(dvol_Z^{L^2}) \approx {\rm const}\cdot d^4 a \wedge dvol_N
= {\rm const}\cdot \l^3d\l\wedge \ dvol_{S^3} \wedge dvol_N.
\ee
where the approximation becomes exact as $\l\to 0$ (and the constant
is of course nonzero).  Our chief use of (\ref{comparevols}) will be
to help estimate the integrals of the non-local terms in
$\mu_d(p)\wedge\mu_d(q)$.  For this purpose, we don't actually need
``$\approx$'' in (\ref{comparevols}); ``$\leq$'' would suffice.  Thus
hypothesis ({\bf Z4}) can be weakened.



\section{Localizing $\mu_d(p)\wedge\mu_d(q)$.}\label{sectlocal}
\setcounter{equation}{0}

    From now on we assume there is a fiber $Z$ with properties ({\bf
Z1}--{\bf Z5}).  To motivate the leading-term calculation in
(\ref{stfails}), suppose for the moment that for $[A]\in Z$ the
tangent space $T_{[A]}Z$ is precisely, rather than approximately, the
space $\pi_A({\cal H}_A)\subset T_{[A]}\M$.  Then if we pull $\mu_d(p)$
back to $Z$, we need only apply (\ref{mud}) to arguments of the form
$\pi_A\tx_A$ with $X\in\mathfrak{h}_A$.  Recalling the definition of
$\xi_X$ in (\ref{defxix}), we then have
\be\label{expansion0}
G^A_0\{\pi^A\tx,\pi^A\ty\}=G^A_0 \left( \{ \tx,\ty \} +Rem_1(X,Y;A) \right),
\ee
where
\be\label{defrem1}
Rem_1(X,Y)=\{\tx,\xi_Y\}-\{\ty,\xi_X\}+\{\xi_X,\xi_Y\}.
\ee
(We omit writing most of the $A$-dependence in these
formulas explicitly.) Here $\{\cdot,\cdot\}$ is a universal, local,
antisymmetric bilinear pairing that takes two $\adp$-valued 1-forms
and produces an $\adp$-valued zero-form.  Note that $Rem_1(X,Y)$ is
antisymmetric in $X$ and $Y$.

    In [G1, Proposition 2.1], it was shown how to expand several of
the expressions appearing in (\ref{expansion0}--\ref{defrem1}) as a
leading-order local term plus a non-local remainder, smaller in
appropriate norms.  In particular, for any vector fields $X,Y$ on $N$,
we have
\be\label{expansion2}
G^A_0\{ \tx,\ty \} =
-\frac{1}{2}F(X,Y)+G^A_0(R''(X,Y)),
\ee
where $F=F_A$ and where
\bearray
\non
2R''(X,Y)
&=& {\cal R}(F )(X,Y)+ F (\lap
X,Y)+F (X,\lap Y) \\ 
\label{defr2prime}
&-&2(\na_i^AF )(\na_i X,Y) -2(\na_i^AF )(X,\na_i Y) -2F (\na_i X,\na_i Y).
\eearray
Here $\cal R$ is an endomorphism proportional to the Riemann tensor
whose precise form does not concern us.  As a consequence
of (\ref{expansion2}), $R''(X,Y)$ is antisymmetric in $X$ and $Y$.
The precise way in which the derivatives of $F$ and the derivatives of
$X$ and $Y$ are hooked together in (\ref{defr2prime}) is critical for
certain estimates (Lemma \ref{oldlemma39}).

   Applying (\ref{expansion2}) to the first term in
(\ref{expansion0}), we find
\be
G^A_0\{\pi \tx,\pi \ty\}=-\frac{1}{2}\left(F(X,Y)-Rem_2(X,Y)\right),
\ee
where
\bearray\non
\frac{1}{2}Rem_2(X,Y)&=& G^A_0(R''(X,Y) +Rem_1(X,Y))\\
&=&G^A_0\left( R''(X,Y) +  \non
\{\tx,\xi_Y\}-\{\ty,\xi_X\}+\{\xi_X,\xi_Y\}\right)\\
&:=&G^A_0( Rem_2'(X,Y)).
\label{defrem2}
\eearray
Inserting all this into (\ref{mud}) we find
\bearray
\non
\lefteqn{\mu_d(p) (\pi \tx,\pi \ty, \pi\tilde{V},\pi\tw)=} \\
\non
&=&\frac{1}{4\pi^2} \tr \left(
\left. \left(F(X,Y)F(V,W)+F(X,V)F(W,Y)+F(X,W)F(Y,V)\right)\right|_p \right)  \\
 && + \left. Rem_3(X,Y,V,W)\right|_p + \left. Rem_4(X,Y,V,W) \right|_p,
\label{mumap10}
\eearray
where
\be\label{defrem3}\label{defrem4}
Rem_3= {\rm const} \cdot
\tr(F\wedge Rem_2), \ \ 
Rem_4 = {\rm const}\cdot 
\tr(Rem_2 \wedge Rem_2).
\ee
(In (\ref{defrem3}) we regard $F$ and $Rem_2$ as $\Gamma(\adp)$-valued
two-forms on the space of vector fields.)
\ss
\ni 
The first term in (\ref{mumap10}) is just
$(8\pi^2)^{-1}\tr(F\wedge F)(X,Y,V,W)|_p$, which, since $F$ is ASD,
can be rewritten as $(8\pi^2)^{-1}|F|^2 dvol(X,Y,V,W)|_p$. Thus
if we define
\be\label{defmuloc}
\muloc(p)(\pi \tx,\pi \ty, \pi\tilde{V},\pi\tw) 
=\frac{1}{8\pi^2} |F(p)|^2 dvol(X,Y,V,W)|_p,
\ee
then (\ref{mumap10}) simplifies to
\bearray\non
\mu_d(p) (\pi \tx,\pi \ty, \pi\tilde{V},\pi\tw) 
&=&\muloc(p) (\pi \tx,\pi \ty, \pi\tilde{V},\pi\tw) \\
&&+\non
Rem_3(X,Y,V,W)|_p + Rem_4(X,Y,V,W)|_p.
\\ \label{mainapprox}
\eearray

Of greatest concern to us will be the local part $\mu_{\rm loc}(p)$ of
this expression.  Note that $\mu_{\rm loc}(p)\wedge \mu_{\rm
loc}(p)=0$, since $dvol_p\wedge dvol_p =0$.  However, we will see that
$\lim_{q\to p}\int_Z
\mu_{\rm loc}(p)\wedge\mu_{\rm loc}(q) \neq 0$.  In this integral it
turns out that instantons of scale $\approx {\rm dist}(p,q)$ give
the main contribution to the integral.  Thus the pullback of 
$\mu_d(p)\wedge\mu_d(p)$ to $Z$ can be thought of loosely as a $\d$-form
concentrated on instantons of scale zero.  

   To integrate $\mu_d(p)\wedge\mu_d(q)$ we must still worry about the
non-local remainder terms $Rem_i$, as well as the fact that the
tangent space $T_{[A]}Z$ is not precisely $\pi_A{\cal H}_A$. We will
see later that, as $\l_0\to 0$, the contributions to the integral of
$\mu_d(p)\wedge\mu_d(q)$ over $Z=Z_{\l_0}$ from both of these corrections
tend to zero.  What we wish to compute now is
\be\label{idea}
\lim_{q\to p} \int_{Z_{\l_0}} \muloc(p)\wedge\muloc(q),
\ee
where $p$ and $\l_0$ are fixed. 

For given $p,q$, as we integrate $\muloc(p)\wedge\muloc(q)$ over $Z$,
the center point $p_A$ of $[A]$ in $Z$ moves around, affecting the
support of the vector fields $X,Y,V,W$ in (\ref{defmuloc}).  Thus for
$\muloc(p)\wedge\muloc(q)(\pi\tilde{X}_1,\dots,\pi\tilde{X}_8)$ to be
nonzero, $p_A$ must lie in $ B(p,8n_0\l^{1/2})\intersect
B(q,8n_0\l^{1/2})$. In particular we can restrict $p_A$ to a small
normal-coordinate ball $U$ centered at $p$ (which we can take to be
one of the $U_i$ in ({\bf Z1})) without affecting
$\int_Z\muloc(p)\wedge\muloc(q)$.  Since we are interested in the
limit as $q\to p$, we may also assume $q\in U$.

    Let $2L={\rm dist}(p,q)$; we will later send $L$ to zero.     
Define $Z_1\subset Z$ to be the set of instantons in $Z$ obeying
the two criteria
\be\label{z1scale}
L^{0.1} \geq \l^{1/2} \geq L, 
\ee
\be
\label{z1center}
p_A \in  B(p,n_0\l^{1/2}).
\ee
Note that if $[A]\in Z_1$ then $p_A\in B(q, (n_0+1)\l^{1/2})$, so that
the cutoff $\b$ in the definition of the vector fields $X_i$ equals
1 at both $p$ and $q$.  We will see later that the contribution to
(\ref{idea}) from the complement of $Z_1$ is negligible.

Let $\{x^i_{\rm old}\}$ denote normal coordinates on $U$. We change
coordinates by setting $x_{\rm new}=L^{-1} x_{\rm old}$ and
replace the metric $g_{\rm old}$ on $U$ by $g_{\rm new}
=L^{-2}g_{\rm old}$.  Because of the conformal invariance of
$|F|^2 dvol$, $\muloc(p)\wedge\muloc(q)$ is unaffected by this change.
However, since $\l=\l_{\rm old}$ represented a distance in the old
coordinate system, we now have a rescaled upper cutoff for $\l_{\rm
new}=L^{-1}\l_{\rm old}$ on $Z$, namely $\l_{0,{\rm
new}}=\l_0/L$.  Measuring all distances in the new metric, 
the defining conditions for $Z_1$ become
\be\label{z1c}
L^{-0.8} \geq \l_{\rm new} \geq L, 
\ee
\be
\label{z1d}
p_A \in  B(p,NL^{-1/2}\l_{\rm new}^{1/2}).
\ee
As $L\to 0$, several things happen.  For $A\in Z$, $|F_A|$ becomes
approximately standard on any fixed ball $B(p,K)$; $g_{\rm new}$
approaches the flat metric $\sum (dx^i_{\rm new})^2$; and (in the
rescaled metric and coordinates), $Z_1$ becomes an $SO(3)$-bundle over
monotonically increasing regions of center-scale space that
exhaust $(0,\infty)\times(\bfr^4-\{ 0\})$ as $L\to 0$.  Because of ({\bf
Z1}), we can identify $Z_1$ with ever-increasing subsets $G_L$ of 
$G:=({\bf H}^*\times {\bf H})/{\bf Z}_2$.
Letting $\muloc'$ denote the pullback of $\muloc$ to ${\bf H}^*\times
{\bf H}$, we therefore have
\be\label{half}
\lim_{L\to 0} \int_{Z_1} \muloc(p)\wedge\muloc(q) =
\frac{1}{2} \lim_{L\to 0}
\int_{G_L} \muloc'(p)\wedge\muloc'(q),
\ee
provided this integral converges.  

Let $\overline{\rho}$ be as in ({\bf Z1}--{\bf Z5}).
Write elements of $G$ as pairs $(a,b)$, and write
$A_{(a,b)}=\rho(a,b)$, $F_{(a,b)}=F_{A_{(a,b)}}$ as in ({\bf Z3}). 
If we define
$\muloc'=\overline{\rho}^*\muloc \in \Omega^4(G)$, then
\be \int_{G_L} \muloc'(p)\wedge\muloc'(q) = 
\int_{G_L} \muloc'(p)\wedge\muloc'(q)(\frac{\partial}{\partial a^1},
\dots, \frac{\partial}{\partial b^4}) da^1 \wedge \dots \wedge da^4
\wedge db^1 \wedge \dots \wedge db^4.
\ee

To compute this we need to know ${\overline{\rho}_*}_{(a,b)}
\partial/\partial a^i, {\overline{\rho}_*}_{(a,b)}\partial/\partial
b^i$.  At each $(a,b)$ define ${X}_i,{Y}_i$ to be the vector
fields on $\bfr^4$ induced by $\partial/\partial a^i,
\partial/\partial b^i$ respectively.  Temporarily writing
$b^i=a^{i+4}$ and $Y_i=X_{i+4}$, 
from ({\bf Z4}) there is an 8x8 matrix $C=Id+O(\e_2(\l_{\rm old})))$ 
for which we have 
\be\label{lorenzowantsanumber}
{\overline{\rho}_*}_{(a,b)}{C^j}_i\partial/\partial a^j
=-\pi_{A_{(a,b)}}\i_{X_i}F_{(a,b)}.
\ee
Hence from (\ref{mainapprox}), if not for the correction matrix $C$ we
would have
\bearray\non
\muloc'(p)(\frac{\partial}{\partial a^1}, \dots,
\frac{\partial}{\partial a^4} ) &=&  
\muloc(p)
(\pi_{A_{(a,b)}}\tx_1^A,\dots,\pi_{A_{(a,b)}} \tx_4^A) \\
&=&
(8\pi^2)^{-2}|F_A(p)|^2dvol_p(X_1,\dots,X_4),
\eearray
with a similar formula if we replace any of the $\partial/\partial
a^i$'s by a $\partial/\partial b^i$.  

Let us ignore, for now, (i) the $O(\e_2(\l_{\rm old}))=O(\e_2(L\l_{\rm
new}))$ difference between the matrix $C$ and the identity, and (ii)
the $O(|x_{\rm old}|^2)=O(L^2|x_{\rm new}|^2)$ difference between the
true metric on the rescaled ball and the flat metric; we will make the
corrections later.  Since the Euclidean volume form is
$dvol=dx^1\wedge\cdots\wedge dx^4$, we will write $dvol_p=d^4x_p$,
$dvol_q=d^4x_q$ below. Hence
\be
\muloc'(p)\wedge\muloc'(q)(\frac{\partial}{\partial a^1}, \dots,
\frac{\partial}{\partial b^4} )
=(8\pi^2)^{-2} |F_A(p)|^2|F_A(q)|^2 
 d^4x_p \wedge d^4x_q (X_1, \dots, Y_4).
\ee

So far we have treated $d^4x_p\wedge d^4x_q$ as an 8-form whose
arguments are vector fields, but we may as well consider it as an
8-form on the 8-dimensional space $T_pM\plus T_qN$. Using the
canonical isomorphisms $T_p{\bfr^4}\iso T_q{\bfr^4}\iso {\bfr^4}$ and
our further identification of ${\bfr}^4$ with ${\bf H}$, we can write
each $X_i,Y_j$ in the form $(v,w)\in {\bf
H}\plus {\bf H}$.  In the coordinate system $\{x^i_{\rm new}\}$, the origin
represents $p$, and we may assume that $q$ lies on the real axis with
coordinate $2\in {\bf H}$.  Let
$\tau_1=1\in {\bf H}$ and let $\{\tau_i\}_2^4$ be the quaternions
${\bf i,j,k}$.  Then $X_i(x)=\tau_i a^{-1}(x-b)$ and $Y_i(x)=\tau_i$.
so the corresponding elements in $T_p\bfr^4\plus T_q\bfr^4
\iso {\bf H}\plus{\bf H}$ are $X'_i=(-\tau_i a^{-1}b, \tau_i a^{-1}(2-b))$
and $Y'_i= (\tau_i,\tau_i)$. Modulo the span of the $Y'_i$,
we have $X'_i=(0, 2\tau_i a^{-1}) := X''_i$, so $d^4x_p\wedge d^4x_q
(X'_1, \dots, X'_4, Y'_1, \dots Y'_4) = d^4x_p\wedge d^4x_q (X''_1,
\dots, X''_4, Y'_1, \dots Y'_4)$.  Since $d^4x_p(X_i'', *,*,*)=0$,  
it follows that $d^4x_p\wedge d^4x_q (X''_1, \dots, X''_4, Y'_1, \dots
Y'_4) =d^4 x_p(Y'_1, \dots, Y'_4) \ d^4x_q(X''_1,\dots, X''_4)$.  But
$d^4x_p(Y'_1,\dots,Y'_4)=1$ and $d^4x_q(X_1,\dots, X_4)= 2^4|a|^{-4}$.
Hence
\be
\muloc'(p)\wedge\muloc'(q)|_{(a,b)}(\frac{\partial}{\partial
a^1}, \dots, \frac{\partial}{\partial b^4})
= (8\pi^2)^{-2} |F_{(a,b)}(0)|^2|F_{(a,b)}(2)|^2 2^4|a|^{-4}.
\ee
where $F_{(a,b)}$ is the curvature of the instanton obtained from the
action of $(a,b)$ on a reference connection in our fiber.  Therefore
\be\label{almostthere}
\muloc'(p)\wedge\muloc'(q)|_{(a,b)}
=2^4(8\pi^2)^{-2} |F_{(a,b)}(0)|^2|F_{(a,b)}(2)|^2  |a|^{-4}
d^4a\wedge d^4b.
\ee

Because $\l_0$ is small, there is a reference connection in our fiber
that looks approximately standard on a ball of any fixed large radius,
with the approximation getting better as $\l_0\to 0$ (the rescaling by
$L$ only improves this approximation).  Our next approximation is to
ignore the difference between the true reference connection $A_0$ and
the standard instanton; we will deal with the error later.  The
connections in the limiting $Z_1$ are then the orbit of the standard
instanton $A_1$ under the action of $G$.  Hence
\bearray\non
|F_{(a,b)}|^2(x) &= &|\tilde{\Phi}_{(a,b)}^{-1}F_{A_1}|^2(x)\\
&=& \left|\frac{d(\overline{a^{-1}(x-b)})\wedge d (a^{-1}(x-b))}
{(1+|a^{-1}(x-b)|^2)^2}\right|^2 \non \\ 
&=& 
\frac{48|a|^4}{(|a|^2+|x-b|^2)^4}.
\label{whyb}
\eearray
Thus
\be \label{ipq}		
I_p:=\lim_{L\to 0}
\int_{G_L} \muloc'(p)\wedge\muloc'(q)
= 
36\pi^{-4}\int_{{\bf H}^*\times {\bf H}}
\frac{2^4\ |a|^{4}\  d^4 a \wedge d^4 b}{
(|a|^{2}+|b|^2)^4(|a|^{2}+|2-b|^2)^4},
\ee
and provided the error terms we have so far ignored are truly
ignorable,
\be
\lim_{L\to 0} 
\int_{Z} \muloc(p)\wedge\muloc(q)
 = \frac{1}{2}I_p
\ee
(see \ref{half}). 

\begin{lemma}
$I_p=1.$
\end{lemma}

\pf First
introduce spherical coordinates in $a$-space (with radial variable we
call $\l$) and cylindrical coordinates in $b$-space (with radial
variable $r$).  The integrals over the 3-sphere in $a$-space 
and the 2-sphere in the imaginary subspace of $b$-space are
trivial, contributing factors $2\pi^2$ and $4\pi$ respectively. Thus
\be
\label{int2}
I_p=
36\pi^{-4} \cdot 8 \pi^3 \cdot 2^4 \cdot
\int_{\l=0}^\infty \int_{r=0}^\infty
\int_{z=-\infty}^\infty 
\frac{\l^7 r^2 }{
(\l^{2}+r^2+ z^2)^4(\l^{2}+r^2+(z-2)^2)^4}dz\ dr\ d\l.
\ee
Introducing polar coordinates in the $\l$-$r$ quarter-plane, the
angular integration reduces us to an integral over two real variables.
Using the Residue Theorem to integrate over $z$ leaves us with a
one-dimensional integral that can be computed in closed form, yielding
$I_p=1$.
\qed

    In the local calculation we ignored errors from four sources: (i)
the contribution from the complement of $Z_1$, (ii) the difference
between the flat metric and the true metric on the rescaled ball; (iii)
the difference between $\overline{\rho}_*v$ and $-\pi_A\tx_v$
(i.e. the difference between the matrix $C$ and the identity); and
(iv) the difference between $|F|_{a,b}$ and the standard instanton of
scale $|a|$ and center $b$.

   Let us first deal with (i).  Since the vector fields $X_i$ we feed
into $\muloc$ are cut off at distances $\geq 2n_0\l^{1/2}$ from $p_A$,
the integrand $\muloc(p)\wedge\muloc(q)(X_1,\dots,X_8)$ vanishes for
$p_A$ outside the ball $B(p,2n_0\l^{1/2})$.  For purposes of integration
we therefore need only that portion of $Z$ lying over a ball of fixed
small radius centered at $p$.  Because of (\ref{ptwsfest2}), the
integrand over such a region is bounded by a constant times the
integrand we used in our previous calculation, cut off in certain
regions.  Since the integrand in (\ref{ipq}) is integrable over all of
${\bf H}^*\times{\bf H}$, given any exhaustion $W_1\subset W_2\subset
...$ of ${\bf H}\times{\bf H}$, the integral over the complement of
$W_n$ goes to zero as $n\to \infty$.  As the sets $G_L$ provide such
an exhaustion, the integral of $\muloc(p)\wedge\muloc(q)$ over the
complement of $Z_1$ tends to zero.

   Next we turn to the errors (ii)-(iv) listed above.  In place of the
set $Z_1$ considered in the derivation above, for $K>0$ consider the
sets $Z_{K,L}$ defined by $\{ L^{0.1}\geq \l^{1/2} \geq L, \ \ p_A \in
B(p,K\l)\}$.  After rescaling by $L$ as before these conditions become
$\{ \l_{\rm new} \geq L, \ p_A \in B(p,K\l_{\rm new})\}$.  This time
as $L\to 0$, $Z_{K,L}$ exhausts $Z_{K,0}:=\{ p_A \in B(p,K\l_{\rm
new})\}$, with $\l_{\rm new}$ unrestricted.  But convergence of the
integral $I_p$ implies that given $\e_3>0$ we can choose $K$ large
enough and $L$ small enough that the integral of the integrand in
(\ref{ipq}) over the complement of $Z_{K,L}$ is less than
$\e_3$.  On the interior set $Z_{K,L}$, hypotheses ({\bf
Z3}--{\bf Z4}) imply that given $\e_4>0$, by taking $L$ sufficiently
small we can arrange for the ratio of the true
$\muloc'(p)\wedge\muloc'(q)$ to be within a multiple $(1+\e_4)^{\pm
1}$ of the integrand in (\ref{ipq}) over all of $Z_{K,L}$.  (Error (ii)
gives an $O(\l_{\rm old})\leq O(L^{0.2}))$ contribution to $\e_4$;
error
(iii) a contribution $\e_2(\l_{\rm old})\leq \e_2(L^{0.2})$, through the
matrix $C$.  As for error (iv), in the rescaled metric and coordinates
({\bf Z3}) implies 
\be\label{bfz3_rescaled}
\left| |F_{a,b}|_{g_{\rm new}}-
\frac{\sqrt{48}|a|^2}{(|a|^2+|x-b|^2)^2} \right| \leq \e_1(\l_{\rm
old})\l_{\rm new}^{-2},
\ee
so that this error gives a 
contribution $\e_1(\l_{\rm old}) \leq \e_1(L^{0.2})$ to $\e_4$.)
Hence we
can arrange for the integral of
$\int_{Z_{K,L}}\muloc'(p)\wedge\muloc'(q)$ over lie within $\e_3/2$ of
the integral over $Z_{K,L}$ of the integrand in (\ref{ipq}). It
follows that by choosing $L$ small enough, the errors introduced by
our approximations can be made arbitrarily small.

We have now proven the following.

\begin{prop} \label{proplocal}For any $p\in N$, 
\be
\lim_{q\to p} \int_{Z_{\l_0}} \muloc(p)\wedge\muloc(q)
=\frac{1}{2}.
\ee
\end{prop}
\qed


\section{The nonlocal terms in $\mu_d(p)\wedge\mu_d(q)$.}\label{sectnonlocal}
\setcounter{equation}{0}

  From (\ref{defrem3}--\ref{mainapprox}), $\mu_d(p)\wedge\mu_d(q)$ can be
expanded as $\muloc(p)\wedge\muloc(q)$ plus a remainder.  Our next
task is to show that, as $\l_0\to 0$, the contribution of this
remainder to $\int_{Z} \mu_d(p)\wedge\mu_d(q)$ tends to zero.  This will
follow from the next proposition, whose proof occupies the remainder
of this paper.

\begin{prop}\label{propintrem}
Let $\Omega$ be the restriction to $Z_{\l_0}$ of $\Omega_{\M}:=
\mu_d(p)\wedge\mu_d(q)-\muloc(p)\wedge\muloc(q) \in \Omega^8(\M)$. Assuming
({\bf Z1}--{\bf Z5}), there exists $\d>0$ such that
\be\label{propintrema}
\int_{Z_{\l_0}} \Omega \leq {\rm const}\ \l_0^{\d},
\ee
where the constant is independent of $p$ and $q$. 
\end{prop}

Observe that  Propositions \ref{proplocal} and \ref{propintrem} together 
prove Theorem \ref{thmst}.

Proving Proposition \ref{propintrem} requires some bounds on
$Rem_2(X,Y)$ for $X,Y\in\frak{h}_A$.  Before starting to derive these,
we need some notational simplification. Below we will be computing
many things that are multilinear in data of the form $(\bfb,a,\a)\in
T_{p_A}N\plus\reals\plus\Lambda^2_+T_{p_A}N$.  Given a single vector
field $X$ constructed from such data, we can denote the defining data
of (\ref{defxhat}) by $(\bfb_X,a_X,\a_X)$.  This notation becomes
cumbersome, especially when computing objects that involve more than a
single vector field.  However, because $|X_{(\bfb,a,\a)}| \leq
c(|\bfb|+(|a|+|\a|)\l^{-1}r_A)$, the $a$ and $\a$ data always enter
our bounds with precisely the same weight, so for shorthand we will
generally lump the $a$ and $\a$ terms together, and simply call them
$a$.  Furthermore, for simplicity we will often omit the
subscripts $X,Y, \dots $ in the defining data $(\bfb_X,a_X,\a_X),
(\bfb_Y,a_Y,\a_Y) \dots$; the dependence on $X,Y,\dots$ can
be reconstructed from context.  E.g. if we write
\be
|\mbox{\rm something bilinear in $X,Y\in\frak{h}_A$}|\leq
c_1 b^2 + c_2 ba + c_3 a^2,
\ee
then on the RHS the notation has the following meaning:
\bearray
b^2&=& \non \bdotb, \\
ba&=& \non \baab, \\
a^2&=& \adota.
\eearray
If the bilinear quantity is antisymmetric in $X,Y$ (as in Theorem
\ref{thmrem2} below), then the
estimate factors through the wedge product, in which case we can take
\bearray
b^2&=& \non \bwedgeb, \\
ba&=& \non \baab, \\
a^2&=& \awedgea.
\eearray
   
Finally, the notation $x\lessim y$ means $x\leq cy$ for a constant $c$
that is uniform in all relevant parameters. 

\ss
With this notation in mind,  we have 

\begin{prop}\label{thmrem2} 
(a) There exists $\d>0$ such that
\be\label{thmrem2a}\label{boundrem2}
\lnorm Rem_2(X,Y) \rnorm_\infty \lesssim\l^{-1+\d}(b^2+ba\cdot \l^{-1/2}
+a^2\cdot\l^{-1/2}).
\ee
Furthermore, there exists $c_1>0$ such that, for $r_A\geq
c_1\l^{1/2}$, we have the pointwise decay
\be
| Rem_2(X,Y)|
\lesssim
r_A^{-1}\left( b^2 +ba\cdot
\l^{-1/2}+a^2\cdot\l^{-1/2} \right).
\label{thmrem2b}\label{boundrem2b}
\ee

(b) Let $v,X_v,\xi_v,\xi'_v$ be as in hypothesis ({\bf Z5}) of \S
\ref{sectfiber}. If we alter the definition of $Rem_2(X_v,X_w)$ by replacing
$\xi_v$ with $\xi'_v$, then the bounds above still apply.

\end{prop}

We will prove Proposition \ref{thmrem2} (actually a slightly stronger
version) in \S \ref{sectpfthmrem2}.  Let us assume it for now
and move onto its application, the proof of Proposition
\ref{propintrem}.  The decay estimate (\ref{thmrem2b}) is crucial in
this proof; the global bound (\ref{thmrem2a}) does not suffice.

\bs
{\bf Proof of Proposition \ref{propintrem}.}  By hypothesis ({\bf
Z1}), $\int_Z \Omega \leq \sum_i \rho_i^*\Omega$. Since
$\Omega\in\Omega^8(Z)$ we can write
$\Omega=f\ dvol_Z^{L^2}$, where the function $f$ can be computed at
$[A]\in Z$ from any (positively) oriented $L^2$-orthonormal basis
$\eta_1, \dots, \eta_8$ of $T_{[A]}Z$ by
\be
f([A])= \Omega(\eta_1, \dots, \eta_8).
\ee
Similarly, we define $f'([A])= \Omega_{\cal M}(\eta'_1, \dots, \eta'_8)$, 
where the $\{\eta'_i\}$ are an orthonormal basis for $\pi_A{\cal H}_A$,
and set $\Omega'=f'\ dvol_Z^{L^2} \in \Omega^8(Z)$. 

We will first show that $\int_Z \Omega' \leq c\l_0^\d$ (where $\d$ is
as in Proposition \ref{thmrem2}), and then deduce that the same is
true for $\int_Z \Omega$.

We proceed to estimate $\Omega'$.  By Proposition \ref{lemmad1}, an
approximately orthonormal basis of $\pi_A{\cal H}_A$, up to a scale
factor $(8\pi^2)^{1/2}$, is given by
$\{\eta'_n=\pi_A\tx_{(\bfb_n,a_n,\a_n)}:=\pi\tx_n\}_1^8$ as
$(\bfb_n,a_n,\a_n)$ run over an orthonormal basis of
$T_{p_A}N\plus\bfr\plus\Lambda^2_+(T_{p_A}N)$. Applying Gram-Schmidt
to $\{ \pi\tx_n \}$, it follows that $f'\leq {\rm const}\cdot
\Omega(\pi\tx_1, \dots, \pi\tx_8) dvol_Z^{L^2}$.  Hence from
(\ref{comparevols}),
\be\label{psiwprime}
\rho_i^*(\Omega')\leq c\ \Omega(\pi\tx_1, \dots, \pi\tx_8) \l^3d\l\wedge
dvol_{S^3}\wedge dvol_N.
\ee
Symbolically we can write $\Omega'=\sum_{i=0}^3\Omega'_i$ as a sum of
terms of the form $F^{i}\wedge Rem_2^{4-i}$, $0\leq i \leq 3$.  We
estimate the integrals of $\Omega'_i$ one case at a time.  Only the
the completely nonlocal term $\Omega'_0$ requires the pointwise decay
estimate (\ref{thmrem2b}); for the remaining terms the uniform bound
(\ref{thmrem2a}) suffices.  Bounding $\int \Omega_3'$ requires some care
but we shall see that the integrals of $\Omega_1'$ and $\Omega_2'$ can
be estimated heavy-handedly.

\ms\ni \underline{Case 1: Terms of the form $Rem_2^4$.}

Let $Z_2\subset Z$ denote the subset of connections for which both $
{\rm dist}(p,p_A)$ and ${\rm dist}(q,p_A)$ are $\geq c_1\l^{1/2}$,
where $c_1$ is as in Proposition \ref{thmrem2}, and let $Z_1=Z-Z_2$.  The
sets $Z_1,Z_2$ are the inverse images of sets $W_1,W_2\subset
(0,\l_0)\times N$ under the map sending a connection to its scale and
center. If for each $\l\in (0,\l_0)$ we define
$W_{1,\l}:=\{ y\in N\mid (\l,y)\in W_1
\}$, then $W_{1,\l}$ is contained in the union of a ball of radius
$\lesssim \l^{1/2}$ centered at $p$ and a similar ball centered at
$q$, so ${\rm Vol}(W_{1,\l})\lessim\l^2$.

For the orthonormal set $\{ (\bfb_n,a_n,\a_n) \}$ we may choose four
elements of the type $(\bfb,0,0)$ and four of the type $(0,*,*)$, all
normalized to unit length.  Then from (\ref{thmrem2b}) on $Z_1$ we have
\bearray\non
|Rem_2^4(X_1,\dots,X_n)| &\lessim& \l^{-4+4\d}\cdot \{ \mbox{\rm
coefficient of}\ b^4a^4\ {\rm in}\ (b^2 +ba \l^{-1/2}+ a^2\l^{-1/2})^4 \}
\\
&\lessim& \l^{-6+4\d}.
\eearray
Hence from (\ref{psiwprime})
\be
\int_{\rho_i^{-1}(Z_1)} {\rho_i^*\Omega'_0} \lessim \int_{W_1}
\l^{-6+4\d} \l^3 d\l \ dvol_N
\lessim \int_0^{\l_0} (\l^{-3+4\d}\ {\rm vol}(W_{1,\l}) )d\l 
\lessim \l_0^{4\d};
\ee
the integral over the gluing-parameter space $S^3$ gives a constant
factor. 

Similarly, on $Z_2$, $ |Rem_2^4(\tx_1,\dots,\tx_n)| \lessim
\l^{-2}r_A(p)^{-2}r_A(q)^{-2}$; the two distances $r_A(p), r_A(q)$
enter this way because in the $Rem_2^4$ term in
$\mu(p)\wedge\mu(q)$, two of the $Rem_2$'s are evaluated at $p$ and
two at $q$ (see (\ref{mainapprox})).  Since
$r_A(p)^{-2}r_A(q)^{-2}\leq r_A(p)^{-4}+ r_A(q)^{-4}$ and 
on $W_2$ both $r_A(p)$ and $r_A(q)$ are $\geq c\l^{1/2}$, we have
\be
\int_{\rho_i^{-1}(Z_2)} {\rho_i^*\Omega'_0} 
\lessim \int_0^{\l_0} \l^3 d\l \left( \int_{c\l^{1/2}}^{{\rm
diam}(N)} \l^{-2}r^{-4}r^3dr\right) 
\lessim \int_0^{\l_0} \l|\log\l|\ d\l
\lessim \l_0^{1.99}.
\ee
Combining this with the integral over $Z_1$ and summing over $i$, 
\be
\int_{Z}\Omega'_0 \lessim
\l_0^{4\d}.
\ee

\bs\ni \underline{Case 2: Terms of the form $F\wedge Rem_2^3$.}

In this and the remaining cases, $F(X_i,X_j)$ is computed either at $p$ or
at $q$, and since the $X_i$ are cut off outside a ball of radius $\sim
\l^{1/2}$ centered at $p_A$, for $i\geq 1$ terms of the form $F^i\wedge
Rem_2^{4-i}(X_1, \dots, X_8)|_{p,q}$ vanish unless $(\l,p_A)$ lies in
the set $Z_1$ defined in Case 1.  All points $p_A$ in the
remaining computations can thus be assumed to lie in a single one of our
sets $U_j$, and $\int_{\rho_j^{-1}(Z)}\rho^*\Omega'_i=
\int_Z\Omega'_i$. 

Note that all vector fields $X,Y\in \mathfrak{h}_A$ satisfy $|X|,|Y|\leq
\b(b+a\l^{-2}r_A)\leq b+a\l^{-1/2}$, and hence $|F(X,Y)|\leq
|F|(b+a\l^{-1/2})^2$.  
Using the uniform bound
(\ref{thmrem2a}) to  estimate the three $Rem_2$ terms, we obtain the
pointwise bound 
\bearray\non
|F\wedge Rem_2^3(X_1,\dots,X_8)| &\lessim& |F|
\cdot \{ \mbox{\rm
coefficient of}\ b^4a^4 \ {\rm in} \\
&& \non \ \ \ 
 (b+a\l^{-1/2})^2\l^{-3+3\d} (b^2 +ba
\l^{-1/2}+ a^2\l^{-1/2})^3 \}\\ &\lessim& |F|
\l^{-5+3\d}
\eearray
where $F$ is evaluated at either $p$ or $q$.  Because of the cutoff
in $X_i$ we may assume that $p_A$ is a distance $\lessim c\l^{1/2}$
from whichever of these points at which we evaluate.  Hence, using
(\ref{ptwsfest2}),
\be
\int_{Z} \Omega'_1
\lessim 
\int_0^{\l_0} \l^3 d\l \left(
\int_0^{c\l^{1/2}}  \frac{\l^{-5+3\d}\l^2}{(\l^2+r^2)^2} r^3 dr \right)
\lessim
\int_0^{\l_0} \l^{3\d}|\log\l| d\l 
\lessim \l_0^{1+2\d}.
\ee

\bs\ni \underline{Case 3: Terms of the form $F^2\wedge Rem_2^2$.}

Here there are two sub-cases, depending on where the points at which
$F$ and $Rem_2$ are evaluated; we can have terms of type
$F(p)F(p)Rem_2(q)Rem_2(q)$ or of type $F(p)F(q)Rem_2(p)Rem_2(q)$. In
each sub-case we bound the $Rem_2$ terms using
(\ref{thmrem2a}). At whichever point $F(X_i,X_j)$ is evaluated, we
can again assume $r_A\lessim \l^{1/2}$, so that $|X_i|\lessim b +
a\l^{-1/2}$. Letting $p',p''$ denote either of $p,q$, we then have
\bearray
|F^2\wedge Rem_2^2(X_1,\dots,X_8)| &\lessim& \non |F|(p')|F|(p'')
\cdot \{ \mbox{\rm
coefficient of}\ b^4a^4 \ {\rm in}\\
&& \non \ \ \ 
 (b+a\l^{-1/2})^4\l^{-2+2\d} (b^2 +ba
\l^{-1/2}+ a^2\l^{-1/2})^2 \}\\ &\lessim& (|F|^2(p')+|F|^2(p'') )
\l^{-4+2\d}.
\eearray
Hence the integral of the different types of $F^2Rem_2^2$ terms can all be
bounded by the integral of $|F|^2(p)\l^{-4+2\d}$:
\be
\int_{Z} \Omega'_2
\lessim
\int_0^{\l_0} \l^3 d\l \left(
\int_0^{c\l^{1/2}}  \frac{\l^{-4+2\d}\l^4}{(\l^2+r^2)^4} r^3 dr \right)
\int_0^{\l_0} \l^{-1+2\d}|\log\l| d\l 
\lessim \l_0^{\d}.
\ee

\bs\ni \underline{Case 4: Terms of the form $F^3\wedge Rem_2$.}

In the previous two cases we were rather wasteful in bounding $|X_i|$
pointwise; this time we must be more economical.

Since $p$ and $q$ enter the problem symmetrically it suffices to deal
with terms of the form $F(p)F(p)F(q)Rem_2(q)$.  Temporarily write
$r_p=r_A(p)$, $r_q=r_A(q), F_p=|F|(p), F_q=|F|(q)$.  Note that for our
term to be nonzero, both $r_p$ and $r_q$ must be $\leq c\l^{1/2}$.
Using this fact several times we find
\bearray
|F^3\wedge Rem_2(X_1,\dots,X_8)| &\lessim& \non F^2_pF_q
\cdot\l^{-1+\d} \cdot \{ \mbox{\rm
coefficient of}\ b^4a^4 \ {\rm in}\\ && \non
(b+a\l^{-1}r_p)^4(b+a\l^{-1}r_q)^2(b^2 +ba
\l^{-1/2}+ a^2\l^{-1/2}) \} \non \\
&\lessim& \l^{-4+\d}(F^2_pF_q r_p^2 +F^2_pF_q r_q^2)
\non \\
&\lessim&
\l^{-4+\d}\left(
F^4_pr_p^4 + F_q^2 + \l^{2} F^3_p
+\l^{-4}F^3_qr_q^6
\right).\label{case4}
\eearray
We can now replace $p$ by $q$ and integrate over the region $\{
(\l,p_A) \mid 0<\l\leq \l_0, \ 0\leq r_A(p)\leq c\l^{1/2}$ as in cases
2 and 3.  For each of the four terms $\l^i|F|^jr^k$ in parentheses in
(\ref{case4}), one finds $\int_0^{c\l^{1/2}} \l^i
(\l^2/(\l^2+r^2)^2)^j r^k r^3 dr\leq {\rm const}$, so
\be
\int_{Z} \Omega'_3
\lessim
\int_0^{\l_0} \l^3 \l^{-4+\d}d\l 
\lessim \l_0^{\d}.
\ee

Combining the four cases, this proves that $\int_Z\Omega'\leq c\l_0^\d$
(assuming Proposition \ref{thmrem2}).  

Now define $Rem_2^{\rm true}(X,Y)$ to be the right-hand side of
(\ref{defrem2}), but with $\xi_X,\xi_Y$ replaced by the objects
$\xi'_X, \xi'_Y$ of ({\bf Z5}).  The form $\Omega$ is obtained from
$\Omega'$ by replacing each occurrence of $Rem_2$ with $Rem_2^{\rm
true}$.  Hence $\Omega-\Omega'$ can be expressed as a sum of terms of
the form $F^i (Rem_2^{\rm true}-Rem_2)^jRem_2^k$ for appropriate
$i,j,k$. By part (b) of Proposition \ref{thmrem2}, the bounds on
$|(Rem_2^{\rm true}(X,Y)-Rem_2(X,Y)|$ are of precisely the same form
as in part (a), so the same argument as above shows that 
$\int_Z(\Omega-\Omega')\leq c\l_0^\d$, establishing
(\ref{propintrema}). 
\qed


\section{The proof of Proposition \ref{thmrem2}}\label{sectpfthmrem2}
\setcounter{equation}{0}

     The proof of Proposition \ref{thmrem2} is long, so we outline the
strategy.  To obtain (\ref{thmrem2a}), we need a pointwise bound on
$G^A_0(Rem_2'(X,Y))$ (see (\ref{defrem2})).  If there were a
four-dimensional Sobolev embedding $L_2^2\embed C^0$, then modulo
extra terms arising from Weitzenb\"{o}ck identities that occur when
comparing objects of the form $\lnorm
\na^A\na^A\phi\rnorm_2$ to objects of the form $\lap^A\phi$, we could 
get a $C^0$ bound on $Rem_2(X,Y)$ from an $L^2$ bound on
$Rem_2'(X,Y)$. (This, in turn, would require some $L^p$ and/or
pointwise bounds on $\xi$.  )

Of course there is no embedding $L_2^2\embed C^0$, but since the
failure is borderline, {\em any} stronger Sobolev-type norm should
give an embedding into $C^0$.  The most efficient Sobolev inequality
for our purposes is the following one. This inequality is not
surprising, but may not be widely known, so we prove it in the
appendix (Corollary
\ref{mainsobolev}). 

\begin{lemma} \label{sobolev} ({\bf Sobolev Embedding Lemma.}) 
Let $E$ be a vector bundle over a compact
4-dimensional manifold $N$.  For $p\in N$, let $r_p$ denote distance
to $p$.  Then for any $\d>0$
there exists a constant $c(\d)$ such that for all connections $\na$
on $E$, all $\phi\in \Gamma(E)$,
and all $p\in N$,
\be\label{soboleva}
|\phi(p)|\leq c(\d)(\lnorm\phi\rnorm_2 + \lnorm
r_p^{-\d}\na\na\phi\rnorm_2).
\ee
Hence
\be\label{sobolevb}
\lnorm \phi \rnorm_\infty \leq c(\d)\sup_{p\in N}(\lnorm\phi\rnorm_2 + \lnorm
r_p^{-\d}\na\na\phi\rnorm_2).
\ee
\end{lemma}

   We will use this lemma to get pointwise bounds on
$\phi=G^A_0(Rem_2'(X,Y))$. Hence we will need to estimate $\lnorm
G^A_0\w\rnorm_2$ and $\lnorm r_p^{-\d}\na^A\na^A G^A_0\w\rnorm_2$ for
$\w=Rem_2'(X,Y)$.  For general $\w$, Proposition \ref{lemma3}
estimates these in terms of weighted $L^2$ norms of $\w$, providing
bounds whose only dependence on the connection is explicitly through
the center point and scale.  (This type of uniformity in the
connection is the hard part of all our elliptic estimates. Uniformity
is important since to estimate an integral over a family of
connections we cannot use any bounds that depend on the connection in
an uncontrolled way.)  Proposition \ref{lemma3} also provides similar
estimates of objects $\xi$ of the form appearing in (\ref{defxix}),
which we need for reasons discussed below.

The pointwise estimates of $G^A_0\w$ in terms of weighted $L^2$ norms
of general $\w$ will be summarized (and generalized) as part of
Proposition \ref{lemma3}, specifically the first half of
(\ref{lemma3e}).  To apply these general estimates to $\w=
Rem_2'(X,Y)$ we still need to bound the weighted $L^2$ norms of
$Rem_2'(X,Y)$.  To understand what this entails, write
$Rem_2'(X,Y)=Rem'_{2,{\rm loc}}+Rem'_{2,{\rm semiloc}}+Rem'_{2,{\rm
nonloc}}$, where
\be\label{deflocnonloc}
Rem'_{2,{\rm loc}}=R''(X,Y), \ \ 
Rem'_{2,\rm semiloc}=\{\tx,\xi_Y\} -\{\tilde{Y},\xi_X\}, \ \
Rem'_{2,{\rm nonloc}}=\{\xi_X,\xi_Y\}
\ee
(see (\ref{defr2prime})).  Because of the cutoffs in $X$ and $Y$, the
expressions $Rem'_{2,{\rm loc}}$ and $Rem'_{2,{\rm semiloc}}$ are
supported in $B(p_A,2\e)$, but $Rem'_{2,{\rm nonloc}}$ is not. Thus
among the estimates we need are weighted $L^2$ bounds on $R''(X,Y)$.
By Lemma \ref{oldlemma39}, below, pointwise we find
\be\label{locterm}
|R''(X,Y)|\leq c |\hat{X}| |\hat{Y}|
(\b+\e^{-2}\chi)(|F|+r_A|\Na F|)
\ee
(recall that $\hat{X}$ is the object that the cutoff $\b$ multiplies
in the definition of $X$).  Here $\chi$ is the characteristic function
of the annulus $\e \leq r_A\leq 2\e$. Thus to apply the estimate
(\ref{lemma3e}) of Proposition \ref{lemma3}
to obtain bounds on $G^A_0(Rem_2'(X,Y))$, we need to estimate certain
expressions of the form $\lnorm \b r_p^{-\d}r_A^m F \rnorm_2,$ and
similar expressions with $F$ replaced by $\Na F$ and/or with
$\b$ replaced by $\chi$.  This will be accomplished in Lemma
\ref{cor10_simp}, where we will list all the purely local estimates we need.

Weighted $L^2$-norm bounds on $Rem'_{2,{\rm semiloc}}$ and
$Rem'_{2,{\rm nonloc}}$ can be obtained from weighted $L^4$-norm
bounds on $\tx$ and $\xi$.  The first of these is another purely local
estimate. The second will be 
achieved in Proposition \ref{pigcor11}, where we will use
the basic elliptic tools in Proposition \ref{lemma3} to turn the
problem into a local estimate again.

   Till now we have made no mention of the role the point $p$ plays in
affecting the weighted norms.  If we compute these weighted norms as
above and take the supremum over $p\in N$ as in (\ref{sobolevb}), we
obtain only the sup-norm bound (\ref{thmrem2a}) for $Rem_2(X,Y)$.  To
prove Proposition \ref{propintrem} we additionally need the pointwise
decay bound (\ref{thmrem2b}).  Since the local quantities we bound
are supported near the center point $p_A$ of $A$, decay is only an
issue for the nonlocal quantities, but these are built out of Green
operators applied to quantities supported near $p_A$.  Thus one
expects that as the distance between $p$ and $p_A$ increases, the
bounds on our non-local quantities should decrease. This turns out to
be true (at least for ${\rm dist}(p,p_A)\geq {\rm const}\cdot
\l^{1/2}$); we simply have to work harder, establishing some general
pointwise bounds in Proposition \ref{estlemma5}.  Our basic estimates
in Proposition
\ref{lemma3} are most useful for $p$ close to $p_A$; to get 
the bounds that lead to (\ref{thmrem2b}), in which $p$ is farther from
$p_A$, we will apply Proposition \ref{estlemma5}.

   To establish (\ref{thmrem2b}) we again break up $Rem_2'(X,Y)$ into
its local, semi-local, and non-local pieces as in
(\ref{deflocnonloc}).  In the cases of $Rem'_{2,{\rm loc}}$ and
$Rem'_{2,{\rm semiloc}}$, Proposition \ref{estlemma5} again reduces
our work to weighted $L^p$ bounds of purely local quantities. For
$Rem'_{2,{\rm nonloc}}$, however, Proposition
\ref{estlemma6} leaves us with bounding an expression of the form
$\|r_{A}^{1+\d'}\{\xi_X,\xi_Y\}\|_2$, and the obvious
approach---H\"{o}lder's inequality and the weighted $L^4$ bounds
already obtained---does not give us a strong enough bound for an
adequate decay rate in (\ref{thmrem2b}).  We will circumvent this by
obtaining a pointwise decay estimate for $\xi$, which in turn gives us
a satisfactory decay rate for $Rem_2$. (In fact, with the pointwise
estimate on $\xi$ in hand, it turns out that the contribution of 
$Rem'_{2,{\rm nonloc}}$ to $Rem_2$ 
%
%
is much smaller than the bounds we obtain from the other two terms.)

   With this discussion behind us, our procedure is clear.  First we
will fill our elliptic toolbox by proving Propositions \ref{lemma3} and 
\ref{estlemma5}.  To apply these we need to compute weighted $L^p$
norms of various quantities appearing in $Rem_2'$, which is our next
step. The final step is then a matter of bookkeeping, applying the
general elliptic tools to bound $G^A_0$ of $Rem'_{2,{\rm
loc}},Rem'_{2,{\rm semiloc}}$, and $Rem'_{2,{\rm nonloc}}$.

   To avoid writing similar hypotheses over and over, and for
notational simplicity, for the rest of this section we impose the
following

\bs \ni{\bf Blanket hypotheses and notation.}  $A$ always denotes a
connection with $[A]\in \M_{k+1,\l_0}^{\k,\nu}$ (see Definition
\ref{restricted_M}).  Every Proposition, Lemma, etc. has an implicit
hypothesis ``for $\l_0$ sufficiently small and for all $ [A]\in
\M'_{k+1,\l_0}\intersect \M^{\k,\nu}$''.  Constants $c$ are uniform in
$A$ and in all other relevant parameters not explicitly shown (though
some would depend on $\k$ and $\nu$, if these were not fixed);
e.g. $c(\d)$ depends only on $\d$.  Constants are continually updated,
and when a hitherto unnamed $c$ appears, there is an implicit ``for
some constant $c$''.  The notation ``$x\lessim y$'' means $x\leq cy$.
$F$ always denotes the curvature of the conection $A$, and $\na=\Na$
denotes the full covariant derivative on $\Gamma(\adp \tensor
\Lambda^*T^*N)$ (the tensor product connection determined by 
$A$ and the Levi-Civita connection).  Given any vector field $X$ on
$N$, we write $\tx=\i_XF$ (thus there is an $A$-dependence we
suppress).  We write $p_A$ for the center point of $A$ and $\l$ for
$\l(A)$.  For any $p\in N$, we let $r_p$ denote distance to $p$, and
write $r_A$ for $r_{p_A}$.  When a point $p$ appears in a hypothesis,
the letter $d$ always means $r_A(p)={\rm dist}(p,p_A)$.  The scale $\e={\rm
const}\ \l^{1/2}$ and cutoff $\b=\b_{\rm std}(r_A/e)$ are always as in
(\ref{frakha}), and $\chi$ denotes the characteristic
function of the annulus $\{ \e \leq r_A \leq 2\e \}$ containing the
support of $d\b$.  We also define the operators
$\D=\D^A:\Omega^1(\adp)\to \Omega^0(\adp)\plus\Omega^2_+(\adp)$ by
\be\label{defd}
\D^A\eta=\left((d^A)^*\eta,\ \sqrt{2}d^A_+\eta\right);
\ee 
thus $\ker(\D^A)=H^1_A$, the harmonic space in the middle of the
elliptic complex
\be\label{fes}
0  \to \Omega^0(\adp) {\buildrel d^A \over \longrightarrow}
 \Omega^1(\adp) {\buildrel \sqrt{2}d^A_+ \over \longrightarrow} 
\Omega^2_+(\adp) \to 0 .
\ee
Define $\lap^A_0,\lap^A_1, \lap^A_+$ to be the Laplacians on
zero-forms, 1-forms, and SD two-forms, respectively, constructed from
this complex, and let $G^A_0,G^A_+$ be the inverses of
$\lap^A_0,\lap^A_+$.  Also define $\lapplus,
\gplus$ on $\Omega^0(\adp)\plus\Omega^2_+(\adp)$ by
$\lapplus=\lap^A_0\plus \lap^A_+$, $\gplus=G^A_0\plus G^A_+$. 
Note that 
\be\non\label{defdstar}
(\D^A)^*(\phi_0,\phi_+) = d^A\phi_0\ +\ \sqrt{2}(d^A_+)^*\phi_+,
\ee
so that the quantity $\xi_X=\xi^A_X$ of (\ref{defxix}) can be written
as
\be\label{xi2}
\xi_X=(\D^A)^*\gplus\D^A\tx.
\ee
Finally, observe that
\be\label{ddstar}\label{dstard}
\D^A(\D^A)^*=\lapplus, \ \ 
(\D^A)^*\D^A =\lap^A_1.
\ee

Now we can finally begin proving Proposition \ref{thmrem2}. In the
following proposition, what drives the estimates are two facts:
(i) $\lapplus$ is uniformly bounded below, and (ii) in the
Weitzenb\"{o}ck identity for $\lapplus$, only Riemannian curvature
terms appear; $F$ does not enter.

\begin{prop} \label{lemma3}\label{estlemma3}
For $\d_0>0$ sufficiently small and any $\d,\d',\d''$ (possibly zero) of
absolute value less than $\d_0$, 
such that for any $p\in N$ and any $\w\in
\Omega^0(\adp)\plus\Omega^2_+(\adp)$

\be\label{lemma3a} 
\|\gplus\w\|_2\lessim\|r_p^{1+\d'}\w\|_2
\ee
and 
\be\label{lemma3b}
\lnorm r_p^{-\d} \Na\Na \gplus \w \rnorm_2
\lessim 
\|r_p^{-\d}\w\|_2 + \l^{\d'-1}\|r_A^{1-\d-\d'}\w\|_2.
\ee

\ss
Furthermore, if 
$\xi=(\D^A)^*\gplus\w\in\Omega^1(\adp)$ (cf. (\ref{xi2})), then
\bearray
\label{lemma3c}
\lnorm \xi \rnorm_2 &\lessim& \lnorm r_A^{1+\d} \w \rnorm_2,\\
\lnorm r_A^{1-\d}\xi\rnorm_4+\lnorm r_A^{1-\d}\na^A\xi\rnorm_2
&\lesssim&
\lnorm r_A^{1-\d} \w\rnorm_2,
\label{lemma3d3}\label{pigtechpropb}\\
\lnorm r_p^{-1-\d}\xi\rnorm_2 + 
\lnorm r_p^{-\d} \xi\rnorm_4+\lnorm r_p^{-\d} \na^A\xi\rnorm_2
&\lesssim&
\lnorm r_p^{-\d}\w\rnorm_2
+\l^{\d'-1}\lnorm r_A^{1-\d-\d'}\w\rnorm_2, 
\label{lemma3d2}\label{pigtechpropa}
\eearray
and
\be\label{lemma3d}
\lnorm r_p^{-\d}\Na\xi\rnorm_4 +\lnorm r_p^{-\d}\Na\Na\xi\rnorm_2 \lessim 
\lnorm r_p^{-\d} (\D^A)^*\w \rnorm_2 + 
\l^{\d'-1}\lnorm r_A^{-\d-\d'}\w\rnorm_2+
\l^{\d'-2}\lnorm r_A^{1-\d-\d'}\w\rnorm_2.
\ee
\ss As a corollary of Lemma \ref{sobolev} 
and (\ref{lemma3c}--\ref{lemma3d}), if $\d_0>0$ is sufficiently small
and $0<\d<\d_0)$
then
\be\label{lemma3e}\label{lemma3f}
|\gplus\w|(p)\leq c(\d)\cdot\left\{\mbox{\rm RHS of (\ref{lemma3b})}\right\}
 \ \ \mbox{\rm and}
\ \ 
|\xi(p)|\lessim c(\d)\cdot\left\{\mbox{\rm RHS of (\ref{lemma3d})}\right\}.
\ee

\end{prop}

We remark that in (\ref{lemma3b}), (\ref{lemma3d2}), and (\ref{lemma3d})
it is important that $r_A$ appear
where it does, rather than $r_p$, or we would not get strong enough
estimates in our applications.  The fact that both $r_A$ and $r_p$
appear together in Proposition \ref{lemma3} complicates its proof.

\ms\ni\pf A slightly less general set of bounds was derived in Lemma 3.3
of \cite{Gppt} for $\lap^A_+$, the Laplacian on SD two-forms only, but
for the reasons mentioned prior to stating the proposition,
essentially the same proof works here.  The only differences are that
(i) in \cite{Gppt} the decay (\ref{ptwsfest2}) was true on all of $N$,
not merely in $B(p_A,2\e)$, and (ii) \cite{Gppt} dealt only with the
case $p=p_A$.  Since the cited proof is rather long, we will not
repeat the parts that require only minor modifications, and will jump
to the points of departure.

To establish (\ref{lemma3a}), the proof of Lemma 3.3a in \cite{Gppt}
works verbatim to show that
\be\label{estlemma3a}
\|r_p^{-\d-1}\gplus\w\|_2+\|r_p^{-\d}\gplus\w\|_4+
\|r_p^{-\d}\na \gplus\w\|_2
\lessim\|r_p^{-\d+1}\w\|_2.
\ee 
Note that $\d$ need not be positive here.
Since $|\gplus\w|\lessim |r_p^{\d-1}\gplus\w|$, (\ref{lemma3a})
follows.

Moving to (\ref{lemma3b}), let $\eta\in\Omega^*(\adp)$ be a form of
arbitrary degree. The procedure in \cite{Gppt} for proving its Lemma
3.3b,c---squaring, integrating by parts, commuting a covariant
derivative past a trace Laplacian $(\Na)^*\Na=\na^*\na$, and juggling
terms---leads to
\bearray\non
\|r_p^{-\d}\na\eta\|_4+\|r_p^{-\d}\na\na\eta\|_2&\leq& c \left(
\|r_p^{-\d}\na^*\na \eta \|_2 +
\|r_p^{-\d}\eta\|_2 +\|r_p^{-\d}\na\eta\|_2 
\right. 
\\ && \left.+
\|r_A^{\d+\d'}r_p^{-\d}F_A\|_4\|r_A^{-\d-\d'}\na\eta\|_2
\right); \label{estlemma3j}\label{twoderivs2}
\eearray
here the smallness of $|\d|$ has also been used to ensure that the
term $|\d| \|r_p^{-\d-1}\na\eta\|_2$ that initially comes up on the
right-hand side is $\lesssim |\d|\|r_p^{-\d}\na\na\eta\|_2$; see
\cite{Gppt}, Lemma 3.2. 
(In \cite{Gppt} there was no need to insert $r_A^{\pm (\d+\d')}$.) 
First consider the case
$\eta=\gplus\w$, where $\w\in \Omega^0(\adp)\plus\Omega^2_+(\adp)$.
The Weitzenb\"{o}ck formula gives
$\lapa\eta=\w+\R(\gplus\w)$, where $\R$ is an endomorphism
proportional to the Riemann tensor.  Moreover we will see in Lemma
\ref{cor10_simp}b below that for $\d,\d'$ sufficiently small,
$\|r_A^{\d+\d'}r_p^{-\d}F_A\|_4\lessim \l^{\d'-1}$.  Inserting these
facts into (\ref{estlemma3j}), one can continue the argument as in
\cite{Gppt} and arrive at an extended version of (\ref{lemma3b}):
\bearray\non
\|r_p^{-\d-1}\na \gplus\omega\|_2+
\|r_p^{-\d}\na \gplus\omega\|_4 + \|r_p^{-\d}\na\na \gplus\omega\|_2 
&\lessim&
\|r_p^{-\d}\w\|_2 + \l^{\d'-1}\|r_A^{1-\d-\d'}\w\|_2.\\
\label{estlem3c}
\eearray

As for (\ref{lemma3c}), since $|\D^*\eta|\leq c |\na\eta|$, the
desired estimate follows from (\ref{estlemma3a}). 

By similar manipulations one can also establish 
\be
\|r_A^{-\d+1}\na \gplus\omega\|_4+ \|r_A^{-\d+1}\na\na \gplus\omega\|_2
\lessim
\|r_A^{1-\d}\w\|_2.
\label{estlem3b} 
\ee
Since $\D$ is $\Na$ followed by a covariantly constant
projection, the same bounds hold with $\na\gplus\w$ replaced by
$\D^*\gplus\w=\xi$, yielding (\ref{lemma3d3}). For the same reason,
(\ref{lemma3d2}) follows from (\ref{estlem3c}).

   Finally, to establish (\ref{lemma3d}), return to (\ref{estlemma3j})
and use the Weitzenb\"{o}ck formula for 1-forms,
\be
\na^*\na\xi=\lap^A_1\xi +\F(\xi) +\R(\xi).
\ee
Here $\F$ is an endomorphism proportional to $F$. 
Since $\lap^A_1\xi=\D^*\D (\D^* \gplus\w) = \D^*\w$ (see (\ref{ddstar})),
we have 
\be
\|r^{-\d}\na^*\na \xi \|_2\leq c\left(
\lnorm r^{-\d}\D^*\w \rnorm_2
+\lnorm r^{-\d} \xi \rnorm_2 +\lnorm r_A^{\d+\d'}r^{-\d}F\rnorm_4
\lnorm r_A^{-\d-\d'} \xi \rnorm_4\right),
\ee
Hence
\bearray\non
\|r^{-\d}\na\xi\|_4+\|r^{-\d}\na\na\xi\|_2&\leq& c \left(
\lnorm r^{-\d} (\D^A)^*\w \rnorm_2
+\lnorm r^{-\d}\xi\rnorm_2 +\lnorm r^{-\d}\na\xi\rnorm_2 \right. \\ &&
\left.  +\lnorm r_A^{\d+\d'}r^{-\d}F^A \rnorm_4 \left(\lnorm
r_A^{-\d-\d'} \xi \rnorm_4 +\lnorm r_A^{-\d-\d'}\na \xi
\rnorm_2\right)
\right).\non \\
\label{estlem3m}
\eearray
Once again $\|r_A^{\d+\d'}r^{-\d}F\|_4\lessim \l^{\d'-1}$, and
(\ref{estlem3c}) (with $\na\gplus\w$ replaced by $\xi$)
implies 
\be
\left(\lnorm r_A^{-\d-\d'} \xi \rnorm_4
+\lnorm r_A^{-\d-\d'}\na \xi \rnorm_2\right)
\leq c\left(
\|r_A^{-\d-\d'}\w\|_2+ \l^{-1+\d''}\|r_A^{-\d-\d'-\d''}\w\|_2\right).
\ee
Using (\ref{estlem3b}--\ref{estlem3c}) to bound the other terms in
(\ref{estlem3m}), the bound (\ref{lemma3d}) follows.
\qed

Proposition \ref{lemma3} gives the same bounds for all $p\in N$; to
obtain (\ref{thmrem2b}), we need estimates that show decay as
$d=r_A(p)$ grows.  The following proposition provides these
estimates. We separate the estimates into cases (a) and (b) below
because for many purposes the only $\w$'s for which we need to
estimate the quantities in Proposition \ref{lemma3} are compactly
supported in a $2\e$-ball around $p_A$, and we get sharper estimates
in this case. Part (a) will thus be used to bound the terms
$G^A_0R''(X,Y)$ and $G^A_0\{ \tx,\xi_Y
\}$  in $Rem_2(X,Y)$; part (b) will be used to bound
$G^A_0\{\xi_X,\xi_Y\}$.

\begin{prop}\label{estlemma5}\label{estlemma6}
Notation as in Proposition
\ref{lemma3}. There exists $\d_0>0$ such that the following are true.

(a) Suppose that for some $\e_0$ (not necessarily related to
$\e=c\l^{1/2}$, and allowed to depend on $\w$), (i) $\supp(\w)\subset
B(p_A,\e_0)$, (ii) $d=r_A(p)={\rm dist}(p,p_A)\geq 2\e_0$, and (iii)
$|F^A|\leq B$ on the complement of $\supp(\w)$.  Let $\tilde{\b}$ be a
cutoff function of the form $\b_{\rm std}(4r_p/d)$ (so that
$\supp(\tilde{\b})\subset B(p,d/2)$).  Then for any $\d'$ with
$|\d'|\leq \d_0$, and any $\d\in (0,\d_0)$, 
\be\label{estlemma5a6}
|\gplus\w|(p) \leq c(\d) d^{-1-\d-\d'}(1+B^{1/2}) \lnorm
r_A^{1+\d'}\w\rnorm_2 
\ee
and 
\be\label{estlemma5a7}
|\xi|(p) \leq c(\d)d^{-\d-\d'}(d^{-2}+B) \lnorm
r_A^{1+\d'}\w\rnorm_2 
\ee
Thus, if ${\rm supp}(\w)\subset B(p_A,2\e)$ and $d\geq 4\e=c\l^{1/2}$,
then using (\ref{ptwsfest}),
\be\label{estlemma5c}
|\gplus\w|(p)\leq c(\d) d^{-1-\d-\d'}
\lnorm r_A^{1+\d'}\w\rnorm_2
\ee
and
\be\label{estlemma5d}
|\xi|(p) \leq c(\d)d^{-2-\d-\d'} \lnorm
r_A^{1+\d'}\w\rnorm_2 
\ee

\ms
(b) Suppose only that $|F^A|\leq B$ on $B(p,d/2)$, where $d={\rm
dist}(p,p_A)>0$; suppose nothing about the 
support of $\w$. Let $\tilde{\b}$ be as in (a). 
Then for
all $\d'$ with $|\d'|\leq \d_0$, we have
\be\label{estlemma6c}
|\gplus\w|(p) \leq c(\d) (1+B^{1/2})\left( d^{-1-\d-\d'}
\lnorm r_A^{1+\d'}\w\rnorm_2 + \lnorm r_p^{-\d}\tilde{\b}\w
\rnorm_2\right). 
\ee

Thus, if $d\geq c\l^{1/2}$, then 
\be\label{estlemma6d}
|\gplus\w|(p) \leq c(\d) \left( d^{-1-\d-\d'}
\lnorm r_A^{1+\d'}\w\rnorm_2 + \lnorm r_p^{-\d}\tilde{\b}\w
\rnorm_2\right). 
\ee

\end{prop}

\pf 
(a) We will apply the Sobolev inequality (\ref{soboleva}), but first
we must bound $\|r_p^{-\d}\na\na (\tilde{\b} \gplus\omega)\|_2, \
\lnorm \tilde{\b} \gplus\w\rnorm_2$, and similar expressions with
$\gplus\w$ replaced by $\xi$.

(i) First we will show that
\bearray
\label{estlemma5a}
\|r_p^{-\d}\na\na (\tilde{\b} \gplus\omega)\|_2 
&\lessim d^{-1-\d-\d'}
(1+B^{1/2})\lnorm r_A^{1+\d'}\w\rnorm_2.
\label{estlemma5a2} \eearray

Let $\eta\in \Omega^0(\adp)\plus\Omega^2_+(\adp)$.  Proceed as in the
proof of Proposition \ref{lemma3b}---squaring, integrating by parts,
etc.---but this time leave the term proportional to $F$ (which arises
from commuting $\Na$ past a trace-Laplacian) in integrated form.  One
arrives at
\be
\|r_p^{-\d}\na\na\eta\|_2^2 \lessim
\|r_p^{-\d}\lap\eta\|_2^2 +
\|r_p^{-\d}\eta\|_2^2 +\|r_p^{-\d}\na\eta\|_2^2
\int r_p^{-2\d}|F||\na\eta|^2
\label{estlemma5b}
\ee
where $\lap=(\Na)^*\Na$.
Now replace $\eta$ by $\tilde{\b}\eta$. In the integral we have
$|F|\leq B$, so 
\be
\|r_p^{-\d}\na\na(\tilde{\b}\eta)\|_2^2\leq c \left(
\|r_p^{-\d}\lap(\tilde{\b}\eta)\|_2^2 +
\|r_p^{-\d}\tilde{\b}\eta\|_2^2 +(1+B)\|r_p^{-\d}\na(\tilde{\b}\eta)\|_2^2
\right). \label{estlemma5c2}
\ee
An integration by parts, plus various steps already seen in the proof
of Lemma \ref{lemma3}, gives 
\be
\|r_p^{-\d}\na(\tilde{\b}\eta)\|_2^2
\lessim \lnorm r_p^{-\d}\lap(\tilde{\b}\eta)\rnorm_2 \ \lnorm
r_p^{-\d}\tilde{\b}\eta\rnorm_2
\leq c\left(k \lnorm r_p^{-\d}\lap(\tilde{\b}\eta)\rnorm_2^2
+k^{-1}\lnorm r_p^{-\d}\tilde{\b}\eta\rnorm_2^2\right)
\ee
for arbitrary $k$.  Inserting this into (\ref{estlemma5c2}), with
$k\lessim (1+B)^{-1}$, we find
\be
\|r_p^{-\d}\na\na(\tilde{\b}\eta)\|_2^2\lessim
\|r_p^{-\d}\lap(\tilde{\b}\eta)\|_2^2 +
(1+B)\|r_p^{-\d}\tilde{\b}\eta\|_2^2. \label{estlemma5e}
\ee
Using the Weitzenb\"{o}ck formula 
as in the proof of Lemma \ref{estlemma3}, we
can replace $\lap$ by $\lapplus$, absorbing the zeroeth-order term into
$(1+B)\lnorm r_p^{-\d}\tilde{\b}\eta\rnorm_2^2$. Additionally, 
by (\ref{estlemma3a}) we have $\lnorm r_p^{-\d}\tilde{\b}\eta\rnorm_2
\lessim \lnorm r_p^{-\d}\lapplus(\tilde{\b}\eta)\rnorm_2$.  Hence
\be
\|r_p^{-\d}\na\na(\tilde{\b}\eta)\|_2\lessim (1+B^{1/2})
\|r_p^{-\d}\lapplus(\tilde{\b}\eta)\|_2. \label{estlemma5e2}
\ee

   Next, note that for any function $f$,
\be\label{estlemma5e3}
|\lapplus(f\eta)-f\lapplus\eta|\lessim (|\na\na f|
|\eta| + |\na f| |\na \eta|).
\ee
Apply this with $f=\tilde{\b}$ and $\eta=\gplus\w$, noting that by the
hypothesis on the support of $\w$ we have
$\tilde{\b}\lapplus\eta=\tilde{\b}\eta\ident 0$.  Since
$|\na^j\tilde{\b}|\leq cd^{-j}$, and since 
on the
support of $\na\tilde{\b}$ we have both 
$d/2\leq r\leq d$ and $d/2\leq r_A\leq 3d/2$, we obtain
\bearray\non
|r_p^{-\d}\lapplus(\tilde{\b}\eta)| &\leq & c
r_p^{-\d}\tilde{\chi}(d^{-2}|\eta|+ d^{-1}|\na\eta|)\\ &\leq&
cd^{-1-\d-\d'}\tilde{\chi}(r_A^{-1+\d'}|\eta|+
r_A^{\d'}|\na\eta| )
\label{estlemma5e4}
\eearray
where $\tilde{\chi}$ denotes the characteristic function of the
annulus $d/4\leq r\leq d/2$.
Inserting this into (\ref{estlemma5e2}), we have
\bearray\non
\|r_p^{-\d}\na\na(\tilde{\b}\eta)\|_2 &\leq& (1+B^{1/2})
cd^{-1-\d-\d'}\left( \lnorm\tilde{\chi}
r_A^{-1+\d'}\eta\rnorm_2+
\lnorm \tilde{\chi} r_A^{\d'}\na\eta\rnorm_2  
\right) \label{estlemma5e5}\\
&\leq& (1+B^{1/2})
cd^{-1-\d-\d'}\left( \lnorm
r_A^{-1+\d'}\eta\rnorm_2+
\lnorm  r_A^{\d'}\na\eta\rnorm_2  
\right)\non\\ \label{estlemma5e6}
\eearray
Now apply (\ref{estlemma3a}) to obtain
\be\label{estlemma5e7}
\|r_p^{-\d}\na\na(\tilde{\b}\eta)\|_2 \leq (1+B^{1/2})
cd^{-1-\d-\d'}\lnorm r_A^{1+\d'}\w \rnorm_2,
\ee
which leads to (\ref{estlemma5a}).  

\ss
Moving on to $\lnorm \tilde{\b} \gplus\w\rnorm_2$, and 
repeating some of the steps in the proof of (a) with $\d=0$, we have
\bearray
\lnorm \tilde{\b}\gplus\w\rnorm_2
\lessim\lnorm \lapplus(\tilde{\b}\gplus\w)\rnorm_2
&\lessim& \non  d^{-1-\d'}\left(
\lnorm \tilde{\chi} r_A^{-1+\d'}\gplus\w\rnorm_2
+\lnorm \tilde{\chi} r_A^{\d'}\na \gplus\w \rnorm_2\right)\\
&\lessim&  d^{-1-\d'}\lnorm r_A^{1+\d'}\w\rnorm_2.
\label{estlemma5a4}
\eearray
This is smaller than the bound (\ref{estlemma5a}), so (\ref{soboleva})
gives (\ref{estlemma5a6}). 

\ms
   The bound (\ref{estlemma5a7}) is derived by methods similar to the
preceding and those used in Proposition \ref{lemma3}.  We leave the
details to the reader.

\ms
(b) Proceed as in (a); the only change is that no we no longer have
$\tilde{\b}\w\ident 0$. The first effect of this change occurs in
(\ref{estlemma5e4}), where we have to add $|r^{-\a}\tilde{\b}\w|$ to
the RHS.  The effect of this term is to add $(1+B^{1/2})\lnorm
r^{-\a}\tilde{\b}\w \rnorm_2$ to the RHS of (\ref{estlemma5e6}) and
(\ref{estlemma5e7}), hence to (\ref{estlemma5a}).  There is a similar
change in the bound on $\|\tilde{\b}\gplus\w\|_2$, but its effect is
smaller than the preceding one.
\qed

    To apply Propositions \ref{lemma3} and \ref{estlemma5} to estimate
$Rem_2$, we need to estimate expresssions of the form $\|r_p^m\w\|_2$ for
various $m$, where $\w=Rem'_{2, {\rm loc}}, Rem'_{2, {\rm semiloc}},$ or
$Rem'_{2, {\rm nonloc}}$ (see {\ref{deflocnonloc}}). First we deal
with the purely local object $Rem'_{2, {\rm loc}}(X,Y)=R''(X,Y)$.  To
start, we need a pointwise estimate, given by the next lemma. The
conclusion of the lemma is deceptively simple; the way in which the
derivatives of $X$ and $Y$ are coupled to each other and to $\na F$ in
the definition of $Rem_2''$ is crucial.  

\begin{lemma}\label{oldlemma39}  For $X,Y\in \mathfrak{h}_A$, 
\be
\label{piglemma7a}
|R''(X,Y)|\lessim  |\hat{X}| |\hat{Y}|
(\b+\e^{-2}\chi)(|F|+r_A|\na F|).
\ee
(Here $\hat{X},\hat{Y}$ are the un-cutoff versions of $X,Y$; see
(\ref{frakha}).)
\end{lemma}

We remark that for general vector fields, this lemma would be false.

\pf Let $\phi=R''(X,Y)$.  From (\ref{defr2prime}) we have
$R''(X,Y)=\b^2R''(\hat{X},\hat{Y}) + $ terms involving the derivative
of $\b$.  The latter are easily dealt with, giving the terms
proportional to $\chi$ in (\ref{piglemma7a}).  For
$R''(\hat{X},\hat{Y})$, the first three terms in (\ref{defr2prime})
have norm bounded by $|F|(|\hat{X}||\hat{Y}|+|\lap \hat{X}||\hat{Y}|
+|\hat{X}||\lap\hat{Y}|)$, and an easy computation shows that for
$X\in \mathfrak{h}_A$, $|\lap \hat{X}|\lessim|\hat{X}|$.  Furthermore,
because $F$ is ASD (and hence Yang-Mills as well) and the
``rotational'' parts of $\hat{X},\hat{Y}$ are SD, the remaining three
terms in $R''(\hat{X},\hat{Y})$ would vanish if the metric on $N$ were
Euclidean.  When we do the bookkeeping necessary for the $O(r_A^2)$
difference between the metric coefficients $g_{ij}$ and $\d_{ij}$, the
we obtain contributions bounded by $|\hat{X}||\hat{Y}|(|F|+r_A|\na
F|)$.
\qed

Thus, bounding $|G^A_0R''(X,Y)|$ pointwise boils down to estimates of
the form in the following lemma.

\begin{lemma}\label{cor10_simp} 
Let $p\in N$ be arbitrary and let $d=\dist(p,p_A)$.  Then we have the
following estimates.

(a) Assume  $0\leq \d<2$  and $n>-2+\d$. Then
\be
\label{cor10f_simp}
\label{cor10g_simp}
\lnorm\b r_p^{-\d} r_A^n F\rnorm_2 + 
\lnorm \b r_p^{-\d} r_A^{n+1} F\rnorm_4+
\lnorm\b r_p^{-\d} r_A^{n+1} \na F\rnorm_2  \lessim 
\left\{ \begin{array}{lc}
\l^{n-\d}, & n-\d<2,
\\
\l^2|\log\l|^{1/2} & n-\d=2,
\\
\l^{1+(n-\d)/2}, & n-\d>2
\end{array}\right. 
\ee
and for all $n$
\be
\label{cor10h1_simp}
\label{cor10h2_simp}
\lnorm\chi r_p^{-\d}r_A^n F\rnorm_2 +
\lnorm\chi r_p^{-\d} r_A^{n+1} \na F\rnorm_2 
\lessim \l^{1-\d+n/2}.
\ee

(b) Let $\e_0>\l_0^{1/2}$ be some fixed number.
For $0<\d<1$ and $-1+\d< n< 3$ we have
\be\label{cor10h5_simp}
\lnorm r_p^{-\d} r_A^n F\rnorm_4 \lesssim 1+
\left\{ \begin{array}{lr}
\l^{n-1-\d} & \mbox{\rm if}\ d\lessim \e_0, \\
\l^{n-1} & \mbox{\rm if}\ \e_0\lessim d.
\end{array}\right.
\ee

\end{lemma}

\pf (a) First consider the case $\d=0$. From (\ref{ptwsfest2})
one quickly finds 
\be
\label{cor10a_simp}\label{2norm_beta_rnf_simp}
\label{cor10b_simp}\label{4norm_beta_rnf_simp}
\lnorm \b r_A^nF \rnorm_2  + \lnorm \b r_A^{n+1}F \rnorm_4 
\lesssim 
\left\{ \begin{array}{rr}
\l^n & -2<n<2,\\ \l^2\sqll &  n=2,\\
\l^{1+n/2} & n>2
\end{array}\right.
\ee
and for all $n$, 
\be
\label{cor10d_simp}\label{2norm_chi_rnf_simp}
\lnorm\chi r_A^n F\rnorm_2 \lessim \l^{1+n/2}.
\ee
As for $\|\b r_A^{n+1}\na F\|$, the same argument as in the proof of
[G1, Lemma 3.3b] shows that
\be\label{cor10i_simp}
\lnorm \b r_A^{n+1} \na F\rnorm_2^2\lesssim
n\lnorm \b r_A^{n} F\rnorm_2^2+ 
\lnorm \b r_A^{n+1} F\rnorm_2^2+ 
\lnorm |d\b| r_A^{n+1} F\rnorm_2^2
+\int \b^2 r_A^{2n+2}|F|^3.
\ee
Since $|d\b|\lessim \e^{-1}\chi$, we have $|d\b| r_A^{n+1} \leq
c\e^{n}\chi$.  Thus
\be
\lnorm \b r_A^{n+1} \na F\rnorm_2\lesssim
\lnorm \b r_A^{n} F\rnorm_2+
\l^{n/2}\lnorm \chi F\rnorm_2
+\left(\int \b^2 r_A^{2n+2}|F|^3\right)^{1/2}.
\ee
 From (\ref{ptwsfest2}) one can deduce 
\be
\left(\int \b^2 r_A^{2n+2}|F|^3\right)^{1/2}\lesssim 
\left\{ \begin{array}{rr}
\l^{n}, & n<3,\\ \l^3\sqll, &  n=3,\\
\l^{3n/2 -3/2}, &  n>3.
\end{array}\right.
\ee
Combining this with our previous bounds we find that 
\be\label{cor10c_simp}\label{2norm_beta_rn_gradf_simp}
\lnorm \b r_A^{n+1}\na F\rnorm_2\lessim
\mbox{\rm RHS of} (\ref{cor10a_simp})
\ee
To bound $\lnorm\chi r_A^{n+1} \na F\rnorm_2$, again use the
analysis leading to (\ref{cor10i_simp}), but with $\b$ replaced by a
smooth extension of $\chi$ of the form $f(r_A/\e)$ with $f$ supported
in $[1/2,3]$.  (It is simplest first to note that since $r_A\leq 2\e$
on ${\rm supp}(\chi)$, $\lnorm\chi r_A^{n+1} \na F\rnorm_2\lessim
\l^{(n+1)/2}\lnorm\chi \na F\rnorm_2$,)
Then analysis similar to the above leads to 
\be\label{cor10e_simp}
\lnorm\chi r_A^{n+1} \na F\rnorm_2 \lessim
\mbox{\rm RHS of} (\ref{cor10d_simp}).
\ee
   This completes the case $\d=0$ and we move on to the general case. 

We first bound $\|\b r_p^{-\d}r_A^nF\|_2$; the method for bounding
$\|\b r_p^{-\d}r_A^{n+1}F\|_4$ is identical.
Break the ball $B(p_A,2\e)$ into two pieces: an inner region
$B_{in}=B(p,d/2)\intersect B(p_A,2\e) $ and an outer region
$B_{out}=B(p_A,2\e)-B_{in}$.  On $B_{out}$, 
we have $r_p\geq d/2$, and hence $r_A/r_p\leq
(r_p+d)/r_p \leq 3$.  Thus
\be r_p^{-\d}r_A^n |F|=(r_A/r_p)^\d r_A^{n-\d}|F|\lessim 
r_A^{n-\d}|F|,
\ee
implying
\bearray\non
\|\b r_p^{-\d}r_A^nF\|_{L^2(B_{out})}\lesssim 
\|\b r_A^{n-\d}F\|_{L^2(B_{out})}&\lesssim &
\|\b r_A^{n-\d}F\|_{L^2(B(p_A,2\e))} 
\\
&\lesssim &
\left\{ \begin{array}{rr}
\l^{n-\d}, & -2<n-\d<2,\\ \l^2\sqll, &  n-\d =2,\\
\l^{1+(n-\d)/2}, & n-\d>2
\end{array}\right. .\non \\
\label{cor10l_simp}
\eearray

For the integral over $B_{in}$, first suppose $n-\d\leq 2$ and 
separately consider the cases
$d\leq \l$, $d\geq \l$. In both cases note that $d/2\leq r_A\leq 3d/2$
in this region.  When $d\leq \l$, we then have $r_A\lessim\l$ and
$|F|\lessim\l^{-2}$ on $B_{in}$, so
\be
\|\b r_p^{-\d}r_A^nF\|_{L^2(B_{in})}\lessim \l^{n-2}
\| r_p^{-\d} \|_{L^2(B_{in})}\lessim \l^{n-2}d^{2-\d}\lessim \l^{n-\d}.
\ee
On the other hand if $d\geq\l$, then since $r_A/d$ is bounded above
and below on $B_{in}$, (\ref{ptwsfest2}) implies $r_A^n|F|\lessim \l^2
r_A^{n-4}\lessim \l^2 d^{n-4}$.  Hence 
\be
\|\b r_p^{-\d}r_A^nF\|_{L^2(B_{in})}\lessim \l^2d^{n-4}
\| r_p^{-\d} \|_{L^2(B_{in})}\lessim \l^2d^{n-\d-2}\leq\l^{n-\d}
\ee
since $n-\d\leq 2$.  Combining this with the estimate for $B_{out}$ we
obtain the top two lines of (\ref{cor10f_simp}).

If $n-\d>2$, separately consider the cases $d\leq 4\e$ and $d\geq
4\e$. If $d\leq 4\e$, the procedure for the case $\d\geq \l$ above
yields
\be
\|\b r_p^{-\d}r_A^nF\|_{L^2(B_{in})}\lessim \l^{n-2}
\| r_p^{-\d} \|_{L^2(B_{in})}\lessim \l^2d^{n-\d-2}
\lessim\l^{1+(n-\d)/2},
\ee
the same bound as on $B_{in}$. 
If $d\geq 4\e$ then 
on the support of $\b$ we have $r_p\geq \e$, so
$r_p^{-\d}r_A^n|F|\lessim
\l^{-\d/2}r_A^n|F|$.  Thus (\ref{cor10a_simp}) yields the remaining
case of 
(\ref{cor10f_simp}) for the bound on $\|r_p^{-\d}r_A^nF\|_2$.

   The method for bounding $\|\b r_p^{-\d}r_A^{n+1}\na F\|$, is
essentially identical to the method for bounding $\|\b
r_p^{-\d}r_A^{n+1}\na F\|_2$, except that for the estimates over $B_{in}$,
first multiply by a cutoff function of the form
$\b_{\rm std}(2r_p/d)$, and then integrate by parts as in
(\ref{cor10i_simp}).  

To bound $\lnorm \chi r_p^{-\d}r_A^nF\rnorm_2$, note that on
$\supp(\chi)$ we have $|F|\leq {\rm const}$ and $r_A\leq \e$, so
\be
\lnorm \chi r_p^{-\d}r_A^nF\rnorm \lesssim \l^{n/2}\lnorm
\chi r_p^{-\d}\rnorm_2
\lesssim \l^{n/2}\lnorm
r_p^{-\d}\rnorm_{L^2(B(p,2\e))}
\lesssim \l^{1+(n-\d)/2}.
\ee
Similarly
$\lnorm \chi r_p^{-\d}r_A^n\na F\rnorm \lesssim \l^{n/2}
\lnorm \chi  r_p^{-\d}\na F\rnorm_2$, and the same procedure as for
$\d=0$ completes the work. 

(b) First write
\be
\lnorm r_p^{1-\d} r_A^n F\rnorm_4
\leq
\lnorm (1-\tilde{\b})r_p^{1-\d} r_A^n F\rnorm_4+
\lnorm \tilde{\b} r_p^{1-\d} r_A^n F\rnorm_4,
\ee
where $\tilde{\b}$ is a cutoff of scale $\e_0$ centered at $p_A$.
On the support of $1-\tilde{\b}$ we have $|F|\leq {\rm const}$, so the
first term on the RHS is bounded by a constant.  The second term can
be estimated as in the proof of (b).
\qed

We are now in a position to bound $|G^A_0Rem'_{2,{\rm loc}}|$
pointwise, but we postpone this until we have collected the estimates
needed to bound the semi-local and non-local contributions to
$G^A_0Rem_2'$.  These require bounds on norms of
$\xi_X=\D^*\gplus\D\tx$, which in turn require pointwise bounds on
$\D\tx$.

\begin{lemma}\label{doglemma1} 
For any vector field $X$ on $N$, and any ASD connection $A$, we
have the pointwise formulas
\be\label{doglemma1a}
(d^A)^*\i_XF_A\ =\  \langle d_+X^*,F_A\rangle
\ee and
\be\label{doglemma1b}
d^A_+(\i_X F_A)=\mbox{\rm Sym}^2_0(\na X^*)\sharp F_A,
\ee
\ssn 
where $X^*$ is the metric dual of $X$, $d_+X^*$ is the self-dual part
of $dX^*$, $\mbox{\rm Sym}^2_0(T)$ denotes the
traceless symmetric part of a rank-two tensor field $T\in
\Gamma(T^*N\tensor T^*N)$, and, in a local
orthonormal basis ${\th^i}$ of the cotangent bundle, $\lb T, F \rb =
\frac{1}{2} T_{ij}F_{ij}\in $ and 
$T\sharp F=T_{ij} F_{jk} \th^i\wedge\th^k$.  Hence
\be\label{doglemma1c}
|\D^A(\i_XF_A)|\lessim (|d_+X^*|+|\mbox{\rm Sym}^2_0(\na X^*)|)|F_A|
\ee
and 
\bearray\non
|\Na\D^A(\i_XF_A)|& \leq & c\left(
(|\na(d_+X^*)|+|\na(\mbox{\rm Sym}^2_0(\na X^*))|)|F_A| \right. \\
&& \left. +(|d_+X^*|+|\mbox{\rm Sym}^2_0(\na X^*)|)|F_A|\right).
\label{doglemma1d}
\eearray

Hence if $X\in\mathfrak{h}_A$, then
\be\label{doglemma2a}
|\D^A\tX| \lesssim
(r_A\b +\e^{-1}\chi)|\hat{X}|\ |F_A|
\ee
and
\be\label{doglemma2a2}
\|(\D^A)^*\D^A\tx|\leq c |\Na\D^A\tX| \lesssim
(\b+\e^{-2}\chi)|\hat{X}|\ (|F_A|+r_{A}|\Na F_A|).
\ee

\end{lemma}

\pf 
 Using the facts that that $d^*=-*d*$, $d^AF=0$, and $*F=-F$, we
have 
\be
(d^A)^*(\i_XF)=-*d^A(*\i_X(*F))= *d^A(X^*\wedge F) =*(dX^*\wedge F) =
(dX^*,F);
\ee
this gives (\ref{doglemma1a}).  Now fix $p\in N$. Calculating in a
local orthonormal frame $\{e_i\}$ of $TN$ and dual coframe $\{\th^i\}$
with $\na e_i|_p=0$,
\bearray\nonumber
d^A_+(i_XF) &=& p_+(\sum \th^i\wedge\i_{\nabla_i X}F)\\ \nonumber &=&
\sum (\na_iX_j) p_+(\th^i\wedge\i_{e_j}F)\\ \label{rep1} &=&
p_+\left(\mbox{\rm Sym}^2(\na X)\sharp F\right)
\eearray
by Lemma 2.3 of [G1].  Since for any symmetric 2-tensor $T$, the
pure-trace part of $T$ yields a {\em self}-dual 2-form under the
operation $\sharp F$, we may replace $\mbox{\rm Sym}^2$ by $\mbox{\rm
Sym}^2_0$ in (\ref{rep1}), and by simple representation theory, the
$p_+$ in (\ref{rep1}) is redundant.

(\ref{doglemma2a}--\ref{doglemma2a2}) follow from 
Lemma \ref{doglemma1} and a pointwise computation of $d_+X^*$, ${\rm
Sym}^2_0(\na X^*)$ that we leave to the reader.
\qed

\begin{cor}\label{doglemma2} For all $X\in\mathfrak{h}_A$, and all
$p\in N$, the 
elements $\tX\in {\cal H}_A$ satisfy the following integral bounds.

(a) If $-1<m<2$, then 
\be
\lnorm r_A^m\tx\rnorm_4\leq c(m)\l^{m-1}.
\ee

(b) If $-1<m\leq 0$, or if $d={\rm dist}(p,p_A)\lessim\e$ and $-1<m<2$, then
\be
\label{doglemma2b}\label{doglemma2b2}\label{doglemma2c}\label{doglemma2c2}
\label{doglemma2c3}
\lnorm r_p^m \D^A\tx\rnorm_2 \leq c(m)
\l^{m/2}(b\cdot \l^{1/2}+a).
\ee

(c) For all $\d\in (0,1)$, 
\be
\label{doglemma2d}
\lnorm r_{p}^{-\d} (\D^A)^*\D^A\tX \rnorm_2 
\leq c\lnorm r_{p}^{-\d} \Na \D^A\tX \rnorm_2
\lesssim \l^{-\d}(b+a\cdot \l^{-1/2}) .
\ee
\end{cor}

\ssn\pf Using (\ref{doglemma2a}--\ref{doglemma2a2}), plus $|{X}|
\leq b+a\l^{-1}r_A$, most of these bounds follow directly from Lemma
\ref{cor10_simp}.  The exception is (\ref{doglemma2b}) in the case
$m>0$, for which one must also use the triangle inequality $r_p\leq
r_A+d\lessim r_A+\l^{m/2}$.   \qed

We are now in a position to derive our final estimates on the norms of
$\xi$ needed to bound $|G^A_0Rem'_{2,{\rm semiloc}}|$ and
$|G^A_0Rem'_{2,{\rm nonloc}}|$ pointwise.  We also use the opportunity
to prove (\ref{lemmad1b}).

\begin{prop}\label{dogprop2} \label{pigcor11}
There exists $\d_0>0$ such that if 
$0<\d<\d_0$ and $0\leq \d'<\d_0$, then for $\xi=(\D^A)^*\gplus
\D^A\tx$ (with $X\in \mathfrak{h}_A$) we have the following.

(a) If $0<\d<\d_0$ then
\bearray
|\xi(p)| &\lessim\l^{\d'}\left(
b\cdot \l^{-1}+a\cdot \l^{-3/2}\right).
\label{dogprop2a}
\eearray
If furthermore  $d=r_A(p)\geq 4\e=c\l^{1/2}$, then
\be
|\xi(p)| \leq c(\d) r_A(p)^{-2-\d-\d'}\l^{\d'}
(b\cdot\l+a\cdot\l^{1/2}).
\label{dogprop2c}
\ee

\ss(b)
\be\label{lemmad1c}
\|\xi\|_2 \lessim \l^{\d'/2}(b\cdot \l + a\cdot \l^{1/2})
\ee
Since $\xi_X=\tx-\pi_A\tx$, this implies (\ref{lemmad1b}). 

\ss(c) If $0\leq \d\leq \d_0$ then
\bearray
\lnorm r_p^{-\d}\xi\rnorm_4
&\lesssim&
\l^{\d'}(b+ a\cdot \l^{-1/2}).
\label{pigcor11a}
\label{pigcor11b}
\eearray

\ss(d) 
If $|\d|<\d_0$ then
\bearray\lnorm r_A^{1+\d}\xi\rnorm_4
&\lessim& \l^{\d/2}(b\cdot \l + a\cdot \l^{1/2}).
\label{pigcor11c}
\label{pigcor11d}
\label{pigcor11d2}
\eearray

\end{prop}

\pf (a) We will omit writing the $\d$-dependence of the constants.
 From (\ref{lemma3d}) given $\d,\d_0$ as above there exists $\d'>0$
such that
\bearray
|\xi(p)|&\lessim& 
\lnorm r_p^{-\d} (d^A_+)^*\w \rnorm_2 + 
\l^{2\d'-1}\lnorm r_A^{-\d-2\d'}\w\rnorm_2+
\l^{2\d'+\d-2}\lnorm r_A^{1-2\d-2\d'}\w\rnorm_2,
\label{estcor8a_copy}
\label{dogprop2e}
\eearray
where $\w={\cal D}^A\tx$.
Using Lemma \ref{doglemma2}, we compute
\be
\l^{2\d'-1}\lnorm r_A^{-\d-2\d'}\w\rnorm_2+
\l^{2\d'+\d-2}\lnorm r_A^{1-2\d-2\d'}\w\rnorm_2
\lessim
\l^{\d'}(b\cdot\l^{-1}+ a\cdot\l^{-3/2}).
\label{dogprop2f}
\ee  
The bound on $\lnorm r_p^{-\d} (d^A_+)^*\w \rnorm_2$ from Lemma
\ref{doglemma2} is smaller than this, so we obtain (\ref{dogprop2a}).

   For (\ref{dogprop2c}), apply (\ref{estlemma5d}) and Lemma
\ref{doglemma2}.

(b), (c), and (d). Apply Proposition \ref{lemma3} and Corollary
\ref{doglemma2}. 
\qed

   We remark that by using the pointwise decay estimate
(\ref{dogprop2e}) one can obtain the weighted $L^4$ decay
\be\label{l4decay}
\|r_p^{-\d}\|_4\lessim\l^{-\d}r_A(p)^{-\d}(b + a\cdot \l^{-1/2})
\ee
for $d\geq c\l^{1/2}$, but this is of no help to us.

   We're now ready to collate all the estimates needed to prove Theorem
\ref{thmrem2}. 

\begin{cor}\label{finalests}
(a) There exist $\d_0>0,\d'>0$ such that for
$0 \leq \d<\d_0$, the following are true.

\bearray
\lnorm  r_p^{-\d}Rem'_{2,{\rm loc}}(X,Y) \rnorm_2 
&\lesssim& \l^{-\d/2}( b^2+ba\cdot\l^{-1/2} +a^2\cdot\l^{-1}).
\label{est1a_simp}\label{est1}\label{est1b_simp}\\
\lnorm r_A^{1\pm\d} Rem'_{2,{\rm loc}}(X,Y) \rnorm_2 &\lesssim&
\l^{1/2\pm\d/2}(b^2+ ba\cdot\l^{-1/2} +a^2\cdot\l^{-1}).\\
\label{est1c1_simp} 
\lnorm r_p^{-\d} Rem'_{2,{\rm semiloc}} \rnorm_2 
&\lessim&
\l^{-1+\d'}(b^2+ba\cdot\l^{-1/2} +a^2\cdot\l^{-1/2}).
\label{est2}\label{est2b_simp} \label{est2a_simp} \\
\lnorm r_A^{1\pm\d} Rem'_{2,{\rm semiloc}} \rnorm_2 &\lesssim& 
\l^{\d'}(b^2 +ba\cdot \l^{-1/2} + a^2\cdot \l^{-1/2}).
\label{est2c_simp} \\
\lnorm r_p^{-\d}Rem'_{2,{\rm nonloc}}(X,Y) \rnorm_2 
&\lesssim&  \l^{\d'}(b^2+ba\cdot\l^{-1/2} +a^2\cdot\l^{-1})
\label{est3b_simp}
\label{est3a_simp}\label{est3}.\\
\lnorm r_A^{1\pm\d}Rem'_{2,{\rm nonloc}}(X,Y) \rnorm_2 &\lesssim& 
\l^{1+\d'}(b^2+ba\cdot\l^{-1/2} +a^2\cdot\l^{-1})
\label{est3c_simp}
\eearray

(b) There exists $\d'>0$ such that for all $p\in N$, the following are
true. 
\bearray
|G^A_0 Rem'_{2,{\rm loc}}(X,Y)|(p) &\lessim &
\l^{-1/2+\d'}(b^2+ ba\cdot\l^{-1/2} +a^2\cdot\l^{-1}).
\label{est1e_simp}\\
|G^A_0 Rem'_{2,{\rm semiloc}}(X,Y)|(p) &\lessim & 
\l^{-1+\d'}(b^2+ba\cdot\l^{-1/2} +a^2\cdot\l^{-1/2}).
\label{est2d_simp}\\
|G^A_0 Rem'_{2,{\rm nonloc}}(X,Y)|(p) &\lessim & 
\l^{\d'}(b^2+ba\cdot\l^{-1/2} +a^2\cdot\l^{-1}).
\label{est3d_simp}
\eearray

(c) There exists $\d_0>0,\d'>0$ such that if $0<\d<\d_0$
and $d\geq 4\e=c\l^{1/2}$ the following are true.
\bearray 
|G^A_0 Rem'_{2,{\rm loc}}(X,Y)|(p) &\lessim& d^{-1-\d-\d'}\l^{\d'}
(b^2\cdot \l^{1/2}+ ba +a^2\cdot\l^{-1/2})
\label{est1e2_simp}\\
|G^A_0 Rem'_{2,{\rm semiloc}}(X,Y)|(p) &\lessim& 
d^{-1-\d-\d'}\l^{\d'}(b^2+ba\cdot\l^{-1/2} +a^2\cdot\l^{-1/2}).
\label{est2e_simp}\\
|G^A_0 Rem'_{2,{\rm nonloc}}(X,Y)|(p) &\lessim&
d^{-1-\d-\d'}\l^{\d'}
\left( b^2\cdot \l +ba \cdot \l^{1/2}+ a^2 \right).
\label{est3e_simp}
\eearray

\end{cor}

\ss\pf 
(a)
These bounds follow directly from Lemma \ref{oldlemma39}, Lemma
Corollary \ref{doglemma2} the $L^4$ bounds in Proposition
\ref{dogprop2}, and H\"{o}lder's inequality.  

\ss
(b) Use part (a) and Lemma \ref{estlemma5}.

\ss
(c) Since $Rem'_{2,{\rm loc}}(X,Y)$ and $Rem'_{2,{\rm semiloc}}(X,Y)$
are supported in $B(p_A,2\e)$, for these terms we can apply
(\ref{estlemma5c}) and the corresponding bounds in (a).  As
$Rem'_{2,{\rm semiloc}}(X,Y)$ is not locally supported, we appeal
instead to (\ref{estlemma6d}):
\be\label{estlemma6c_copy}
|G^A_0 \{\xi,\xi\}|(p) \lessim
\left( d^{-1-\d-\d'}
\lnorm r_A^{1+\d'}\{\xi,\xi\}\rnorm_2 + 
\lnorm r_p^{-\d}\tilde{\b}\{\xi,\xi\}
\rnorm_2\right),
\ee
where $\tilde{\b}$ is a cutoff of scale $d/2$ as in Lemma
\ref{estlemma5}. 

If we estimate $\lnorm r_A^{1+\d'}\{\xi,\xi\}\rnorm_2$ using
(\ref{est3c_simp}) we obtain the right-hand side of (\ref{est3e_simp}).
Were we next to estimate $\lnorm r_p^{-\d}\tilde{\b}\{\xi,\xi\}\rnorm_2$,
the resulting bound would be to large to be of use.
Instead, since $d\leq r_A\leq 3d$ on the support of
$\tilde{\b}$, we can use the pointwise bound (\ref{dogprop2c}) to find
\bearray
\lnorm r_p^{-\d}\tilde{\b}\{\xi,\xi\}\rnorm_2 &\lessim& \non
\lnorm r_p^{-\d}\tilde{\b}\rnorm_2\left(
d^{-2-\d-\d'}\l{\d'}
\left(b\cdot\l+a\cdot\l^{1/2}
\right) \right)^2 \\
&\lessim& \non
d^{-2-3\d-2\d'}\l^{2\d'}(b^2\cdot \l^2+ba\cdot \l^{3/2}+a^2\cdot\l)\\
&\lesssim&
d^{-1-\d-\d'}\l^{3/2-\d''}(b^2\cdot +ba\cdot \l^{-1/2}+a^2\cdot\l^{-1}),
\eearray
which is much smaller than our bound on $\lnorm
r_A^{1+\d'}\{\xi,\xi\}\rnorm_2$. Thus (\ref{est3e_simp}) follows.
\qed

Finally, we have the

\ms\ni{\bf Proof of Proposition \ref{thmrem2}.}
(a) Add the bounds (\ref{est1e_simp}--\ref{est3d_simp}) to obtain
(\ref{thmrem2a}). If we
add the bounds 
(\ref{est1e2_simp}--\ref{est3e_simp}) we obtain a stronger bound than
(\ref{thmrem2b}): 
\be
| Rem_2(X,Y)|
\lesssim
r_A^{-1-\d'}\l^{\d'}\left( b^2 +ba\cdot
\l^{-1/2}+a^2\cdot\l^{-1/2} \right).
\ee

(b) In the proof of part (a), the only way in which $\xi$ entered was
through the $L^4$ bounds on $\|r_p^{-\d}\xi_4\|_4$,
$\|r_A^{1+\d'}\xi_4\|_4$, and the pointwise decay (\ref{dogprop2c}).
Hence our assertion follows from the hypothesis ({\bf Z5}) of section
\ref{sectfiber}.
\qed


\section{Appendix}
\setcounter{equation}{0}

The point of the following weighted Sobolev inequality 
is that on an $m$-dimensional manifold there is no Sobolev
embedding $L_1^m\embed L^\infty$, but the failure is borderline.  Thus
by introducing an arbitrarily small weight into the Sobolev norm we
are able to obtain an embedding.

\begin{lemma}\label{estlemma4}
Let $E\to N$ be a Riemannian vector bundle with metric-compatible
connection $\na$, where $N$ is compact, Riemannian, and $m$-dimensional
($m>1$).  Given $p\in N$ and $R_2>R_1>0$, let $\Omega(p;R_1,R_2)$
denote the annulus $\{R_1\leq r_p \leq R_2 \}$, where $r_p$ is
distance to $p$.  There exists a constant $c$ such that for any $\d>0,
R_2>R_1>0$ (but smaller than the injectivity radius), any $p\in N$,
and any $\phi\in \Gamma(E)$ we have
\be
|\phi(p)|
\leq 
c\d^{-(1-1/m)}  R_2^\d \left(
\frac{1}{R_2-R_1}\lnorm r_p^{-\d} \phi \rnorm_{L^m(\Omega(p;R_1,R_2))} +
\lnorm r_p^{-\d} \na \phi \rnorm_{L^m(B_{R_2}(p))}\right).
\label{estlemma4a1}
\ee
Consequently,

\be
|\phi(p)| \leq 
\d^{-(1-1/m)} \left(
\lnorm \phi \rnorm_{L^m(N)}+
\lnorm r_p^{-\d} \na \phi \rnorm_{L^m(N)}\right) 
\label{estlemma4a1.5}
\ee
and 
\be
\lnorm \phi \rnorm_{L^\infty(N)} 
\leq 
c \d^{-(1-1/m)} \left(
\lnorm \phi \rnorm_{L^m(N)}+
\sup_{p\in N}\left(\lnorm r_p^{-\d} \na \phi \rnorm_{L^m(N)}\right) \right).
\label{estlemma4a2}
\ee
\end{lemma}

\ss\pf By Kato's inequality, it suffices to prove this for the trivial
real line bundle, i.e. for functions on $N$.

   First replace $N$ by $\bfr^m$ and consider a compactly supported
function $f\in C_0^\infty(B(0,R))$.  Let $\theta\in S^{m-1}$.  Then,
using polar coordinates on $B(0,R)$, we have
\be
|f(0)| = \left|\int_0^R \frac{\partial f}{\partial r}(r,\theta) dr \right|
\leq \int_0^R |\na f|(r,\theta) dr,
\ee
implying
\be
\vol(S^{m-1})|f(0)|\leq \int_{S^{m-1}} d\theta \left(
\int_0^R |\na f|(r,\theta) dr\right) = \int_{B(0,R)}|\na f|\cdot
r^{1-m} dvol.
\ee
Applying the same argument on a normal-coordinate ball $B(p,R)$ in $N$
(where $f\in C_0^\infty(B(p,R))$), using the compactness of $N$ to get
uniformity in the constants below, we obtain
\bearray
|f(p)| \leq \non c \int_{B(p,R)} |\na f| r_p^{1-m} dvol
&=& \non c \int_{B(p,R)} r^{-\d}|\na f| r_p^{1-m+\d} dvol\\
&\leq& \non c\lnorm r_p^{-\d} \na f \rnorm_{L^m(B(p,R))} 
\lnorm r_p^{1-m+\d} \rnorm_{L^{m/(m-1)}(B(p,R))} \\
&\leq& c \d^{-(1-1/m)} R^\d \lnorm r_p^{-\d} \na f
\rnorm_{L^m(B(p,R))}. \label{estlemma4b}
\eearray

Now remove the assumption that $f$ is supported inside a normal
coordinate ball.  Replace $f$ in the preceding argument by $\b(r) f$,
where $\b$ is a cutoff function identically 1 for $r\leq R_1$ and
vanishing for $r\geq R_2$; thus $|\na\b|\leq c/(R_2-R_1)$. We then
have
\bearray
|f(p)| &\leq& \non
c\d^{-(1-1/m)} R_2^\d \lnorm r_p^{-\d} \na (\b f)
\rnorm_{L^m}\\
&\leq& \non
c\d^{-(1-1/m)} R_2^\d \left(
\lnorm r_p^{-\d} (\na \b)f \rnorm_{L^m} +
\lnorm \b r_p^{-\d} \na f \rnorm_{L^m}\right)\\
&\leq& \non
c\d^{-(1-1/m)} R_2^\d \left(
\frac{1}{R_2-R_1}\lnorm r_p^{-\d} f \rnorm_{L^m(\Omega(R_1,R_2,p))} +
\lnorm r_p^{-\d} \na f \rnorm_{L^m(B_{R_2}(p))}\right),\\
\label{estlemma4c}
\eearray
yielding (\ref{estlemma4a1}).  Taking $R_2=2R_1$ to be, say, half the
injectivity radius of $N$, we obtain (\ref{estlemma4a2}).
\qed

As a corollary, we have

\begin{cor}\label{mainsobolev}
 Let $E,N,\na$ be as in Lemma \ref{estlemma4}, and assume dim($N$)=4.
Then for all $\d\in (0,1)$, there exist constants $c(\d)$ such that for all
$\phi\in\Gamma(E)$,
\be\label{estlemma4a10}
|\phi(p)|
\leq 
c(\d)\left(
\lnorm  \phi \rnorm_{L^2(N)}+
\lnorm r_p^{-\d} \na \na \phi \rnorm_{L^2(N)}\right),
\ee
and hence 
\be
\lnorm \phi \rnorm_{\infty} 
\leq 
c(\d)\sup_{p\in N}\left(
\lnorm  \phi \rnorm_{2}+
\lnorm r_p^{-\d} \na \na \phi \rnorm_{2}\right).
\label{estlemma4a9}
\ee 
\end{cor}

\ni\pf Applying Lemma \ref{estlemma4} with $m=4$, we have
\be
|\phi(p)| \leq c(\d) \left(
\lnorm  \phi \rnorm_{4}+ 
\lnorm r_p^{-\d} \na \phi \rnorm_{4}\right)
\label{estlemma4a3.5}
\ee
Using the Sobolev embedding $L_1^2(N)\embed L^4(N)$ we then find
\bearray
|\phi(p)| 
&\leq &
c(\d) \left(
\lnorm  \phi \rnorm_{2}+
\lnorm r_p^{-\d-1} \na \phi \rnorm_{2}+
\lnorm r_p^{-\d} \na \na \phi \rnorm_{2}\right)
\label{estlemma4a6}
\eearray
But since $\d<1$, we also have the weighted Sobolev inequality of
Lemma 3.1 of \cite{Gppt}:
\be
\|r^{-1-\d}\psi \|_2 \leq c( \| r^{-\d} \psi \|_2+ \| r^{-\d} \na \psi
\|_2)
\ee 
(the proof is again a polar-coordinate computation).  Using this we
can bootstrap (\ref{estlemma4a6}) into the form (\ref{estlemma4a10}) .
\qed





\end{document}